\documentclass[11pt]{amsart}

\usepackage{amssymb}

\makeindex

\newcommand{\POD}[1]{{\cD}_{#1}}

\def\proof{\noindent {\bf Proof}}
\def\qed{$\blacksquare$}

\def\beq{\begin{equation} }
\def\eeq{\end{equation} }
\def\eps{\varepsilon}

\def\benum{\begin{enumerate} }
\def\eenum{\end{enumerate} }

\def\pk{p^{(k)}}

\def\RR{ {\mathbb{R}} }

\def\Rnn{ \RR^{n \times n} }

\def\RRx{\RR \langle x \rangle}

\def\RRxxS{\RR\langle x,x^\ast \rangle}

\def\FF{\RRxxS}

\def\bbD{\mathbb D}

\def\cA{ {\mathcal A} }
\def\cB{ {\mathcal B} }
\def\cC{ {\mathcal C} }
\def\cD{ {\mathcal D} }
\def\cE{ {\mathcal K} }

\def\cI{ {\mathcal I} }

\def\cL{ {\mathcal L} }
\def\cM{{\mathcal M}}
\def\cN{ {\mathcal N} }

\def\cP{{\mathcal P}}

\def\cS{{\mathcal S} }

\def\cV{{\mathcal V}}
\def\cW{{\mathcal W}}
\def\cU{\mathcal U}
\def\cX{\mathcal X}
\def\cY{\mathcal Y}

\def\tX{\tilde X}

\def\j1tog{ j= 1, \ldots, g }

\def\bbR{\mathbb R}

\def\cI{{\mathcal I}}

\def\L0t{\ (L_0 \otimes I_n ) \ }

\def\PLINE1{\beta}

\newcommand{\CC}{{\mathbb C}}

\newtheorem{thm}{Theorem}[section]
\newtheorem{cor}[thm]{Corollary}
\newtheorem{lem}[thm]{Lemma}
\newtheorem{prop}[thm]{Proposition}
\theoremstyle{definition}

\theoremstyle{remark}

\numberwithin{equation}{section}

\newtheorem{exa}[thm]{Example}
\newenvironment{example}%
         {\begin{exa}}
             {{\hfill $\Box ~~$}\end{exa}}

\newcounter{Inc}

\begin{document}
\setcounter{page}{1}
\title[Positive polynomials and the spectral theorem]
{Positive polynomials  in scalar and matrix variables, the spectral theorem and optimization}
\dedicatory{Tibi Constantinescu, in memoriam}

\author{ J. William Helton}

\address{Department of Mathematics,
   University of
California at San Diego, La Jolla CA 92093}
\email{helton@math.ucsd.edu}

\author{Mihai Putinar}
\thanks{Partially supported by grants from the National Science Foundation
and the Ford Motor Co.}

\address{Department of Mathematics, University of California, Santa
Barbara, CA 93106}

\email{mputinar@math.ucsb.edu}

\today

\begin{abstract}
We follow a stream of the history of positive matrices and positive
functionals, as
applied to
algebraic sums of squares
decompositions, with emphasis on the interaction between classical
moment problems, function theory of one or several complex
variables and modern operator theory.
The second part of the
survey focuses on recently discovered connections between real
algebraic geometry and optimization
as well as   polynomials in matrix variables and
some control theory problems.
These new applications have
prompted a series of recent studies devoted to the structure of
positivity and convexity
in a free $*$-algebra, the appropriate setting
for analyzing inequalities on polynomials having matrix variables.
We sketch  some of these
developments, add to them and comment on the rapidly growing literature.
\end{abstract}

\maketitle

\section{Introduction}

This is an essay, addressed to non-experts, on the structure of
positive polynomials on semi-algebraic sets, various facets of the
spectral theorem for Hilbert space operators,
inequalities and sharp constraints
for elements of a free $\ast-$algebra,
and some
recent applications of all of these to polynomial optimization and
engineering.
The circle of ideas exposed below is becoming
increasingly popular but not known in detail outside the
traditional groups of workers in functional analysis or real
algebra who have developed parts of it.
For instance, it is not
yet clear how to teach and facilitate the access of beginners to
this beautiful emerging field.
The exposition of topics  below may
provide  elementary ingredients for such a course.

The unifying concept behind all the apparently diverging topics
mentioned above is the fact that universal positive functions (in
appropriate rings) are sums of squares. Indeed, when we prove
inequalities we essentially complete squares, and on the other
hand when we do spectral analysis we decompose a symmetric or a
hermitian form into a weighted (possibly continuous) sum or
difference of squares. There are of course technical difficulties
on each side, but they do not obscure the common root of algebraic
versus analytical positivity.

We will encounter quite a few positivity criteria, expressed in
terms of: matrices, kernels, forms, values of functions,
parameters of continued fractions, asymptotic expansions and
algebraic certificates. Dual to sums of squares and the main
positive objects we study are the power moments of positive
measures, rapidly decaying at infinity. These moments will be
regarded as discrete data given by fixed coordinate frames in the
correspondence between an algebra (of polynomials or operators)
and its spectrum, with restrictions on its location. Both concepts
of real spectrum (in algebraic geometry) and joint spectrum (in
operator theory) are naturally connected in this way to moment
problems. From the practitioner's point of view, moments represent
observable/computable numerical manifestations of more complicated
entities.

It is not a coincidence that the genius of Hilbert presides over
all aspects of positivity we will touch. We owe him the origins
and basic concepts related to: the spectral theorem, real algebra,
algebraic geometry and mathematical logic. As ubiquitous as it is,
a Hilbert space will show up unexpectedly and necessarily in the
proofs of certain purely algebraic statements. On the other hand
our limited survey does not aim at offering a comprehensive
picture of Hilbert's much wider legacy.

Not unexpected, or, better later than never,
the real algebraist's
positivity and the classical  analyst's positive
definiteness have recently merged into a powerful framework;
this is needed
and shaped by several applied fields of mathematics.
We will
bring into our discussion one principal customer: control theory.
The dominant development in linear systems engineering in the
1990's was matrix inequalities and many tricks and ad hoc
techniques for making complicated matrix expressions into tame
ones, indeed into the
 Linear Matrix Inequalities, LMIs, loved by all
 who can obtain them.
Since matrices do not commute a large portion of the subject could
be viewed as manipulation of polynomials and rational functions of
non-commuting (free) variables, and so a beginning toward  helpful
mathematical theory would be a  semi-algebraic geometry for free
$*$-algebras, especially its implications for convexity.
Such
ventures sprung to life within the last five years and this
article attempts to introduce, survey and fill in some gaps in
 this rapidly expanding
area of noncommutative semi-algebraic geometry.

The table of contents offers an idea of the topics we touch in the
survey and what we left outside. We are well aware that in a
limited space while  viewing  a wide angle, as captives of our
background and preferences, we have omitted key aspects. We
apologize in advance for all our omissions in this territory, and
for  inaccuracies when stepping on outer domains; they are all
non-intentional and reflect our limitations. Fortunately, the
reader will have the choice of expanding and complementing our
article with several recent excellent surveys and monographs
(mentioned throughout the text and some recapitulated in the last
section).

The authors thank the American Institute of Mathematics, Palo
Alto, CA, for the unique opportunity (during a 2005 workshop) to
interact with several key contributors to the recent theory of
positive polynomials. They also thank the organizers of the ``Real
Algebra Fest, 2005", University of the Saskatchewan, Canada, for
their interest and enthusiasm. The second author thanks the Real
Algebra Group at the University of Konstanz, Germany, for offering
him the possibility to expose and discuss the first sections of
the material presented below.

We dedicate these pages to Tibi Constantinescu, old time friend
and colleague, master of all aspects of matrix positivity.\footnote {
 {\bf Advice to the reader.} {\it
Although the ordering of the material below follows a logic
derived from a general theoretical and historical perspective,
most of the sections, and sometimes even subsections, can be
 read independently. We have tried to keep to a minimum the number of
cross references and we have repeated definitions.

For instance, the reader oriented towards optimization and
engineering applications of the decompositions of polynomials into
sums of squares can start with sections on the Positivstellensatz
\S \ref{sec:PosSS}, \S \ref{sec:posSScompact},
then turn to sections on global optimization and engineering
\S \ref{sec:applSAG}
and sections \S  \ref{sec:convexity},
\S \ref{sec:engLMI}  on systems whose  structure
does not depend on the size of the system (is dimension free).
After this we suggest a tour
 for the bigger picture through  neighboring sections.

The  operator theorist should find most of the paper straightforwrd
 to read
with the exception of
 Section \ref{sec:logic} which relates this very functional
analytic topic to logic. Also an experienced operator theorist
could start reading in many places, for example, in \S \ref{sec:nonCom}
or in  \S \ref{sec:moments} or in  \S \ref{sec:applSAG} or
even at the beginning.
We reassure everyone that the important
 general Positivstellensatz in \S \ref{sec:PosSS}
whose proof  requires logic,
is stated in a self-contained way.

Intriguing for the algebraist and logician might be Sections
\ref{sec:specThm}, \ref{sec:moments}, and \ref{sec:complexVar}
which contain analytic material (mostly derived from the spectral
theorem and its many facets) which has cousins and even closer
relatives in algebra.

In any case, all readers should be aware of the modular structure
of the text, and try non-sequential orderings to access it. }}

{\footnotesize \tableofcontents}

\section{The spectral theorem}
\label{sec:specThm}

The modern proof of the spectral
theorem for self-adjoint or unitary operators uses commutative
Banach algebra techniques, cf. for instance \cite{Douglas}. This
perspective departs from the older, and more constructive approach
imposed by the original study of special classes of integral
operators. In this direction, we reproduce below an early idea of
F. Riesz \cite{Riesz} for defining the spectral scale of a
self-adjoint operator from a minimal set of simple observations,
one of them being the structure of positive polynomials on a real
interval.

\subsection{Self-adjoint operators} Let $H$ be a separable,
complex Hilbert space and let $A \in \mathcal L(H)$ be a linear,
continuous operator acting on $H$. We call $A$ {\it self-adjoint}
\index{self-adjoint}
if $A=A^\ast$, that is $\langle Ax,x\rangle
\in \mathbb R$ for all vectors $x \in H$. The continuity
assumption implies the existence of bounds \begin{equation}
\label{bounds}
  m \| x\|^2 \leq \langle Ax,x\rangle
  \leq M \|x\|^2, \ \ x \in H.
\end{equation}
  The operator $A$ is called
{\it non-negative},
\index{non-negative operator}
denoted in short
$A \geq 0$, if
$$\langle Ax,x\rangle \geq 0, \ \  x \in H.$$
The operator $A$ is {\it positive}
\index{positive operator}
if it
is non-negative and $(\langle Ax,x\rangle = 0) \ \Rightarrow \
(x=0).$

We need a couple of basic observations, see \S 104 of \cite{RN}.
The real algebraists should enjoy comparing
 these facts with the axioms of an
order in an arbitrary ring.\\

{\bf a).} {\it A bounded monotonic sequence of self-adjoint
operators converges (in the strong operator topology) to a
self-adjoint
operator.}

\smallskip

\noindent
Indeed, assume $0 \leq A_1 \leq A_2 \leq ...\leq I$ and take $B =
A_{n+k} -A_n$ for some fixed values of $n,k \in \mathbb N$.
Observe that $0 \leq B \leq I$,
so Cauchy-Schwarz' inequality holds for the bilinear form
 $\langle Bx,y \rangle$. Use this to get:
$ \langle B x, Bx \rangle ^2 \leq \langle Bx,x
\rangle \langle B^2x, Bx \rangle
\leq \langle Bx,x \rangle  \langle Bx, Bx \rangle,
$
from which
$$ \| B x\|^2 = \langle B x, Bx \rangle  \leq  \langle Bx,x \rangle
$$
Thus, for every vector $x \in H$:
$$\|A_{n+k}x -A_nx\|^2 \ \leq \ \langle A_{n+k}x,x \rangle  - \langle
A_nx, x \rangle.$$
Since the sequence $\langle A_n x, x \rangle $ is bounded
and monotonic, it has a limit. Hence $\lim_n A_n x$ exists for
every $x\in H$, which proves the statement.\\

{\bf b).} {\it Every non-negative operator $A$ admits a unique
non-negative square root $\sqrt{A}$: $(\sqrt{A})^2 = A$.}\\

\smallskip

\noindent
For the proof one can normalize $A$, so that $0 \leq A \leq I$ and
use a convergent series decomposition for $\sqrt{x} =
\sqrt{1-(1-x)}$, in conjunction with the above remark. See for
details \S 104 of \cite{RN}.

Conversely, if $T \in L(H)$, then $T^\ast T \geq 0$.\\

{\bf c).} {\it Let $A,B$ be two commuting non-negative (linear
bounded) operators. Then $AB$ is also non-negative.}

\smallskip

\noindent
Note that, if $AB=BA$, the above proof implies $\sqrt{B} A =
A\sqrt{B}$. For the proof we compute directly
$$ \langle ABx,x\rangle = \langle A \sqrt{B} \sqrt{B}x,x\rangle
=$$
$$ \langle \sqrt{B} A \sqrt{B} x, x \rangle = \langle A
\sqrt{B} x,  \sqrt{B} x\rangle \geq 0.$$\\

With the above observations
we can enhance the polynomial functional
calculus of a self-adjoint operator.
Let
$\mathbb C[t], \mathbb R[t]$
\index{${\mathbb{C}}[t], \RR[t]$ polynomial algebras}
denote the algebra of polynomials with complex, respectively
real,  coefficients in one variable and let $A=A^\ast$ be a
self-adjoint operator with bounds (\ref{bounds}).
The expression $p(A)$ makes sense for
every $p \in \mathbb C[t]$, and the
{polynomial functional calculus for $A$}
\index{polynomial functional calculus}
which is the map $\phi$
$$p \stackrel{\phi}{\mapsto} p(A)$$
 is
obviously linear, multiplicative and unital (1 maps to I).
Less obvious is the  key fact that that $\phi$ is positivity
preserving:

\begin{prop} If the polynomial $p \in \mathbb R[t]$
satisfies $p(t) \geq 0$ for all  $t$ in $ [m, M]$ and the self-adjoint
operator $A$ satisfies $mI \leq A \leq MI$, then $p(A) \geq 0$.
\end{prop}

\begin{proof}.
 A decomposition of the real polynomial $p$ into
irreducible, real factors yields:
$$ p(t) = c \prod_i (t-\alpha_i) \; \prod_j (\beta_j - t)
\;  \prod_k
[(t-\gamma_k)^2 + \delta_k^2], $$
with
$c>0, \ \alpha_i \leq m \leq M \leq \beta_j$
and $\gamma_k \in \bbR  ,\delta_k \in \mathbb R.$
According to
the observation c) above, we find $p(A) \geq 0$.  \ \qed
\end{proof}

The proposition  immediately implies
\begin{cor}
The homomorphism $\phi$ on $\CC[t]$ extends to $C[m,M]$ and beyond.
Moreover,
$$\| p(A) \| \leq \sup_{[m,M]}|p|= :\|p \|_\infty .$$
\end{cor}

\begin{proof}.
The inequality follows because
$ \sup_{[m,M]}|p| \pm p $
is a  polynomial  nonnegative on $[m,M]$,
so $\|p \|_\infty I \geq \pm p(A)$
which gives the required inequality.
Thus $\phi$ is sup norm continuous
and extends by continuity to the completion of the polynomials,
which is of course the algebra $C[m,M]$ of the continuous functions.
\end{proof}

The Spectral Theorem \index{spectral theorem} immediately follows.

\begin{thm}
  If the self adjoint bounded operator $A$ on $H$
  has a cyclic vector $\xi$,
  then there is a positive Borel measure $\mu$ on $[m,M]$
  and a  unitary operator $U: H \mapsto L^2(\mu)$
  identifying $H$ with $L^2(\mu)$ such that
    $$ U A U^* = M_x .$$
    Here for any $g$ in $L^\infty$
    the multiplication operator     $M_g$ is defined by
      $M_g f =gf$ on all $f \in  L^2(\mu)$.
\end{thm}
\noindent
The vector $\xi$ {\it cyclic}\index{cyclic}  means
$$
span \; \{ A^k \xi : \ k= 0, 1, 2. \cdots \}=
\{ p(A) \xi : \ p \ a \ polynomial \ \}
$$
is dense in $H$.

\begin{proof} \
Define a linear functional
$L: C([m,M]) \mapsto \CC $  by
$$
L(f):= \langle f(A) \xi, \xi \rangle \qquad
for \ all \ f \in C([m,M]).
$$
The Representation Theorem (see Proposition \ref{Riesz representation}
for more detail) for such $L$ says
there is  a Borel measure  $\mu$ such that
$$ L(f) = \int_{[m,M]} f d\mu ;$$
moreover, $\mu$ is a positive measure because
if $f \geq 0$ on $[m,M]$, then $L(f)\geq 0$.
A critical feature is
 \beq
 \label{eq:HL2}
 \int p \overline{q} d \mu
  =
\langle p(A)  \xi,  q(A) \xi \rangle
\eeq
which holds, since $= L( p \overline{q})  =
\langle p(A) \overline{q}(A) \xi,  \xi \rangle$.
We have built our representing space
(using a formula which haunts the rest of this paper)
and now we identify $H$ with this space.

Define $U$ by $U p(A)\xi = p$ which specifies
it on a dense set (by the cyclic assumption)
provided  $U p_1(A)\xi=  U p_2(A)\xi$
implies $e(A)\xi:= p_1(A)\xi - p_2(A)\xi=0$;
in  other words,
$ 0=
\langle e(A)  \xi,  q(A) \xi \rangle
$ for all polynomials $q$.
Thus $0= \int e \overline{q} d \mu$,
so $e=0$ a.e. wrt $\mu$.
Now to properties of $U$:
\begin{enumerate}
  \item
  $U$ is isometric. (That is what (\ref{eq:HL2}) says.)
Thus $U$ extends to $H$ and has closed range.
\item
The range of $U$ is dense since it contains the polynomials.
\item
$U A p(A) \xi = x p(x) = x U p(A) \xi $ for all polynomials $p$.
By the density imposed by cyclicity for any $v$ in $H$ we have
$$ U A v = M_x Uv. $$
Note  the constrction gives $ U\xi  = 1$. \qed
\end{enumerate}
\end{proof}

\subsection{A bigger functional calculus and spectral measures}
Our
next aim is to consider a bounded, increasing sequence $p_n$
of real polynomial functions on the interval $[m,M]$ and define,
according to observation a):
$$ f(A)x = \lim p_n(A)x, \ \ x \in H,$$
where $f$ is a point-wise limit of $p_n$.
A standard argument
shows that, if $q_n$ is another sequence of polynomials,
monotonically converging on $[m,M]$ to $f$, then
$$ \lim q_n(A)x = \lim_n p_n(A)x,\ \ x \in H.$$
See for details \S 106 of \cite{RN}.
The new calculus $f \mapsto f(A)$
remains linear and multiplicative.

In particular, we can apply the above definition to the step
functions
$$ \chi_s (t) = \left\{ \begin{array}{cc}
                       1,& t \leq s,\\
                       0,& t>s.
                       \end{array}\right.$$
                       This yields a monotonic, operator valued
                       function
$$ F_A(s) = \chi_s(A),$$
with the additional properties $F_A(s) = F_A(s)^\ast = F_A(s)^2$
and
$$ F_A(s) = \left\{ \begin{array}{cc}
                       0,& s<m,\\
                       I,& s\geq M.
                       \end{array}\right.$$
With the aid of this {\it spectral scale}
\index{spectral scale}
one can interpret the functional calculus as an operator valued
Riemann-Stieltjes integral
$$ f(A) = \int_m^M f(t) d F_A(t).$$ The {\it spectral measure}
\index{spectral measure} $E_A$ of $A$ is the operator valued
measure associated to the monotonic function $F_A$, that is, after
extending the integral to Borel sets $\sigma$,
$$ E_A (\sigma) = \int_{\sigma \cap [m,M]} dF_A(t).$$
Thus $E_A(\sigma)$ is a family of mutually commuting orthogonal
projections, subject to the multiplicativity constraint
$$E_A(\sigma \cap \tau) = E_A(\sigma) E_A(\tau).$$
As a matter of
notation, we have then for every bounded, Borel measurable
function $f$:
\beq
\label{eq:specThm}
f(A) = \int_m^M f(t) E_A(dt).
\eeq
This is a form of the Spectral Theorem \index{spectral theorem}
which does not
assume cyclicity.

A good exercise for the reader is to identify the above objects in
the case of a finite dimensional Hilbert space $H$ and a
self-adjoint linear transformation $A$ acting on it.
A typical
infinite dimensional example will be discussed later in connection
with the moment problem.

\subsection{Unitary operators} The spectral theorem for a unitary
transformation \index{unitary operator}$U \in L(H), \ U^\ast U = U
U^\ast = I,$ can be derived in a very similar manner.

The needed structure of positive polynomials is contained in the
following classical result.

\begin{lem}[Riesz-Fej\'er] A non-negative trigonometric polynomial
is the modulus square of a trigonometric
polynomial.
\index{Riesz-Fej\'er Lemma}
\end{lem}

\begin{proof}. Let $p(e^{i\theta}) = \sum_{-d}^d c_j e^{ij\theta}$ and
assume that $p(e^{i\theta}) \geq 0, \ \ \theta \in [0,2 \pi]$.
Then necessarily $c_{-j} = \overline{c_j}$. By passing to complex
coordinates, the rational function $p(z) = \sum_{-d}^{d} c_j z^j$
must be identical to $\overline{p(1/\overline{z})}$. That is its
zeros and poles are symmetrical (in the sense of Schwarz) with
respect to the unit circle.

Write $z^d p(z) = q(z)$, so that $q$ is a polynomial of degree
$2d$. One finds, in view of the mentioned symmetry:
$$ q(z) = c z^\nu \prod_j (z-\lambda_j)^2 \prod_k
(z-\mu_k)(z-1/\overline{\mu_k}),$$ where $c \neq 0$ is a constant,
$|\lambda_j| = 1$ and $0<|\mu_k|<1$.

For $z = e^{i \theta}$ we obtain
$$ p(e^{i \theta}) = |p(e^{i \theta})| = |q(e^{i \theta}| = $$ $$
|c| \prod_j |e^{i \theta}-\lambda_j|^2 \prod_k \frac{|e^{i
\theta}-\mu_k|^2}{|\mu_k|^2}.$$ \qed
\end{proof}

Returning to the unitary operator $U$ we infer, for $p \in \mathbb
C[z]$,
$$ \Re p(e^{i \theta}) \geq 0 \ \ \Rightarrow \Re p(U) \geq 0.$$
Indeed, according to the above Lemma, $\Re p(e^{i\theta}) =
|q(e^{i\theta})|^2,$ whence $$\Re p(U) = q(U)^\ast q(U) \geq 0.$$
Then, exactly as in the preceding section one constructs the
spectral scale and spectral measure of $U$.

For an operator $T$ we denote its
``real part" and ``imaginary part" by
$\Re T = (T+T^\ast)/2$ and $\Im T =(T-T^\ast)/{2i}$.
\index{$ \Re T, \Im T$ real, imaginary part}

The reader will find other elementary facts (\`a la
Riesz-Fej\'er's Lemma) about the decompositions of non-negative
polynomials into sums of squares in the second volume of Polya and
Szeg\"o's problem book \cite{PolyaS}. This particular collection
of observations about positive polynomials reflects, from the
mathematical analyst point of view, the importance of the subject
in the first two decades of the XX-th century.

\subsection{Riesz-Herglotz formula}
\label{sec:ReiszHerg}
The practitioners of spectral analysis know that the strength and
beauty of the spectral theorem lies in the effective dictionary it
establishes between matrices, measures and analytic functions. In
the particular case of unitary operators, these correspondences
also go back to F. Riesz. The classical Riesz-Herglotz formula
\index{Riesz-Herglotz formula} is incorporated below in a more
general statement. To keep the spirit of positivity of the last
sections, we are interested below in the {\it additive} (rather
than multiplicative) structure of polynomials (or more general
functions) satisfying Riesz-Fej\'er's condition:
$$ \Re p(z) \geq 0, \ \ |z| <1.$$

We denote by $\mathbb D$ the unit disk in the complex plane.
Given a set $X$ by a
{\it positive semi-definite kernel} \index{positive semi-definite
kernel} we mean a function $K : X \times X \longrightarrow \mathbb
C$ satisfying
$$ \sum_{i,j=1}^N K(x_i, x_j) c_i \overline{c_j} \geq 0,$$
for every finite selection of points $x_1,...,x_N \in X$ and
complex scalars $c_1,...,c_N$.

\begin{thm}\label{RH} Let $f: \mathbb D \longrightarrow \mathbb C$ be an
analytic function. The following statements are equivalent:\\

a). $\Re f(z) \geq 0, \ \ z \in \mathbb D$,\\

b). (Riesz-Herglotz formula). There exists a positive Borel
measure $\mu$ on $[-\pi, \pi]$ and a real constant $C$, such that:
$$ f(z) = iC + \int_{-\pi}^\pi \frac{e^{it}+z}{e^{it}-z} d\mu(t),\
\ z \in \mathbb D,$$\\

c). The kernel $K_f : \mathbb D \times \mathbb D \longrightarrow
\mathbb C$,
$$ K_f(z,w) = \frac{f(z) + \overline{f(w)}}{1-z\overline{w}},\ \ z,w
\in \mathbb D,$$
is positive semi-definite,\\

d). There exists a unitary operator $U \in \mathcal L(H)$, a
vector $\xi \in H$ and a constant $a \in \mathbb C, \ \Re a \geq
0$, such that:
$$ f(z) = a + z\langle (U-z)^{-1}\xi, \xi\rangle, \ \ z \in
\mathbb D.$$
\end{thm}

\begin{proof}. We merely sketch the main ideas in the proof.
The reader can consult for details the monograph \cite{AM}.

$a) \Rightarrow b).$ Let $r<1$. As a consequence of Cauchy's
formula:
$$ f(z) = i \Im f(0) + \frac{1}{2\pi}\int_{-\pi}^\pi
\frac{re^{it}+z}{re^{it}-z}
\Re f(r e^{it}) dt,\ \ |z|<r.$$ Since the positive measures
$\frac{1}{2\pi} \Re f(r e^{it}) dt$ have constant mass on $[-\pi,
\pi]$:
$$ \frac{1}{2\pi} \int_{-\pi}^\pi \Re f(r e^{it}) dt = \Re f(0),\ \
r<1,$$ they form a weak$-\ast$ relatively compact family (in the
space of finite measure). Any weak$-\ast$ limit will satisfy the
identity in b) (hence all limit points coincide).

$b) \Rightarrow c)$. A direct computation yields:
\begin{equation}
\label{NP kernel} K_f(z,w) = \int_{-\pi}^\pi
\frac{2}{(e^{it}-z)(e^{-it}-\overline{w})} d\mu(t), \ \ z,w \in
\mathbb D.
\end{equation}
Since for a fixed value of $t$, the integrand is positive
semi-definite, and we average over a positive measure, the whole
kernel will turn out to be positive semi-definite.

$c) \Rightarrow a).$ Follows by evaluating $K_f$ on the diagonal:
$$ 2 \Re f(z) = (1-|z|^2) K_f(z,z) \geq 0.$$

$b) \Rightarrow d).$ Let $H = L^2(\mu)$ and $Uf(t)  = e^{it}f(t)$.
Then $U$ is a unitary operator, and the constant function $\xi =
\sqrt{2}$ yields the representation d).

$d) \Rightarrow b).$ In view of the spectral theorem, we can
evaluate the spectral measure $E_U$ on the vector $\xi$ and obtain
a positive measure $\mu$ satisfying:
$$f(z) = a + z\langle (U-z)^{-1}\xi, \xi\rangle = a + z
\int_{-\pi}^\pi \frac{d\mu(t)}{e^{it} -z}=
$$ $$ a+ \frac{1}{2} \int_{-\pi}^\pi \frac{e^{it}+z}{e^{it}-z}
d\mu(t) -\frac{1}{2} \int_{-\pi}^\pi d\mu(t) , \ \ z \in \mathbb
D.$$ By identifying the constants we obtain, up to the factor $2$,
conclusion b).  \ \qed
\end{proof}
\bigskip

The theorem above has far reaching consequences in quite divergent
directions: function theory, operator theory and control theory of
linear systems, see for instance \cite{AM, FF,M03, RR}. We confine
ourselves to describe only a generic consequence.

First, we recall that, exactly as in the case of finite matrices,
a positive semi-definite kernel can be written as a sum of
squares. Indeed, if $K : X \times X \longrightarrow \mathbb C$ is
positive semi-definite, one can define a sesqui-linear form on the
vector space $\oplus_{x \in X} \mathbb C$,  with basis $e(x), x
\in X$, by
$$ \| \sum_i c_i e(x_i) \|^2 = \sum_{i,j=1}^N K(x_i, x_j) c_i \overline{c_j}.$$
This is a positive semi-definite inner product. The associated
separated (i.e. Hausdorff) Hilbert space completion $H$ carries
the classes of the vectors $[e(x)] \in H$. They factor $K$ into a
sum of squares:
$$ K(x,y) = \langle [e(x)], [e(y)]\rangle = \sum_k \langle [e(x)], f_k \rangle
\langle f_k, [e(y)]\rangle,$$ where $(f_k)$ is any orthonormal
basis of $H$. For details, see for instance the Appendix to
\cite{RN}.

The following result represents the quintessential bounded
analytic interpolation theorem.

\begin{thm}[Nevanlinna-Pick]
\index{Nevanlinna-Pick Theorem}
Let $\{a_i \in \mathbb D; \ i \in
I\}$ be a set of points in the unit disk, and let $\{ c_i \in
\mathbb C; \ \Re c_i \geq 0, \  i \in I\}$ be a collection of
points in the right half-plane, indexed over the same set.

There exists an analytic function $f$ in the unit disk, with $\Re
f(z) \geq 0, \ |z|<1,$  and $f(a_i) = c_i, \  i \in I, $ if and
only if the kernel
$$ \frac{c_i + \overline{c_j}}{1-a_i \overline{a_j}}, \ \ i,j \in
I,$$ is positive semi-definite.
\end{thm}

\begin{proof}.
Point c) in the preceding Theorem shows that the condition is
necessary.

A Moebius transform in the range $(f \mapsto g = (f-1)/(f+1))$
will change the statement into:
$$ g: \mathbb D \longrightarrow \mathbb D, \ \ g(a_i) = d_i,$$
if and only if the kernel
$$ \frac{1-d_i \overline{d_j}}{1-a_i \overline{a_j}}, \ \ i,j \in
I,$$ is positive semi-definite.

To prove that the condition in the statement is also sufficient,
assume that the latter kernel is positive semi-definite. As
before, factor it (into a sum of squares):
$$\frac{1-d_i \overline{d_j}}{1-a_i \overline{a_j}} = \langle
h(i), h(j) \rangle, \ \ i,j \in I,$$ where $h:I \longrightarrow H$
is a function with values in an auxiliary Hilbert space $H$.

Then $$ 1 + \langle a_i h(i), a_j h(j) \rangle = d_i\overline{d_j}
+ \langle h(i), h(j) \rangle, \ \ i,j \in I.$$ The preceding
identity can be interpreted as an equality between scalar products
in $\mathbb C \oplus H$:
$$ \langle \left( \begin{array}{c} 1\\
                                    a_i h(i)\\
                                    \end{array}
                                    \right) , \left(
                                    \begin{array}{c} 1\\
                                    a_j h(ij)\\
                                    \end{array}
                                    \right) \rangle =
\langle \left( \begin{array}{c} d_i\\
                                    h(i)\\
                                    \end{array}
                                    \right), \left(
                                    \begin{array}{c} d_j\\
                                    h(j)\\
                                    \end{array}
                                    \right) \rangle, \ \ i,j \in
                                    I.$$
Let $ H_1 \subset \mathbb C \oplus H$ be the linear span of the
vectors $(1, a_i h(i))^T,\ \ i \in I$. The map
$$ V \left( \begin{array}{c} 1\\
                                    a_i h(i)\\
                                    \end{array}
                                    \right) = \left( \begin{array}{c} d_i\\
                                    h(i)\\
                                    \end{array}
                                    \right)$$
                                    extends then by linearity to
                                    an isometric transformation
                                    $V:H_1 \longrightarrow H$.
Since the linear isometry $V$ can be extended (for instance by
zero on the orthogonal complement of $H_1$)  to a contractive
linear operator $T : \mathbb C \oplus H \longrightarrow \mathbb C
\oplus H$, we obtain a block matrix decomposition of $T$
satisfying:
$$ \left[ \begin{array}{cc}
                      A & B\\
                      C & D\\
                      \end{array} \right]\left( \begin{array}{c} 1\\
                                    a_i h(i)\\
                                    \end{array}
                                    \right) = \left( \begin{array}{c} d_i\\
                                    h(i)\\
                                    \end{array}
                                    \right).$$
Since $\| D\| \leq 1$, the operator $I-zD$ is invertible for all
$z \in \mathbb D$. From the above equations we find, after
identifying $A$ with a scalar:
$$ h(i) = (I-a_i D)^{-1} C 1,  \ \ \ d_i = A + a_i B h(i).$$

We define the analytic function
$$ g(z) = A + z B (I-zD)^{-1}C 1, \ \ |z|<1.$$
It satisfies, as requested: $g(a_i) = d_i, \  i \in I$.

By reversing the above reasoning we infer, with $h(z) =
(I-zD)^{-1}C 1 \in H$:
$$\left[ \begin{array}{cc}
                      A & B\\
                      C & D\\
                      \end{array} \right]\left( \begin{array}{c} 1\\
                                    z h(z)\\
                                    \end{array}
                                    \right) = \left( \begin{array}{c} g(z)\\
                                    h(z)\\
                                    \end{array}
                                    \right).$$
Since $T$ is a contraction,
$$ \| g(z) \|^2 + \| h(z) \|^2 \leq 1 + \| z h(z)\|^2 \leq 1 + \|
h(z) \|^2, \ \ |z|<1,$$ whence
$$ |g(z)| \leq 1, \ \  |z|<1.$$  \ \qed
\end{proof}
\bigskip

The above proof contains the germ of what experts in control
theory call ``realization theory". For the present survey it is
illustrative as a constructive link between matrices and analytic
functions with bounds; it will also be useful as a model to follow
in more general, non-commutative settings.

A great deal of research was done in the last two decades on
analogs of Riesz-Herglotz type formulas in several complex
variables. As expected, when generalizing to $\mathbb C^n$, there
are complications and surprises on the road. See for instance
\cite{AM, BT, CW, EP} and in several non-commuting variables
\cite{BGM05,K05}. We will return to some of these topics from
the perspective of positive polynomials and moment sequences.

\subsection{von Neumann's inequality}
\label{sec:vonN}
We have just seen that the
heart of the spectral theorem for self-adjoint or unitary
operators was the positivity of the polynomial functional
calculus. A surprisingly general inequality, of the same type,
applicable to an arbitrary bounded operator, was discovered by von
Neumann \cite{vN2}.

\begin{thm}
\index{von Neumann's inequality}
Let $T \in \mathcal L(H), \| T \| \leq 1,$ be a contractive
operator. If a polynomial $p \in \mathbb C[z]$ satisfies $ \Re
p(z) \geq 0, \ z \in \mathbb D,$ then $\Re p(T) \geq 0$.
\end{thm}

\begin{proof}.
According to Riesz-Herglotz formula we can write
$$ p(z) = iC + \int_{-\pi}^\pi \frac{e^{it}+z}{e^{it}-z} d\mu(t),
\ |z|<1,$$ where $C \in \mathbb R$ and $\mu$ is a positive
measure.

Fix $r<1$, close to $1$, and evaluate the above representation at
$z=rT$:
$$ p(rT) = iC + \int_{-\pi}^\pi (e^{it}+ rT)(e^{it}-rT)^{-1}
d\mu(t).$$ Therefore
$$ p(rT) + p(rT)^\ast = $$ $$\int_{-\pi}^\pi (e^{it}-rT)^{-1}[(e^{it}+
rT)(e^{-it}- rT^\ast)+(e^{it}- rT)(e^{-it}+
rT^\ast)](e^{-it}-rT^\ast)^{-1} d\mu(t) = $$
$$2 \int_{-\pi}^\pi (e^{it}-rT)^{-1}[I - r^2 T
T^\ast](e^{-it}-rT^\ast)^{-1} d\mu(t)
\geq 0.$$

Letting $r \rightarrow 1$ we find $\Re p(T) \geq 0$. \  \ \qed
\end{proof}
\bigskip

A Moebius transform argument, as in the proof of Nevanlinna-Pick
Theorem, yields the equivalent statement (for a contractive linear
operator $T$):
$$ (|p(z)| \leq 1, \ \ |z|<1)\  \Rightarrow \ \| p(T) \| \leq 1.$$

Von Neumann's original proof relied on the continued fraction
structure of the analytic functions from the disk to the disk. The
recursive construction of the continued fraction goes back to
Schur \cite{Schur} and can be explained in a few lines.\\

{\bf Schur's algorithm.}
\index{Schur's Algorithm}
Let $f: \mathbb D \longrightarrow \mathbb D$
be an analytic function.
Then, in
view of Schwarz Lemma, there exists an analytic function $f_1
:\mathbb D \longrightarrow \mathbb D$ with the property:
$$ \frac{f(z) - f(0)}{1-\overline{f(0)}f(z)} = z f_1(z),$$
or equivalently, writing $s_0 = f(0)$:
$$ f(z) = \frac{s_0 + zf_1(z)}{1+\overline{s_0} zf_1(z)}.$$
In its turn,
$$ f_1(z) = \frac{s_1 + zf_2(z)}{1+\overline{s_1} zf_2(z)},$$
with  an analytic $f_2 :\mathbb D \longrightarrow \mathbb D$, and
so on.

This algorithm terminates after finitely many iterations for finite Blashcke products
$$ f(z) = \prod_{k=1}^N \frac{z-\lambda_k}{1-\overline{\lambda_k}z},
\ \ |\lambda_k|<1.$$
Its importance lies in the fact that the finite section of Schur
parameters $(s_0, s_1,...,s_n)$ depends via universal expressions
on the first section (same number) of Taylor coefficients of $f$
at $z=0$. Thus, the conditions
$$ |s_0(c_0)|\leq 1, \ |s_1(c_0,c_1)| \leq 1, \ldots $$
characterize which power series
$$c_0 + c_1z + c_2 z^2 +...,$$
are associated to analytic functions from the disk to the disk.
For details and a variety of applications, see \cite{C, FF, RR}.

One notable application  is to solve the classical
\index{Carath\'eodory-Fej\'er interpolation problem}
Carath\'eodory-Fej\'er interpolation problem, a close relative of
the Nevanlinna-Pick problem we presented earlier. Here one
specifies complex numbers $c_0, \cdots, c_m$ and seeks $f: \bbD
\to \bbD$ analytic for which
$$
\frac{1}{j!} \frac{d^j f}{dz^j}(0) = c_j,  \qquad j= 0, \cdots, m.
$$
The Schur Algorithm constructs such a function and in the same
time gives a simple criterion when the solution exists.
Alternatively, a special type of matrix $(c_{n-m})_{n,m=0}^{m}$,
with zero entries under the diagonal $(c_j =0, \ j<0)$, called a Toeplitz matrix,
based on $c_0, \cdots, c_m$ is a contraction if and only if a
solution to the Carath\'eodory-Fej\'er problem exists. A version
of this fact in the right half plane (rather than the disk) is
proved in Theorem \ref{trig problem}.

As another application, we can derive (also following
Schur)
\index{separation of zeros}
an effective criterion for
deciding whether a polynomial has all roots inside the unit disk.
Let
$$ p(z) = c_d z^d + c_{d-1} z^{d-1} + ...+ c_0 \in \mathbb C[z],$$
and define
$$ p^\flat(z) = z^d \overline{p(1/\overline{z})} = \overline{c_0}
z^d + \overline{c_1} z^{d-1} + ...+ \overline{c_d}.$$ It is clear
that
$$ |p(e^{it})| = |p^\flat(e^{it})|, \ \ t \in [-\pi, \pi],$$
and that the roots of $p^\flat$ are symmetric with respect to the
unit circle to the roots of $p$. Therefore, $p$ has all roots
contained in the open unit disk if and only if $\frac{p}{p^\flat}$
is an analytic function from the disk to the disk, that is, if and
only if the kernel
$$ \frac{p^\flat(z) \overline{p^\flat(w)} - p(z)
\overline{p(w)}}{1-z\overline{w}}, \ \ z,w \in \mathbb D,$$ is
positive definite. As a matter of fact $\frac{p}{p^\flat}$ is a
finite Blashcke product, and Schur's algorithm terminates in this case after finitely many
iterations.

In general, regarded as a Hermitian form, evaluated to the
variables $Z_i = z^i, 0 \leq i \leq d$, the signature of the
above kernel (that is the number of zeros, negative and positive
squares in its canonical decomposition) counts how many roots the
polynomial $p$ has inside the disk, and on its boundary. For many
more details see the beautiful survey \cite{KN}.

\section{Moment problems}
\label{sec:moments}
In this section we return to Hilbert space and the
spectral theorem, by unifying the analysis and algebra concepts we
have discussed in the previous sections. This is done in the
context of power moment problems, one of the oldest and still
lively sources of questions and inspiration in mathematical
analysis.

As before, $x =(x_1,...,x_g)$ stands for the coordinates in
$\mathbb R^g$, and, at the same time, for a tuple of commuting
indeterminates. We adopt the multi-index notation $x^\alpha =
x_1^{\alpha_1}... x_g^{\alpha_g}, \ \alpha \in \mathbb N^g$. Let
$\mu$ be a positive, rapidly decreasing measure on $\mathbb R^g$.
The {\it moments} \index{moments} of $\mu$ are the real numbers:
$$ a_\alpha = \int x^\alpha d\mu(x), \ \ \alpha \in \mathbb N^g.$$
For its theoretical importance and wide range of applications, the
correspondence
$$ \{ \mu; \ {\rm positive \  measure} \} \ \longrightarrow \ \{
(a_\alpha); \ {\rm moment \ sequence} \}$$ can be put on an equal
level with the Fourier-Laplace, Radon or wavelet transforms. It is
the positivity of the original measure which makes the analysis of
this category of moment problems interesting and non-trivial, and
appropriate for our discussion. For general aspects and
applications of moment problems (not treated below) the reader can
consult the monographs \cite{A,BCR,FF,ST} and the excellent survey
\cite{Fuglede}. The old article of Marcel Riesz \cite{MRiesz}
remains unsurpassed for the classical aspects of the one variable
theory.

Given a multi-sequence of real numbers $(a_\alpha)_{\alpha \in
\mathbb N^g}$ a linear functional representing the potential
integral of polynomials can be defined as:
$$ L : \mathbb R[x] \longrightarrow \mathbb R, \ \ L(x^\alpha) =
a_\alpha, \  \alpha \in \mathbb N^g,$$
and vice-versa.
When necessary we will complexify $L$ to a complex linear functional on
$\mathbb C[x]$.

  If $(a_\alpha)_{\alpha \in \mathbb N^g}$ are the
moments of a positive measure, then for a polynomial $p \in
\mathbb R[x]$ we have
$$ L (p^2) = \int_{\mathbb R^g} p^2 d\mu \geq 0.$$
Moreover, in the above positivity there is more structure: we can
define on $\mathbb C[x]$ a pre-Hilbert space bracket by:
$$ \langle p, q \rangle = L(p \overline{q}), \ \ p,q \in \mathbb
C[x].$$ The inner product is positive semi-definite, hence
the Cauchy-Schwarz inequality holds:
$$
|\langle p, q \rangle|^2 \leq \| p \|^2 \| q \|^2.
$$
Thus, the set of null-vectors $N = \{ p \in \mathbb C[x];\ \| p \|
=0 \}$ is a linear subspace, invariant under the multiplication by
any polynomial. Let $H$ be the Hilbert space completion of
$\mathbb C[x]/N$ with respect to the induced Hermitian form. Let
$\mathcal D = \mathbb C[x]/N$ be the image of the polynomial
algebra in $H$. It is a dense linear subspace, carrying the
multiplication operators:
$$ M_{x_i} : \cD \longrightarrow \cD, \ \ M_{x_i} p = x_i p.$$
Note that these are well defined, symmetric linear operators:
$$ \langle  M_{x_i}p, q \rangle = L(x_i p \overline{q}) = \langle p,
M_{x_i}q \rangle, \ \ p,q \in \cD,$$ and they commute
$$ M_{x_i} M_{x_j} = M_{x_j} M_{x_i}.$$
Finally the (constant function) vector $\xi = 1$ is
cyclic\index{cyclic vector}, in the sense that $\cD$ is the linear
span of repeated actions of $M_{x_1},...,M_{x_g}$ on $\xi$:
$$ \cD = \bigvee_{\alpha \in \mathbb N^g} M_{x_1}^{\alpha_1} ...
M_{x_g}^{\alpha_g} \xi.$$
We collect these observations into a single statement.

\begin{prop}
\label{general OT correspondence}
There is a bijective
correspondence between all
linear functionals
$$
L \in  \mathbb R[x]', \ \
L|_{\Sigma^2 \mathbb R[x]} \geq 0,$$
and the pairs $(M,\xi)$ of
$g$-tuples $M=(M_1,...,M_g)$ of commuting, symmetric linear
operators with a cyclic vector $\xi$ (acting on a separable
Hilbert space).
The correspondence is given by the relation
$$ L(p) = \langle p(M)\xi, \xi \rangle, \ \ p \in \mathbb R[x].$$
\end{prop}

\noindent
Above the word commuting has to be taken with caution: implicitly
it is understood that we define the span $\cD$ as before, and
remark that every $M_i$ leaves $\cD$ invariant. Then $M_i$
commutes with $M_j$ as endomorphisms of $\cD$.

Having a positive measure $\mu$ represent the functional $L$ adds
in general new constraints in this dictionary.

Let $\cP_+(K)$
\index{$\cP_+(K)$ \ polynomials $\geq 0$ on  $ K$}
be the  set of all
polynomials which are non-negative on the set
$ K \subset \mathbb{R}^g$ and note that this is a convex cone.

\begin{prop}
\label{Riesz representation}
A linear functional $L \in
\mathbb R[x]'$ is
representable by a positive measure $\mu$:
$$ L(p) = \int p d\mu, \ \ p \in \RR[x]$$
if and only if $L|_{\cP_+(\RR^g)} \geq 0.$
\end{prop}

\noindent Although this observation (in several variables) is
attributed to Haviland, see \cite{A}, it is implicitly contained
in Marcel Riesz article \cite{MRiesz}. Again we see exactly
{\it the  gap
$$ \Sigma^2 \mathbb R[x] \subset \cP_+(\mathbb R^g),$$
which we must understand in order to characterize the
moments of positive measures}
(as already outlined in Minkowski's and Hilbert's early works).

\bigskip
\begin{proof}. If the functional $L$ is represented by a positive
measure, then it is obviously non-negative on all non-negative
polynomials.

To prove the converse, assume that $L|_{\cP_+(\mathbb R^g)} \geq
0$. Let $C_{pBd}(\mathbb R^g)$ \index{$C_{pBd}(\mathbb R^g)$}
be the space of continuous functions
$f$ having a polynomial bound at infinity:
$$ |f(x)| \leq C (1+|x|)^N,$$
with the constants $C, N>0$ depending on $f$. We will extend $L$,
following M. Riesz \cite{MRiesz}, to a non-negative functional on
$C_{pBd}(\mathbb R^g)$.

This extension process, \index{M. Riesz extension} parallel and
arguably prior to the Hahn-Banach Theorem, works as follows.
Assume that
$$ \hat{L} : V \longrightarrow \mathbb R$$
is a positive extension of $L$ to a vector subspace $V \subset
C_{pBd}(\mathbb R^g)$. That is:
$$ (h \in V, \    h \geq 0) \ \Rightarrow (\hat{L}(h) \geq
0).$$
Remark that $L$ is defined on all polynomial functions.
Assume $V$ is not  the whole space and choose a
non-zero function $f \in C_{pBd}(\mathbb R^g) \setminus V$. Since $f$
has polynomial growth, there are elements $h_1, h_2 \in V$
satisfying
$$ h_1 \leq f \leq h_2.$$
By the positivity of $\hat{L}$,
we see $\hat{L}h_1 \leq  \hat{L}f \leq \hat{L}h_2$, that
is
$$ \sup_{ h_1 \leq f} \hat{L}(h_1) \leq \inf_{f \leq h_2}
\hat{L}(h_2).$$ Choose any real number $c$ between these limits
and define
$$ L'(h+ \lambda f) = \hat{L}(h) + \lambda c, \ \ h \in V, \
\lambda \in \mathbb R.$$ This will be a positive extension of $L$
to the larger space $V \oplus \mathbb R f$.

By a standard application of Zorn's Lemma, we find a positive
extension of $L$ to the whole space. Finally, F. Riesz
Representation Theorem provides a positive measure $\mu$ on
$\mathbb R^g$, such that $L(p) = \int p d\mu, \ \ p \in \mathbb
R[x].$ \ \qed
\end{proof}

Next we focus on  a few particular contexts (either low
dimensions, or special supporting sets for the measure)
where the structure of the
positive functionals and tuples of operators appearing in our
dictionary can be further understood.

\subsection{The trigonometric moment problem} We specialize to
dimension $n=2$ and to measures supported on the unit circle
(torus) $\mathbb T = \{ z \in \mathbb C; \ \ |z|=1\}.$ The group
structure of $\mathbb T$ identifies our moment problem to the
Fourier transform. It is convenient in this case to work with
complex coordinates $z = x+iy \in \mathbb C = \mathbb R^2$, and
complex valued polynomials. In general, we denote by $\Sigma^2_h
\mathbb C[x]$ the sums of moduli squares (i.e. $|q|^2$) of complex
coefficient polynomials.

The ring of regular functions on the torus is
$$A= \mathbb C[z,
\overline{z}]/(1-z\overline{z}) = \mathbb C[z] \oplus \overline{z}
\mathbb C[\overline{z}],$$
where $(1-z\overline{z})$ denotes the ideal generated by
$1-z\overline{z}$.
A non-negative linear functional $L$ on
$\Sigma^2_h A$ necessarily satisfies
$$ L(\overline{f}) = \overline{L(f)}, \ \ f \in A.$$
Hence
$L$ is determined by the complex moments $L(z^n), \  n \geq
0$.
The following result gives a satisfactory solution to the
trigonometric moment problem on the one dimensional torus.

\begin{thm}\label{trig problem}
\index{Trigonometric moment problem}
  Let $(c_n)_{n=-\infty}^\infty$ be a sequence of complex
numbers subject to the conditions $c_0 \geq 0, \ \ c_{-n} =
\overline{c_n}, \ n \geq 0$.
The following assertions are equivalent:\\

a). There exists a unique positive measure $\mu$ on $\mathbb T$,
such that:
$$ c_n = \int_{\mathbb T} z^n d\mu(z), \ \ n \geq 0;$$\\

b). The Toeplitz matrix \index{Toeplitz matrix}\
$(c_{n-m})_{n,m=0}^\infty$ is positive
semi-definite;\\

c). There exists an analytic function $F : \mathbb D
\longrightarrow \mathbb C, \ \Re F \geq 0,$ such that
$$ F(z) = c_0 + 2 \sum_{k=1}^\infty c_{-k} z^k, \ \ |z|<1;$$\\

d). There exists a unitary operator $U \in L(H)$ and a  vector
$\xi \in H$ cyclic for the pair $(U,U^\ast)$, such that
$$ \langle U^n \xi, \xi \rangle = c_n, \ \ n \geq 0.$$
\end{thm}

\begin{proof}.\
  Let $L : \mathbb C[z,
\overline{z}]/(1-z\overline{z}) \longrightarrow \mathbb C$ be the
linear functional defined by
$$ L(z^n) = c_n, \ n \geq 0.$$
Condition b) is equivalent to
$$ L(|p|^2) \geq 0, \ \ p \in \mathbb C[z,
\overline{z}]/(1-z\overline{z}).$$

Indeed, assume that $p(z) = \sum_{j=0}^g \alpha_j z^j.$ Then,
since $\overline{z} z =1$,
$$ |p(z)|^2 = \sum_{j,k =0}^g \alpha_j \overline{\alpha_k}
z^{j-k},$$ whence
$$ L(|p|^2) = \sum_{j,k=0}^g \alpha_j \overline{\alpha_k} c_{j-k}.$$
Thus $a) \Rightarrow b)$ trivially. In view of the Riesz-Fej\'er
Lemma, the functional $L$ is non-negative on all non-negative
polynomial functions on the torus. Hence, in view of Proposition
\ref{Riesz representation} it is represented by a positive
measure. The uniqueness is assured by the compactness of $\mathbb
T$ and Stone-Weierstrass Theorem (trigonometric polynomials are
uniformly dense in the space of continuous functions on $\mathbb
T$). The rest follows from Theorem \ref{RH}. \ \qed
\end{proof}

Notable in the above Theorem is the fact that the main objects are
in bijective, and constructive, correspondence established
essentially by Riesz-Herglotz formula. Fine properties of the
measure $\mu$ can be transferred in this way into restrictions
imposed on the generating function $F$ or the unitary operator
$U$.

For applications and variations of the above result (for instance
a matrix valued analog of it) the reader can consult \cite{AM, A,
FF, RR}.

\subsection{Hamburger's moment problem} The passage from the
torus to the real line reveals some unexpected turns, due to the
non-compactness of the line. One may argue that the correct analog
on the line would be the continuous Fourier transform. Indeed, we
only recall that Bochner's Theorem provides an elegant
characterization of the Fourier transforms of positive measures.

Instead, we remain consistent and study polynomial functions and
positive measures acting on them.
Specifically, consider an $\mathbb R$-linear functional
$$ L : \mathbb R[x] \longrightarrow \mathbb R,
  \ L|_{\Sigma^2 \mathbb R[x] }
\geq 0.$$
By denoting
$$ c_k = L(x^k), \ \ k \geq 0,$$
the condition $L|_{\Sigma^2 \mathbb R[x] }$
is equivalent to the positive
semi-definiteness of the Hankel matrix \index{Hankel matrix} \  \
$$ (c_{k+l})_{k,l=0}^\infty \geq 0,$$
since
$$
0 \leq \sum_{k,l} f_k c_{k+l} f_l =  \sum_{k,l} L(f_k x^k x^l  f_l)
=   L(  \sum_{k}f_k x^k    \sum_{l} x^l  f_l)
= L ( f(x)^2 ).
$$
Next use that every non-negative polynomial on the line is a sum of
squares of polynomials, to invoke Proposition \ref{Riesz
representation} for the proof of the following classical fact.

\begin{thm}[Hamburger]\index{Hamburger's Theorem} Let
$(c_k)_{k=0}^\infty$ be a sequence of real numbers. There exists a
rapidly decaying, positive measure $\mu$ on the real line, such
that
$$ c_k = \int_{-\infty}^\infty x^k d\mu(x), \ k \geq 0,$$
if and only if the matrix $(c_{k+l})_{k,l=0}^\infty$ is positive
semi-definite.
\end{thm}

Now we sketch  {\it a second proof of Hamburger Theorem},
based on the Hilbert space construction
we have outlined in the previous section.
Namely, start with the positive semi-definite matrix
$(c_{k+l})_{k,l=0}^\infty$ and construct a Hilbert space
(Hausdorff) completion $H$ of $\mathbb C[x]$, satisfying
$$ \langle x^k, x^l \rangle = c_{k+l}, \ \ k,l \geq 0.$$
Let $\mathcal D$ denote as before the  image of the algebra of
polynomials in $H$; the image is dense. The (single)
multiplication operator $$ (M p)(x) = x p(x), \ \ p \in \cD,$$ is
symmetric and maps $\cD$ into itself. Moreover, $M$ commutes with
the complex conjugation symmetry of $H$:
$$ \overline{Mp} = M\overline{p}.$$
By a classical result of von-Neumann \cite{vN1} there exists a
self-adjoint (possibly unbounded) operator $A$ which extends $M$
to a larger domain. Since $A$ possesses a spectral measure $E_A$
(exactly as in the bounded case), we obtain:
$$ c_k = \langle x^k, 1 \rangle = \langle M^k 1, 1 \rangle = $$
$$ \langle A^k 1, 1 \rangle = \int_{-\infty}^\infty x^k \langle
E_A(dx)1, 1\rangle.$$ The measure $\langle E_A(dx)1, 1\rangle$ is
positive and has prescribed moments $(c_k)$. \qed

\medskip
This second proof offers more insight into the uniqueness part of
Hamburger's problem. Every self-adjoint extension $A$ of the
symmetric operator $M$ produces a solution $\mu (dx) = \langle
E_A(dx)1,1 \rangle$. The set $K$ of all positive measures with
prescribed moments $(c_k)$ is  convex and compact in the
weak-$\ast$ topology. The subset of Nevanlinna extremal elements
of $K$ are identified with the measures $\langle E_A(dx)1,1
\rangle$ associated to the self-adjoint extensions $A$ of $M$.
In
particular one proves in this way the following useful uniqueness
criterion.

\begin{prop} Let $(c_k)$ be the moment sequence of a positive
measure $\mu$ on the line. Then a positive measure with the same
moments coincides with $\mu$ if and only if the subspace
$$ (iI+M)\cD \ \ {\rm is \ dense \ in }\  H,$$
or equivalently, there exists a sequence of polynomials $p_n \in
\mathbb C[x]$ satisfying
$$ \lim_{n \rightarrow \infty} \int^\infty_{-\infty}
 |(i+x)p_n(x) - 1|^2 \; d\mu(x) =
0.$$
\end{prop}
Note that both conditions are intrinsic in terms of the initial
data $(c_k)$. For the original function theoretic proof see
\cite{MRiesz}. For the operator theoretic proof see for instance
\cite{A}.

There exists a classical analytic function counterpart of the
above objects, exactly as in the previous case
(see \S \ref{sec:ReiszHerg}, \S \ref{sec:vonN} ) of the unit circle.
Namely, assuming that
$$ c_k = \langle A^k 1, 1\rangle =
\int_{-\infty}^\infty x^k d\mu(x), \ \ k \geq 0,$$
as before, the analytic function
$$F(z) = \int^\infty_{-\infty} \frac{d\mu(x)}{x-z} = \langle (A-z)^{-1} 1, 1
\rangle$$ is well defined in the upper half-plane $\Im z >0$ and
has the asymptotic expansion at infinity (in the sense of
Poincar\'e, uniformly convergent in wedges $0 < \delta < \arg z <
\pi - \delta$):
$$ F(z) \approx-\frac{c_0}z-\frac{c_1}{z^2}
-\cdots,\quad\Im(z)>0.$$

One step further, we have a purely algebraic recursion which
determines the continued fraction development
$$-\frac{c_0}z-\frac{c_1}{z^2}-\cdots=-\cfrac{c_0}{z-\alpha_0
-\cfrac{\beta_0}{z-\alpha_1-\cfrac{\beta_1}{z-\alpha_2-\
\cfrac{\beta_2}{\ddots}}}},\quad\alpha_k \in \mathbb R, \ \beta_k
\geq 0.$$

It was Stieltjes, and then Hamburger, who  originally remarked
that $(c_k)$ is the moment sequence of a positive measure if and
only if the elements $\beta_k$ in the continued fraction
development of the generating (formal) series are non-negative.
Moreover, in this case they proved that there exists a unique
representing measure if and only if the continued fraction
converges in the upper half-plane. For details and a great
collection of classical examples see Perron's monograph
\cite{Perron}. A well known uniqueness criterion was obtained via
this formalism by Carleman \cite{Carleman}. It states that
uniqueness holds if
$$ \sum_1^\infty \frac{1}{c_{2k}^{1/(2k)}} = \infty.$$
The condition is however not necessary for uniqueness.

The alert reader has seen the great kinship between the continued
fraction  recursion
 just elucidated and the recursion called the Schur Algorithm in
\S \ref{sec:vonN}. These are essentially the same thing, but one
is in the disk setting while the other is in the half plane.

\subsubsection{Moments on the semiaxis $[0, \infty]$}
The above picture applies with minor modifications to Stieltjes
problem, \index{Stieltjes moment problem} that is the power moment
problem on the semi-axis $[0,\infty)$.

\begin{example}
We reproduce below an example found by Stieltjes, and refined by
Hamburger. See for details \cite{Perron}. Let $\rho$ and $\delta$
be positive constants, and denote
$$\alpha=\frac1{2+\delta},\quad\gamma=\rho^{-\alpha}.$$
Then
$$a_n=(2+\delta)\rho^{n+1}\Gamma[(2+\delta)(n+1)]
=\int_0^\infty x^n
e^{-\gamma x^\alpha}\>dx,\quad n\ge0,$$ is a moment sequence on
the positive semi-axis. A residue integral argument implies
$$\int_0^\infty x^n\sin\bigl(\gamma x^\alpha\tan(\pi\alpha)\bigr)\,
e^{-\gamma x^\alpha}\>dx=0,\quad n\ge0.$$
 Hence
$$a_n=\int_0^\infty x^n\bigl(1+t\sin(\gamma x^\alpha\tan(\pi\alpha))
\bigr)e^{-\gamma x^\alpha}\>dx,$$ for all $n\ge0$ and
$t\in(-1,1)$. This shows that the moment sequence $(a_n)$
does not uniquely
determine $\mu$ even knowing its  support is $[0,\infty)$.
\end{example}

Summing up the above ideas, we have bijective
correspondences between the following sets ($\mathbb C_+$ stands for
the open upper half plane):\\

  A). {\it Rapidly decaying positive measures $\mu$ on the real line;\\

  B). Analytic functions $F: \mathbb C_+ \longrightarrow
\overline{\mathbb C_+}$,
  satisfying $\sup_{t>1} |t F(it)| < \infty$;\\

  C). Self-adjoint operators $A$ with a cyclic vector $\xi$.}\\

  More precisely:
  $$ F(z) = \langle (A-z)^{-1} \xi, \xi \rangle
  = \int^\infty_{-\infty}
  \frac{d\mu(x)}{x-z}, \ \ z \in \mathbb C_+.
  $$
  The moment sequence $c_k = \int^\infty_{-\infty}
   x^k d\mu(x), k \geq 0,$
  appears in the asymptotic expansion of $F$, at infinity, but it does
not determine $F$, $(A,\xi)$ or $\mu$. For further details about
Hamburger and Stieltjes moment problems see Akhiezer's monograph
\cite{A}.

\subsection{Several variables} The moment problem on $\mathbb R^g,
\ \ g >1$, is considerably more difficult and less understood.
Although we have the general correspondence remarked in
Proposition \ref{general OT correspondence}, the gap between a
commuting tuple of unbounded symmetric operators and a strongly
commuting one (i.e. by definition one possessing a joint spectral
measure) is quite wide. A variety of strong commutativity criteria
came to rescue; a distinguished one, due to Nelson \cite{Nelson},
is worth mentioning in more detail.
\index{Nelson's self-adjoint extension criterion}

Assume that $L : \mathbb R [x_1,...,x_g] \longrightarrow \mathbb
R$ is a functional satisfying (the non-negative Hankel form
condition) $L|_{\Sigma^2 \mathbb R [x]} \geq 0$. We complexify $L$
and associate, as usual by now, the Hilbert space $H$ with inner
product:
$$ \langle p, q \rangle = L(p \overline{q}), \ \ \mathbb C[x].$$
The symmetric multipliers $M_{x_k}$ commute on the common dense
domain $\mathcal D = \mathbb C[x] \subset H$.
Exactly as in the
one variable case, there exists a positive measure $\mu$ on
$\mathbb R^g$ representing $L$ if and only if there are
(possibly unbounded)
self-adjoint extensions $M_{x_k} \subset A_k, \  1 \leq k \leq n,$
commuting at the level of their resolvents:
$$ [ (A_k-z)^{-1}, (A_j-z)^{-1}] :=
(A_k-z)^{-1} (A_j-z)^{-1} \
- \ (A_j-z)^{-1} (A_k-z)^{-1}= 0,$$
$$
\ \ {\rm for} \ \Im z >0, \ 1 \leq j,k \leq n.$$ See for details
\cite{Fuglede}. Although individually every $M_{x_k}$ admits at
least one self-adjoint extension, it is the joint strong
commutativity (in the resolvent  sense) of the extensions needed
to solve the moment problem.

Nelson's theorem gives a sufficient condition in this sense: if
$(1 + x_1^2 +...+ x_g^2) \mathcal D$ is dense in $H$, then the
tuple of multipliers $(M_{x_1},...,{M_{x_g}})$ admits an extension
to a strongly commuting tuple of self-adjoint operators. Moreover,
this insures the uniqueness of the representing measure $\mu$. For
complete proofs and more details see \cite{Berg1, Fuglede}.

A tantalizing open question in this area can be phrased as
follows:\\

{\bf Open problem.}
{\it Let $(c_{\alpha+ \beta})_{\alpha, \beta \in \mathbb N^g}$
be a positive semi-definite Hankel form. Find
effective conditions insuring that $(c_\alpha)$ are the moments of
a positive measure.

Or equivalently, in predual form, find effective criteria (in
terms of the coefficients) for a polynomial to be non-negative on
$\mathbb R^g$}.\\

We know from Tarski's principle that the positivity of a
polynomial is decidable. The term ``effective" above means to find
exact rational expressions in the coefficients which guarantee the
non-negativity of the polynomial.

We do not touch in this discussion a variety of other aspects of
the multivariate moment problem such as uniqueness criteria,
orthogonal polynomials, cubature formulas and the truncated
problem. See for instance \cite{Berg1, Berg2, CF,GV,KM}.

\subsection{Positivstellens\"atze on compact, semi-algebraic sets}
\label{sec:posSScompact}

Now we look at a very popular
 classes of Positivstellens\"atze. The hypotheses are
more  restrictive (by requiring bounded sets)
than the general one, but the conclusion
gives a simpler certificate of positivity.
The techniques of proof are  those used in
 the multivariate moment problem but
 measures with  compact semi-algebraic support
 allow  much more detail.

To state the theorems in this section requires  the notions of
preorder, $PO(F)$ and of quadratic module which we now give, but
the treatment of them in Section \ref{sec:PosSS} on the general
Positivstellensatz
  gives more properties and a different
context than done here.
Let
$F = \{ f_1,...,f_p \}$ denote  a set of real polynomials.
The {\it preordering}\index{preordering, PO}
generated by $F$ is
$$ PO(F) = \{ \sum_{\sigma \in \{ 0,1\}^r} s_\sigma f_1^{\sigma_1}
... f_r^{\sigma_r}; \ \ s_\sigma \in \Sigma^2 R[x] \}.$$
The {\it quadratic module generated by $F$} is defined to be:
$$ QM(F) = \sum_{f \in F \cup \{ 1 \}} f \Sigma^2 \mathbb R[x].$$
\index{quadratic module}

We start with a fundamental result of
Schm\"udgen,  proved in 1991
(\cite{Schmudgen1}), which
makes use in an innovative way of Stengle's
general Positivstellensatz.

\begin{thm}[Schm\"udgen]
\index{Schm\"udgen's Theorem}
\index{Positivstellensatz Schm\"udgen's}
Let
$F = \{ f_1,...,f_p \}$ be a set of real polynomials in $g$
variables, such that the non-negativity set $\POD{F}$ is
compact in $\mathbb R^g$. Then

a). A functional $L \in \mathbb R[x]'$ is representable by a
positive measure supported on $K$ if and only if
$$ L|_{PO(F)} \geq 0.$$

b). Every positive polynomial on $\POD{F}$ belongs to the
preorder $PO(F)$.
\end{thm}

Due to the compactness of the support, and Stone-Weierstrass
Theorem, the representing measure is unique. We will discuss later
the proof of b) in a similar context.

We
call the quadratic module $QM(F)$ {\it
archimedean}\index{archimedean} if there exists $C>0$ such that
$$ C - x_1^2 -...-x_g^2 \in QM(F).$$
This implies in particular that the semi-algebraic set
$\POD{F}$ is contained in
the ball centered at zero, of radius
$\sqrt{C}$. Also, from the convexity theory point of view, this
means that the convex cone $QM(F) \subset \mathbb R[x]$ contains
the constant function $1$ in its algebraic interior (see
\cite{Kothe} for the precise definition).
If
the set $\POD{F}$ is compact, then one can make the associated
quadratic module archimedean by adding to the defining set one
more term, of the form $C - x_1^2 -...-x_g^2$.

The key to Schm\"udgen's Theorem and to a few forthcoming results
in this survey is the following specialization of Proposition
\ref{general OT correspondence}.

\begin{lem}
\cite{P2}\label{joint spectrum}
Let $F$ be a finite set of
polynomials in $\mathbb R[x]$ with
associated quadratic module  $QM(F)$
having the archimedean property. There
exists a bijective correspondence between:\\

a). Commutative $g$-tuples $A$ of   bounded self-adjoint operators
with
cyclic vector $\xi$ and joint spectrum contained in $\POD{F}$;\\

b). Positive measures $\mu$ supported on $\POD{F}$;\\

c). Linear functionals $L \in \mathbb R[x]'$ satisfying
$L|_{QM(F)} \geq 0.$\\

The correspondence is constructive, given by the relations:
$$ L(p) = \langle p(A)\xi, \xi \rangle
= \int_{\POD{F}} p d\mu, \ \qquad \ p \in \mathbb R[x].$$
\end{lem}

\begin{proof}.  Only the implication $c) \Rightarrow a)$ needs an
argument. Assume c) holds and construct the Hilbert space $H$
associated to the functional $L$. Let $M=(M_{x_1},...,M_{x_g})$
denote the tuple of multiplication operators acting on $H$. Due to
the archimedean property,
$$ \langle (C - x_1^2 -...-x_g^2)p, p \rangle \geq 0, \ \ p \in
\mathbb C [x],$$ whence every $M_{x_k}$ is a bounded self-adjoint
operator. Moreover, the condition
$$ \langle f_j p, p \rangle \geq 0, \ \ p \in
\mathbb C [x],$$ assures that $f_j(M) \geq 0$, that is, by the
spectral mapping theorem, the joint spectrum of $M$ lies on
$\POD{F}$. Let $E_M$ be the joint spectral measure of $M$.
Then
$$ L(p) = \int_{\POD{F}} p(x) \langle E_M(dx)1, 1 \rangle,$$
and the proof is complete. \ \qed
\end{proof}

For terminology and general facts about spectral theory in a
commutative Banach algebra see \cite{Douglas}.

With this dictionary between positive linear functionals and
tuples of commuting operators with prescribed joint spectrum we
can improve Schm\"udgen's result.

\begin{thm}[\cite{P2}]
\label{quadratic module} \index{Positivstellensatz}
 Let $F$ be a finite
set of real polynomials
in $g$ variables, such that the associated quadratic module
$QM(F)$ is archimedean.

Then a  polynomial strictly positive on $\POD{F}$ belongs to $QM(F)$.
\end{thm}

\begin{proof}.
Assume by contradiction that $p$ is a positive
polynomial on $\POD{F}$ which does not belong to $QM(F)$. By
a refinement of Minkowski separation theorem due to Eidelheit and
Kakutani (see \cite{Kothe}), there exists a linear functional $L
\in \mathbb R[x]'$ such that $L(1)>0$ and:
$$ L(p) \leq 0 \leq L(q), \ \ q \in QM(F).$$
(Essential here is the fact that the constant function $1$ is in
the algebraic interior of the convex cone $QM(F)$).
Then Lemma
\ref{joint spectrum} provides a positive measure $\mu$ supported
on $\POD{F}$, with the property:
$$ L(p) = \int_{\POD{F}} p d\mu \leq 0.$$
The measure is non-trivial because
$$ L(1) = \mu(\POD{F}) >0,$$
and on the other hand $p>0$ on $\POD{F}$, a contradiction.
 \ \qed
\end{proof}

An algebraic proof of the latter theorem is due to Jacobi and
Prestel, see \cite{PD}.

\section{Complex variables}
\label{sec:complexVar}
The operator theoretic counterpart to
positive functionals described in the previous section becomes
more transparent in the case of complex variables. We present
below, closely following \cite{P3}, a series of generalizations of
Lemma \ref{joint spectrum} and Theorem \ref{quadratic module}. It
is at this point when von Neumann's inequality becomes relevant.

Throughout this section $z=(z_1,...,z_g)$ denote the complex
coordinates in $\mathbb C^g$. Then real coordinates of the
underlying space $\mathbb R^{2d}$ are denoted by $x =
(x_1,...,x_{2n})$, where $z_k = x_k + i x_{n+k}$. We will work as
before in the polynomial algebra $\mathbb C [x]= \mathbb
C[z,\overline{z}]$, and consider there the convex hulls of
non-negative polynomials:
\index{$\Sigma^2_h$}
$$ \Sigma^2 = {\rm co} \{ p^2;\ p \in \mathbb R [x]\},$$
and
$$ \Sigma^2_h = {\rm co} \{ |q|^2;\ q \in \mathbb C [z]\}.$$

It is easy to see that the cone of hermitian positive squares
$\Sigma^2_h$ is  a proper subset of $\Sigma^2$. Indeed, remark
that
$$ \frac{\partial}{\partial z_k}  \frac{\partial}{\partial
\overline{z}_k} |q|^2 \geq 0, \ \
q \in \mathbb C [z],$$ while the same Laplace operator has
negative values on $(y_k-x_k^2)^2$.

Let $F = \{ p_1,...,p_m\}$ be a finite subset of $\mathbb R [x]$
so that the basic semi-algebraic set
$$ \POD{F} =\{ x \in \mathbb R^{2d};\  p_1(x) \geq 0, ...,
p_m(x) \geq 0 \}$$
is compact. Let $p$ be a real polynomial which is positive on $
\POD{F}$.  We will indicate conditions which insure the
decompositions:
\begin{equation} p \in \Sigma^2_h + p_1 \Sigma^2_h + ...+ p_m \Sigma^2_h,
\end{equation}
or
\begin{equation} p \in \Sigma^2 + p_1 \Sigma^2_h + ...+ p_m \Sigma^2_h,
\end{equation}
or
\begin{equation} p \in QM(F) = \Sigma^2 + p_1 \Sigma^2 + ...+ p_m
\Sigma^2.
\end{equation}

The last one is covered by Theorem \ref{quadratic module}. The
other two require stronger assumptions on $p$, as we shall see
below.

We start by recalling an early, unrelated observation about
strictly positive hermitian polynomials \cite{Quillen}.

\index{Quillen's Theorem}
\begin{thm}[Quillen] If a bi-homogeneous polynomial $p \in \mathbb
C[z,\overline{z}]$ satisfies $p(z,\overline{z})>0$ for $z \neq 0$,
then there exists $N \in \mathbb N$ such that:
$$ |z|^{2N} p(z,\overline{z}) \in \Sigma^2_h.$$
\end{thm}
This result was rediscovered, and obtained by a different method,
by Catlin and d'Angelo \cite{CD1}. Their approach proved to be
geometric in its nature and very flexible, leading to a variety of
refinements of Quillen's theorem, see for instance \cite{CD2, dA,
dAV}. This line of research, not developed in the present survey,
culminates in completely removing the strict positivity
assumption. Specifically, the following characterization of
bi-homogeneous quotients of polynomials of $\Sigma^2_h$ was
recently discovered.

\index{Varolin's Theorem}
\begin{thm}[Varolin, \cite{Var}] Let the bi-homogeneous polynomial
 $p \in \mathbb C[z,\overline{z}]$ satisfy $p(z,\overline{z}) \geq 0$ for $z \in \mathbb C^g$.
 Write $p(z,\overline{z}) = \sum_{\alpha, \beta} p_{\alpha, \beta} z^\alpha \overline{z}^\beta$
 and, using the positive and negative spectral subspaces of the Hermitian matrix $(p_{\alpha, \beta})$,
 decompose $p = p_+ - p_-$, with $p_\pm \in \Sigma^2_h$.

 Then there are non-zero polynomials $s,t \in \Sigma^2_h$ with the property $$s p = t$$ if and only if
 there exists a positive constant $C$, such that $$p_+ + p_- \leq C(p_+ - p_-)$$ on
 $\mathbb C^g$.
 \end{thm}

Next we return to the compact semi-algebraic set $\POD{F}
\subset \mathbb C^g$ and the three levels of positivity
certificates described by the preceding convex cones.

We recall that a commutative $g$-tuple of linear bounded operators
$T$ acting on a Hilbert space $H$ is called {\it subnormal}
\index{subnormal $g$-tuple} if there exists a larger Hilbert space
$H \subset K$ and a commutative $g$-tuple of normal operators $N$
acting on $K$, so that every $N_j$ leaves $H$ invariant and
$N_j|_H = T_j, \ 1 \leq j \leq g$.
A commutative $g$ tuple $N$ of
normal operators $N_j = X_j +i X_{g +j}$ consists of $2g$ mutually
commuting self-adjoint operators $X_1,...,X_{2g}$.
Hence $N$
possesses a joint spectral measure $E_N$, supported on the joint
spectrum $\sigma(N) \subset \mathbb C^g$.

Assume from now on that the cone
$\Sigma^2_h + p_1 \Sigma^2_h + ... + p_m \Sigma^2_h$
is archimedean, that is, after a
normalization:
$$ 1 -|z_1|^2 -...-|z_g|^2 \in \Sigma^2_h + p_1 \Sigma^2_h + ...+ p_m
\Sigma^2_h.$$ Let $ L \in \mathbb C[x]'$ be a linear functional
satisfying
$$
L|_{\Sigma^2_h + p_1 \Sigma^2_h + ...+ p_m \Sigma^2_h} \geq
0.$$
Instead of constructing the completion of the whole ring of
real polynomials, we consider the same $L$-inner product, defined
only on complex polynomials $\mathbb C[z]$, in the variables $z$.
Let $H$ be the associated Hilbert space. The multiplication
operators $M_{z_j}$ act and commute on $H$. They are bounded due
to the above mentioned archimedean property:
$$ \| M_{z_j} q(z) \|^2 = L(|z_j q(z)|^2) \leq L(|q(z)|^2) = \| q
\|^2.$$ The only difference to the framework of the multivariate
Hamburger moment problem is that now $M_{z_j}$ are not necessarily
self-adjoint operators and the commutators $[M_{z_j},
M_{z_k}^\ast]$ may not vanish. The constant function vector $1$
remains cyclic, in the sense that the span of the vectors
$M_{z_1}^{\alpha_1} ... M_{z_g}^{\alpha_g} 1$ is the whole Hilbert
space $H$.

Let $M = (M_{z_1},...,M_{z_g})$ be the commutative n-tuple of
operators we have just constructed. For a polynomial
$p(z,\overline{z})$ we define after Colojoara and Foias, cf.
(\cite{AM}), the {\it hereditary functional
calculus}\index{hereditary functional calculus}
$\tilde{p}(M,M^\ast)$ by putting all adjoints $M_{z_k}^\ast$ in
the monomials of $p$ to the left of the powers of $M_{z_j}$'s. For
example,
$$ \widetilde{(|z_1|^2 z_2^2 \overline{z_3})}(M) = M_{z_1}^\ast M_{z_3}^\ast
M_{z_1} M_{z_2}^2.$$

We have thus established the first part of the following
dictionary.

\begin{prop}\label{hereditary calculus} Let $F = \{ p_1,...,p_m\}$ be
a finite set of real polynomials in $2g$
variables, such that
$$ 1 -|z_1|^2 -...-|z_g|^2 \in \Sigma^2_h + p_1 \Sigma^2_h + ...+ p_m
\Sigma^2_h.$$

a). There exists a bijective correspondence between functionals $L
\in \mathbb C[x]'$ which are non-negative on $\Sigma^2_h + p_1
\Sigma^2_h + ...+ p_m \Sigma^2_h$ and commutative $g$-tuples of
linear Hilbert space operators $T$, with a cyclic vector $\xi$,
subject to the conditions
$$ \tilde{p_j}(T,T^\ast) \geq 0, \ 1 \leq j \leq m.$$

b). If, in addition $L$ is non-negative on the larger cone
$\Sigma^2 + p_1 \Sigma^2_h + ...+ p_m \Sigma^2_h,$ then, and only
then, $T$ is also subnormal.

In both cases, the correspondence is given by
$$ L(p(z,\overline{z})) = \langle \tilde{p}(T,T^\ast)\xi, \xi \rangle.$$
\end{prop}

\begin{proof}.\
For the proof of part b) remark that the Hilbert space
completion $H$ of the ring of complex polynomials with respect to the
inner product $\langle p,q\rangle = L(p \overline{q}), \ p,q \in
\mathbb C[z]$ carries the bounded, commutative $g$-tuple $M$ of
multiplication operators with the variables $z_1,...,z_g$.
The positivity of the functional $L$ on $\Sigma^2$ is equivalent to the
multivariate analog of the Halmos-Bram subnormality condition applied to $M$.
See for details \cite{Demanze, P2}.

Conversely, if $T$ is a cyclic subnormal tuple of operators with
normal extension $N$, and $p(z,\overline{z})$ is a polynomial,
then
$$
\langle {|\tilde p|}^2(T,T^\ast) \xi, \xi \rangle
= \langle |p|^2(N,N^\ast) \xi, \xi
\rangle = \| p(N,N^\ast)\xi\|^2 \geq 0.$$
This uses
the very definition of the hereditary calculus,
for example,
$$
\langle
T_2^{*^4} T_1^{*^6} T_1^{6}
T_2^4
\xi, \; \xi \rangle
=
\langle
T_1^6 T_2^4 \xi, \; T_1^{6} T_2^{4} \xi \rangle
=
\langle
T_1^6 N_2^4 \xi, \; T_1^{6} N_2^{4} \xi \rangle
$$
$$
\langle
N_1^6 N_2^4 \xi, \; N_1^{6} N_2^{4} \xi \rangle
=
\| N_1^{6} N_2^{4} \xi \|.
$$
 \ \qed
\end{proof}

The following translation of the proposition shows that the class
of all commutative tuples of operators serves as a better
``spectrum" for the polynomial algebra in the variables
$(z,\overline{z})$.

\begin{cor} Let $F$ be as in the Proposition and let
$p(z,\overline{z})$ be a polynomial. If $\tilde{p}(T,T^\ast) > 0$
for every commutative $g$-tuple of linear Hilbert space operators
$T$, satisfying $ \tilde{p_j}(T,T^\ast) \geq 0, \ 1 \leq j \leq
m,$ then $p$ belongs to $\Sigma^2_h + p_1 \Sigma^2_h + ...+ p_m
\Sigma^2_h.$
\end{cor}

  The cyclic vector condition is not relevant for this statement.

\begin{proof}.\  The proof follows from now a known pattern. Assume by
contradiction that $p \notin \Sigma^2_h + p_1 \Sigma^2_h + ...+
p_m \Sigma^2_h.$ By Minkowski-Eidelheit-Kakutani separation
theorem, there exists a linear functional $L  \in \mathbb C[x]'$
satisfying the conditions of Proposition \ref{hereditary calculus}
and $L(p) \leq 0 < L(1)$. Then the commutative $g$-tuple $M$
associated to the inner-product space defined by $L$ satisfies
$$ \langle p(M,M^\ast)1, 1\rangle = L(p) \leq 0,$$
a contradiction.  \ \qed
\end{proof}

Even the simple case of the unit ball or unit polydisk in $\mathbb
C^g$ is interesting from this perspective. Assume first that
$n=1$. According to von-Neumann's inequality,
$$ \| p(T) \| \leq 1,$$
whenever $T$ is a contraction and $\sup_{z \in \mathbb D} |p(z)|
\leq 1.$ Thus, in view of the above proposition, for every
polynomial $p(z)$ and constant $M > \sup_{z \in \mathbb D}
|p(z)|$, we have
$$ M^2 - |p(z)|^2 \in \Sigma^2_h + (1-|z|^2)\Sigma^2_h.$$
Needless to say that this statement is equivalent to
von-Neumann's inequality.

In complex dimension two, a celebrated theorem of Ando (see
\cite{AM, CW}) asserts that, for every pair of commuting \index{Ando's Theorem}
contractions $(T_1,T_2)$ and every polynomial $p(z_1,z_2)$ one has
a von-Neumann type inequality:
$$ \| p(T_1,T_2)\| \leq \| p \|_{\infty, \mathbb D^2}.$$
And a not less celebrated example of Varopoulos (see again
\cite{AM}) shows that it is no more the case in dimension $n=3$
and higher. Specifically, according to our corollary, for every
polynomial $p(z)$ and $\epsilon>0$, we have
$$  (\| p \|_{\infty, \mathbb D^2} + \epsilon)^2 - |p(z)|^2 \in
\Sigma^2_h + (1-|z_1|^2)\Sigma^2_h + (1-|z_2|^2)\Sigma^2_h,$$ but
the statement is not true (for the unit polydisk) in higher
dimensions.\\

{\bf Open problem.} {\it It would be interesting to find an
algebraic explanation of Ando's Theorem via the above equivalent
sums of squares decomposition.}\\

  On the other hand, by exploiting part b) of Proposition
  \ref{hereditary calculus} one can prove a sharper weighted sums of
  squares decomposition on complex analytic polyhedra.

\begin{thm} Let $S = \{ z \in \mathbb C^g; |p_j(z)| \leq 1, 1 \leq j \leq m\}$
be a compact semi-algebraic set, where $p_j$ are complex
polynomials. Assume that the convex cone
$$C = \Sigma^2 +
(1-|p_j|^2) \Sigma^2_h+... + (1-|p_m|^2) \Sigma^2_h$$
  is archimedean. Then every real polynomial $p$ which is strictly
positive on $S$ belongs to $C$.
\end{thm}

For a proof see \cite{P2,P3}. The article \cite{BGM05} contains a
similar approach to bounded analytic multipliers on
 the Hardy
space of the bi-disk.  The recent note \cite{P3} contains a few
other sums of squares translations of some recently proved
inequalities in operator theory.

Finally, we reproduce from \cite{HMP3} the following general
Nichtnegativstellensatz over the complex affine space. Note that
when evaluating on tuples of commutative matrices, we do not have
to impose the strict positivity of the polynomial to be
decomposed.

\begin{thm}\index{Nichtnegativstellensatz on $\mathbb C^g$} Let $p(z, \overline{z})$ be a real valued polynomial,
where $z \in \mathbb C^g$. Then there are polynomials $q_i \in
\mathbb C[z], \ 1 \leq i \leq k,$ with the property
$$p(z, \overline{z}) = \sum_{i=1}^k |q_i(z)|^2,$$
if and only if, for all tuples of commuting matrices $X =
(X_1,...,X_g) \in M_d(\mathbb C), \ \ d \geq 1,$ we have
$$ p(X,X^\ast) \geq 0.$$
\end{thm}

The proof follows the general scheme outlined in this section and
we omit it. See \cite{HMP3} for full details.

\section{Real algebra and mathematical logic}
\label{sec:logic} Keeping in mind the main theme of our essay
(sums of squares decompositions), we briefly recall below, without
aiming at completeness, some classical facts of real algebra and
mathematical logic. We follow an approximate chronological order.
For a more comprehensive, recent and very authoritative survey of
real algebra and real algebraic geometry aspects of sums of
squares we refer to Scheiderer \cite{Scheiderer2}.

\subsection{Minkowski and Hilbert}
\label{sec:MinkHilb} In the same time to, and even before, the
analysis aspects of sums of squares decompositions we have
discussed have been discovered, similar questions have appeared in
number theory and algebra. Lagrange's famous theorem (that every
positive integer can be written as a sum of squares of four
integers) was the origin of many beautiful studies, see Chapter XX
of Hardy and Wright's monograph \cite{HW}.

According to Hilbert \cite{Hilbert3}, after the teenager Minkowski
\index{Minkowski} won the 1882 ``Grand Prix" of the French Academy
of Sciences on a theme related to Lagrange's four squares theorem,
he has started working in K\"onisberg on his thesis devoted to
quadratic forms of a higher number of variables. It was in
Minkowski's inaugural dissertation, with Hilbert as opponent, that
he remarked that {\it ``it is not probable that every positive
form can be represented as a sum of squares"} \cite{Minkowski}.

The opponent (Hilbert) produced the first (non-explicit) example,
see \cite{Hilbert1}. \index{Hilbert's example} His idea is the
following. Consider nine points $a_1,...,a_9$ in $\mathbb R^2$, as
the base of a pencil of cubics (that is, a family of curves
obtained via a linear combination of their third degree
 defining equations), so
that every cubic polynomial vanishing at the first eight points
$a_1, ...,a_8$ will automatically vanish at $a_9$. By a rather
involved geometric argument, one can prove the existence of a
polynomial $p(x,y)$, of degree six, which is non-negative on
$\mathbb R^2$, vanishes at $a_1,...,a_8$ and satisfies $p(a_9)>0$.
Then clearly $p$ cannot be written as a sum of squares of
polynomials:
$$ p = q_1^2 + ...+ q_N^2,$$
because every $q_i$ would have degree at most three and therefore
would be null at $a_9$, too. Hilbert's argument is reproduced at
the end of Chapter II of Gelfand and Vilenkin's monograph
\cite{GV}.

The first explicit example based on Hilbert idea was constructed
by Robinson in 1969, see \cite{Reznick1} for details. Robinson's
homogenized polynomial is \index{Robinson's example}:
$$ P(x,y,z) = x^6 + y^6 + z^6 - $$ $$-(x^4 y^2 + x^2 y^4 + x^4 z^2 + x^2
z^4 + y^4z^2 + y^2 z^4) + 3 x^2 y^2 z^2.$$ About the same time
(some six dozen years after Hilbert's article) Motzkin has
produced a very simple polynomial (again shown in homogenized
form)\index{Motzkin's example}:
$$ Q(x,y,z) = z^6 + x^2 y^2 (x^2 + y^2 - z^2).$$
The reader will find easily why $Q$ is non-negative, but not a sum
of squares. More examples of non-negative polynomials which are
not sums of squares were discovered by Choi, Lam and Reznick, and
separately Schm\"udgen. We refer to Reznick's monograph
\cite{Reznick1} for more details and for an elegant geometric
duality method (based on the so called Fisher inner product)
adapted to the analysis of the convex cones of such polynomials.

One of Hilbert's celebrated 1900 problems in mathematics was about
the structure of positive polynomials and the logical implications
of the existence of a constructive way of testing
positivity. Here are his words:\\

{\sl Problem 17. Expression of definite forms by squares:}\\

{\it ``...the question arises whether every definite form may not
be expressed as a quotient of sums of squares of forms...

... it is desirable, for certain questions as to the possibility
of certain geometrical constructions, to know whether the
coefficients of the forms to be used in the expression may always
be taken from the realm of rationality given by the coefficients
of the form represented."} \cite{Hilbert3}.\index{Hilbert's 17-th
Problem}\\

Hilbert's intuition proved to be correct on both conjectures
raised by his question. His query was the origin of a series of
remarkable results in algebra and logic, see \cite{BCR, PD,
Scheiderer2}.

\subsection{Real fields}
Hilbert's 17-th problem was solved in the affirmative by E. Artin
\cite{Artin} in 1927, as an application of the theory of real
fields he has developed with Schreier \cite{AS}. For history and
self-contained introductions to real algebra, and complete proofs
of Artin's Theorem, we refer to either one of the following
monographs \cite{BCR, Jacobson, PD}. We merely sketch below the
main ideas of Artin-Schreier theory (as exposed in Jacobson's
algebra book \cite{Jacobson}), to serve as a comparison basis for
the computations we will develop in later sections in the
framework of star algebras.

An {\it ordered field}\index{ordered field} is a characteristic
zero field $F$ with an {\it ordering}\index{ordering}, that is a
prescribed subset $P$ of positive elements, satisfying:
$$ P + P \subset P, \ \ \ P \cdot P \subset P, \ \ \ F = P \cup \{ 0 \}
\cup \{ -P\}.$$ Since $a^2 \in P$ for all $a \neq 0$, if $\sum_i
a_i^2 = 0$, then every $a_i = 0$. Or equivalently, $-1$ cannot be
written as a sum of squares in $F$. By a theorem of Artin and
Schreier, every field with the latter property can be ordered. An
ordered field $R$ is {\it real closed}\index{real closed field} if
every positive element has a square root. In this case, exactly as
in the case of real numbers, the extension $R(\sqrt{-1})$ is
algebraically closed.

A central result in Artin-Schreier theory is the existence and
uniqueness of the real closure $R$ of an ordered field $F$: that
is the extension $F \subset R$ is algebraic and $x \in F$ is
positive in $R$ if and only if is positive in $F$. Interestingly
enough, the proof of this fact uses Sturm's algorithm for the
determination of the number of roots of a polynomial with real
coefficients (or more generally with coefficients in a real closed
field).\\

{\bf Sturm's algorithm.} \index{Sturm's Algorithm} Let $R$ be a
real closed field and let
$$p(x) = a_0 + a_1 x + ...+  x^d \in R[x]$$ be a polynomial.
Let $$ C = 1 + |a_{d-1}| +...+|a_1|+|a_0|.$$
  Define the sequence of polynomials:
$$ p_0 = p, \ \ \ p_1 = p' \ \ ({\rm the \ derivative}),$$
$$ p_{j+1} = p_j q_j - p_{j-1}, \ \ \deg p_{j+1} < \deg p_j.$$
Then for a large $n$, $p_n$ = 0. Sturm's Theorem asserts that the
number of roots of $p$ in $R$ is $N(-C) -N(C)$, where $N(a)$ is
the number of sign changes in the sequence $p_0(a), p_1(a),...,
p_n(a)$.

As an application of the existence of the real closure of an
ordered field, one can prove that an element $x \in F$ is a sum of
squares if and only if it is positive in every order on $F$. Or,
equivalently, if $y \in F$ is not a sum of squares, then there
exists an order on $F$ with respect to which $y<0$.

\begin{thm}[Artin]\index{Artin's Theorem} Let $F$ be a subfield of
$\mathbb R$ which has a unique ordering and let $f$ be a rational
function with coefficients in $F$. If $f(a_1,...,a_g)\geq 0$ for
all $(a_1,...,a_g) \in F^g$ for which $f$ is defined, then $f$ is
a sum of squares of rational functions with coefficients in $F$.
\end{thm}

  The idea of the proof is to admit by contradiction that $f$ is
not a sum of squares, hence it does not belong to an ordering of
the field of rational functions $K = F(x_1,...,x_g)$ in $g$
variables. By completing $K$ to a real closed field $R$, one finds
``ideal" points $b_1,...,b_g \in R$, so that $f(b_1,...,b_g)<0$.
By Sturm's counting theorem one shows then that there are points
$a_1,...,a_g \in F$ with the property $f(a_1,...,a_g)<0$. For
details see for instance \S 11.4 in \cite{Jacobson}.\\

Examples of fields with a unique ordering are $\mathbb Q$ and
$\mathbb R$. Artin's Theorem prompts a series of natural
questions, as for instance : how many squares are necessary, is
there a universal denominator in the decomposition of $f$ as a sum
of squares, are there degree bounds? All these problems were
thoroughly studied during the last decades, \cite{BCR, PD}.

It was Tarski who in the late 1920-ies put Sturm's Algorithm into
a very general and surprising statement. His work had however an
unusually long gestation and has remained unknown to the working
mathematician until mid XX-th Century. His articles are available
now from a variety of parallel sources, see for instance
\cite{Tarski} and the historical notes in \cite{PD}. His main
thesis is contained in the following principle, cited from his
original 1948 RAND publication \cite{Tarski}:\\

  {\bf Tarski's elimination
theory for real closed fields.}\index{Tarski's Principle} ``{\it
To any formula \\
$\phi(x_1, ...,x_g)$ in the vocabulary $\{ 0,1,+,\cdot,<\}$ and
with variables in a real closed field, one can effectively
associate two objects:

(i) a quantifier free formula $\overline{\phi}(x_1, ...,x_g)$ in
the same vocabulary, and

(ii) a proof of the equivalence $\phi \equiv \overline{\phi}$ that
uses only the axioms of real closed fields.}"\\

He aimed this theorem at the completeness of the logical system of
elementary algebra and geometry, very much in the line of
Hilbert's programme in the foundations of mathematics. As a
consequence one obtains the transfer principle alluded above, in
the proof of Artin's Theorem: let $R_1 \subset R_2$ be real closed
fields. A system of polynomial inequalities and equalities with
coefficients in $R_1$ has a solution in $R_2$ if and only if it
has a solution in $R_1$.

Let $R$ be a real closed field. We recall that a {\it
semi-algebraic set}\index{semi-algebraic set} in $R^g$ is a finite
union of finite intersections of sets of the form
$$ \{ x \in R^g; \ p(x)=0\}, \ \  \{ x \in R^g; \ p(x)> 0\},\ \
{\rm or} \ \  \{ x \in R^g; \ p(x) \geq 0\}.$$

A self-contained account of Tarski's theorem can be found in
\cite{Seidenberg}. See also \cite{BCR, Jacobson, PD}. In practice,
the most useful form of Tarski's result is the following theorem.

\begin{thm}[Tarski-Seidenberg] If $R$ is a real closed field and
$S$ is a semi-algebraic set in $R^g \times R^m$, then the
projection of $S$ onto $R^g$ is also semi-algebraic.
\end{thm}

\noindent For a self contained proof of the above theorem see the
Appendix in  \cite{Hormander}.

Applications of Tarski's principle came late, but were
spectacular. We only mention for illustration one of them:\\

{\bf H\"ormander's inequality.} (1955) \index{ H\"ormander's
inequality} : {\it For every
polynomial\\
$f(x_1,...,x_g) \in \mathbb R[x_1,...,x_g]$ there are positive
constants $c,r$ such that
$$ |f(x)| \geq c \ {\rm dist}(x, V(f))^r,\ \  x \in \mathbb R^g, |x|
\leq 1.$$}\\
  Above $V(f)$ stands for the real zero set of $f$. The
inequality was generalized to real analytic functions by
Lojasiewicz in 1964, and served as the origin of fundamental
discoveries in modern analytic geometry and the theory of partial
differential operators, see for instance \cite{BCR, Hormander}.

\subsection{The general Positivstellensatz}
\label{sec:PosSS} A great jewel of real algebraic geometry, which
is now causing excitement in applications, is the
Positivstellen-s\"atze and the real Nullstellensatz (which it
contains). This section states these theorems and can be read
independently of earlier parts of this paper.

The  Positivstellen-s\"atze
 lives
in a polynomial ring with coefficients in a real closed field (as
opposed to the complex numbers) and
 were discovered only in the 1960-ies
 (see \cite{Dubois,Krivine})
 before being rediscovered and refined by Stengle in 1974 (\cite{Stengle}).
The statement of the Nullstellensatz departs from Hilbert's
Nullstellensatz over an algebraically closed field, by imposing an
additional sum of squares term in the characterization of a
radical ideal, as we shall see below. It is interesting to remark
that Stengle's article makes specific references, as origins or
motivations of his investigation, to works in mathematical logic
and mathematical programming.

In order to state Stengle's Positivstellen-s\"atze we need first a
few definitions and conventions. Let $R$ be a real closed field
(many readers  will be happy to think of $R$ as the real numbers)
and denote $x=(x_1,...,x_g) \in R^g$ and also regard $x$ as a
$g$-tuple of commuting indeterminates. Let $\Sigma^2 A$
\index{$\Sigma^2 A$} denote the set of all sums of squares in the
algebra $A$ \index{sums of squares}.

Let $S \subset R[x]$ be a subset, and write
$$ 
\POD{S} = \{ x \in R^g; \ p(x) \geq 0, \ \forall p \in S\},
$$
for the positivity set of of the functions $S$. \index{$\POD{S}$
 non-negativity set of $S$}
If $S = \{ p_1,...,p_r\}$ is finite then $\POD{S}$ is a {\it basic
closed semi-algebraic set} \index{basic closed semi-algebraic set
$\mathcal D (S)$}. The {\it preordering}\index{preordering, PO}
generated by $S$ is
$$ PO(S) = \{ \sum_{\sigma \in \{ 0,1\}^r} s_\sigma p_1^{\sigma_1}
... p_r^{\sigma_r}; \ \ s_\sigma \in \Sigma^2 R[x] \}.$$ The {\it
quadratic module}\index{quadratic module, QM} generated by $S$ is
$$ QM(S) = \Sigma^2 R[x] + p_1 \Sigma^2 R[x] + \ldots + p_r \Sigma^2
R[x].$$

Note that a preordering satisfies conditions similar to an
ordering in a field:
$$  PO(S) + PO(S) \subset PO(S), \ \ PO(S) \cdot PO(S) \subset
PO(S), \ \ \Sigma^2 R[x] \subset PO(S),$$ while the quadratic
module fails to be closed under multiplication, but still
satisfies:
$$ QM(S) + QM(S) \subset QM(S), \ \ \Sigma^2 R[x] \subset QM(S), \ \
\Sigma^2 R[x] \cdot QM(S) \subset QM(S),$$ and clearly,
$$ QM(S) \subset PO(S).$$\\

\begin{thm}[Stengle]
\label{thm:stengle} \index{Stengle's Theorem}
\index{Positivstellensatz, Stengle} Let $R$ be a real closed field
and let $p_1,...,p_r \in R[x_1,...,x_g].$ Let $\cS = \mathcal
D(p_1,...,p_r)$ and let $T = PO(p_1,...,p_r)$ be the preorder
generated by $p_i$. Let $f \in R[x_1,...,x_g]$. Then \\

(a). \index{Positivstellensatz} $f>0$ on $\cS$ if and only if
there are $s,t \in T$ satisfying
$sf = 1+t$;\\

(b). \index{Nichtnegativstellensatz}$f \geq 0$ on $\cS$ if and
only if there are $s,t \in T$ and an
integer $N \geq 0$, with the property $sf = f^{2N} + t$;\\

(c). $f=0$ on $\cS$ if and only if there exists an integer $N \geq
0$ with the property $-f^{2N} \in T$.\\
\end{thm}

We derive a few particular consequences. For instance, the
following real Nullstellensatz is contained in the preceding
result: assume that $p,q \in R[x]$ and
$$ (p(x) = 0) \ \Rightarrow (q(x) = 0).$$
Then point (c) applies to $f=q$ and $p_1 = -p^2$. We infer: there
exists $N \geq 0$ such that
$$ -q^{2N} = s_1 - p^2 s_2, \ \ s_1, s_2 \in \Sigma^2 R[x],$$
therefore:
$$ q^{2N} + s \in (p), \ \ s \in \Sigma^2 R[x],$$
where $(p)$ denotes the ideal generated by $p$. Obviously, if the
latter condition holds, then $q$ vanishes on the zero set of $p$.
Thus, we have proved: \index{Nullstellensatz}
$$ [(p(x) = 0) \ \Rightarrow (q(x) = 0)]
\Leftrightarrow [\exists (N \geq 0, \ s \in \Sigma^2 R[x]): \
q^{2N} + s \in (p)].$$ Variants of this are pleasurable, and we
suggest as an exercise the reader repeat the above but take $f=q$
and $p_1= p, \ p_2=-p$.

As another example, assume that $p(x) \geq 0$ for all $x \in R^g$.
Then the theorem applies to $f=p$ and $p_1 = 1$ and we obtain:
there exists an integer $N \geq 0$ and elements $ s_1, s_2 \in
\Sigma^2 R[x]$, such that
$$ s_1 p = p^{2N} + s_2.$$ In particular,
$$ s_1^2 p \in \Sigma^2 R[x],$$
which is exactly the conclusion of Artin's Theorem.

  The concepts of
(pre)ordering and quadratic module can be defined for an arbitrary
commutative ring with unit; these, together with the important
construct of the real spectrum provide the natural framework for
developing modern real algebra and real algebraic geometry. For
the general versions of the S\"atze outlined in this section the
reader can consult as a guide \cite{Scheiderer2}, and for complete
details and ramifications, the monographs \cite{BCR,PD}.

\section{Applications of semi-algebraic geometry}
\label{sec:applSAG}
The prospect of applying  semi-algebraic geometry to a variety of
areas is the cause of excitement in many communities; and we list
a few of them here.

\subsection{Global optimization of polynomials}

An exciting turn in the unfolding of real algebraic geometry are
applications to optimization. To be consistent with the
non-commutative setting of the subsequent sections we denote below
by $x \in \mathbb R^g$ a generic point in Euclidean space, and in
the same time the $g$-tuple of indeterminates in the polynomial
algebra.

\subsubsection{Minimizing a Polynomial on $\RR^g$}
A classical question is: given a polynomial $q \in \mathbb R[x]$,
find
    $$ \min_{x \in \RR^g} q(x)$$
and the minimizer $x^{opt}$. The goal is to obtain a numerical
solution to this problem and it is daunting even in a modest
dimension such as $g = 15$. Finding a local optimum  is
numerically ``easy" using the many available variations of
gradient descent and Newton's method. However, polynomials are
notorious for having many many local minima.

A naive approach is to grid $\RR^g$, lets say with 64 grid points
per dimension (a fairly course grid), and compare values of $q$ on
this grid. This requires $64^{15}$
$ \sim 10^9 10^7$ function evaluations or something like  10,000
hours to compute. Such prohibitive requirements occur in many high
dimensional spaces and go under the heading of the ``curse of
dimensionality".

The success of sums of squares and Positivstellens\"atze methods
rides on the heels of semi-definite programming, a subject which
effectively goes back a decade and a half ago, and which
effectively allows numerical computation of a sum of squares
decomposition of a given polynomial $q$. The cost of the
computation is determined by the number of terms of the polynomial
$q$ and is less effected by the number $g$ of variables and the
degree of $q$. To be more specific, this approach to optimization
consists of starting with a number $q^{**}$ and numerically solve
$$
q - q^{**} = s,$$ for  $s \in \Sigma^2$. If this is possible,
lower $q^{**}$ according to some algorithm and try again. If not,
raise $q^{**}$ and try again. Hopefully, one obtains $q^{*o}$ at
the transition (between being possible to write $q-q^{**}$ as a
sums of squares and not) and obtains
$$
q - q^{*o} \in \Sigma^2
$$
and conclude that this is an optimum. This method was proposed
first by Shor \cite{Shor} and subsequently refined by Lasserre
\cite{Lasserre1} and by Parrilo \cite{parThesis}.

Parrilo and Sturmfels \cite{PS} reported experiments with a
special class of 10,000 polynomials for which the true global
minimum could be computed explicitly. They found in all cases that
$q^{*o}$ determined by sums of squares optimization equals the
true minimum.

Theoretical evidence supporting this direction is \index{SOS
approximation of polynomials}
  the following observation, see \cite{BCR} \S 9.
\begin{thm}
\label{thm:bdedOpt} Given a polynomial $q \in \mathbb R[x]$, the
following are equivalent:

(1) $q \geq 0$ on the cube $[-1, 1]^g$.

(2) For all $\eps >0$, there is $s \in \Sigma^2$ such that
$$ \| q - s \|_{L^1( [-1, 1]^g )}< \eps.$$
\end{thm}

A refinement of this result was recently obtained by Lasserre and
Netzer \cite{Lasserre-Netzer}. Namely, the two authors prove that
an additive, small perturbations with a fixed polynomial, produces
a sum of squares which is close to the original polynomial in the
$L^1$ norm of the coefficients. We reproduce, without proofs,
their main result.

\begin{thm} \cite{Lasserre-Netzer} \index{SOS approximations via high degree
perturbations} Let $p \in \mathbb R [x_1,...,x_g]$ be a polynomial
of degree $d$, and let
$$ \Theta_r  = 1 + x_1^{2r} + ...+ x_g^{2r},$$
where $r \geq d/2$ is fixed. Define
$$ \epsilon_r^\ast = \min_L \{ L(p);\  L \in \mathbb R_{2r}
[x_1,...,x_g]',  \ L(\Theta_r) \leq 1, \ L|_\Sigma^2 \geq 0 \}.$$

Then $ \epsilon_r^\ast \leq 0$ and the minimum is attained. The
polynomial
$$ p_{\epsilon,  r} = p + \epsilon \Theta_r$$
is a sum of squares if and only if $\epsilon \geq -
\epsilon_r^\ast.$

Moreover, if the polynomial $p$ is non-negative on the unit cube
$[-1,1]^g$, then $\lim_{r \rightarrow \infty}  \epsilon_r^\ast =
0$.

\end{thm}

Variations of the above theorem, with supports on semi-algebraic
sets, relevant examples and an analysis of the degree bounds are
contained in the same article \cite{Lasserre-Netzer}.

For quite a few years by now, Lasserre has emphasized the
tantamount importance of such perturbation results for
optimization using sums of squares (henceforth abbreviated $SOS$)
methods, see \cite{Lasserre1}, in that it suggests that
determining if a given $p$ is nonnegative on a bounded region by
computing a sums of squares has a good probability of being
effective.

We shall not prove the stated perturbation results, but remark
that a free algebra version of them holds, \cite{KSprept}.

In the opposite pessimistic direction there are the precise
computations of Choi-Lam-Reznick (see \cite{Reznick1}) and a
recent result due to Bleckermann \cite{Bleckermann}.

As a backup to the above optimization scheme, if a $ q - q^{*o}
\in \Sigma^2 $ fails to be a sum of squares, then one can pick a
positive integer $m$ and attempt to solve
$$
( 1 + |x|^2)^m (q - q^{*o}) \in \Sigma^2.
$$
Reznick's Theorem \cite{Reznick2} tells us that for
 some $m$ this solves
the optimization problem exactly. Engineers call using the term
with some non zero $m$ ``relaxing the problem", but these days
they call most modifications of almost anything a ``relaxation"
\index{relaxation}.

\subsubsection{Constrained optimization}
Now we give Jean Lasserre's  interpretation of Theorem
\ref{quadratic module}. Let $\cP$ denote a collection of polynomials.
The standard constrained optimization
problem for polynomials is:\\

\centerline{\it  minimize $q(x)$ subject to
$x \in
\POD{\cP} := \{x \in \mathbb R^g;  p(x) \geq 0,
\ p \in \cP \}$.}
\bigskip

Denote the minimum value of $q$ by $q^{ {\tiny{opt}}}$. We
describe the idea when $\cP$ contains but one polynomial $p$.
Assume $\nabla p(x)$ does not vanish for $x \in \partial \mathcal
D_p$.

The standard first order necessary conditions for $x^{
{\tiny{opt}}} \in \partial \POD{\cP} $ to be a local solution to this problem is
$$ \nabla q(x^{ {\tiny{opt}}}) = \lambda \nabla p(x^{ {\tiny{opt}}})$$
with $ \lambda >0$. We emphasize, this is a local condition and
$\lambda $ is called the Lagrange multiplier.

Now we turn to analyzing the global optimum. Suppose that $q$ can
be expressed in the form:
$$ q -q^{**} = s_1 + s_2 p,\ \ s_{1,2} \in \Sigma^2,$$
which implies $q(x) \geq q^{**}$ for all $x \in \POD{p}$. So
$q^{**}$ is a lower bound. This is a stronger form of the
Positivstellensatz than is always true.
Then this optimistic
statement can be interpreted as a global optimality condition when
  $q^{**}=q^{ {\tiny{opt}}}$.
Also it implies the classical Lagrange multiplier linearized condition,
as we now see.
At the {\it global} minimum $x^{
{\tiny{opt}}}$ we have
$$0= q(x^{ {\tiny{opt}}}) -q^{ {\tiny{opt}}}
= s_1 (x^{ {\tiny{opt}}}) + s_2 (x^{ {\tiny{opt}}}) p(x^{
{\tiny{opt}}})$$ which implies $0 = s_1( x^{ {\tiny{opt}}})$ and,
since $s_1$ is a sum of squares, we get
$\nabla s_1 (x^{ {\tiny{opt}}})=0.$
Also  $ s_2 (x^{ {\tiny{opt}}}=0$,
$\nabla s_2 (x^{ {\tiny{opt}}})=0$ whenever $ p(x^{ {\tiny{opt}}})\neq 0$.
Calculate
$$
\nabla q \ = \ \nabla s_1 +  p \nabla s_2  + s_2 \nabla p.
$$
If $ p(x^{ {\tiny{opt}}})=0$,  we get
$$
  \nabla q(x^{ {\tiny{opt}}})  \
  = \  s_2 (x^{ {\tiny{opt}}}) \nabla p(x^{ {\tiny{opt}}})
$$
  and if $ p(x^{ {\tiny{opt}}}) \neq 0$
  we get  $\nabla q(x^{ {\tiny{opt}}})=0$,
  the classic condition for an optimum in the interior.
Set $\lambda = s_2(x^{ {\tiny{opt}}})$ to get $\lambda \nabla
p(x^{ {\tiny{opt}}}) = \nabla q(x^{ {\tiny{opt}}})$
  the classic Lagrange multiplier condition
  as a (weak) consequence of the Positivstellensatz.

The reference for this and more general (finitely many $p_j$ in
terms of the classical Kuhn-Tucker optimality conditions) is
\cite{Lasserre1} Proposition 5.1.

\bigskip

Also regarding constrained optimization
we mention that, at the technical level,
the method of moments has re-entered into polynomial
optimization.
Quite specifically, Lasserre and followers are
relaxing the original problem
$$ \min_{x \in \mathcal D} q(x) $$
as
$$ \min_{\mu} \int_{\mathcal D} q d\mu,$$
where the minimum is taken over all probability measures supported
on $\mathcal D$. They prove that it is a great advantage to work
in the space of moments (as free coordinates), see \cite{Henrion2,
Lasserre1, Lasserre2}.

\subsection{Primal-dual optimality conditions}

In this section we explain in more detail Lasserre's point of view
\cite{Lasserre0, Lasserre1} of linearizing polynomial optimization
via sums of squares decompositions and via moment data seen as
independent variables.

Specifically, we start with a polynomial $f \in {\mathbb R}[x]$
and seek values of the scalar $\lambda$ for which $f - \lambda \in
\Sigma^2$. To this aim we consider a variable linear functional $L
\in {\mathbb R}[x]'$ and denote the corresponding moments
$$ y_\alpha = L(x^\alpha), \ \ |\alpha| \leq 2d.$$
The integer $d$ is fixed throughout the whole section and will not
explicitly appear in all coming formulas. We denote $y =
(y_\alpha)_{|\alpha| \leq 2d}$ and consider the associated Hankel
matrix
$$ M_y = (y_{\alpha + \beta})_{|\alpha|, |\beta| \leq d}.$$
In all these considerations it is important to fix an ordering
(such as the graded lexicographic one) on the multi-indices $\alpha$. Let
$$ V(x) = (1,  x_1,  \cdots,  x_g,  x_1^2,  \cdots )$$
be the "tautological" vector consisting of all monomials of degree
less than or equal to $d$. Let $ \alpha_0 = (0, 0, \cdots, 0)$, so
$x^{\alpha_0} = 1$.

The matrix valued polynomial
$$ \sum_\alpha B_\alpha x^\alpha = V(x) \cdot V(x)^T,$$
produces a sequence of matrix coefficients $B_\alpha$ which carry
a  Hankel type structure
and whose entries are either $0$ or $1$. Write
$$ f = \sum_\alpha f_\alpha x^\alpha.$$

The next lemma equates "minimization" of a polynomial
  via a sum of squares to
a matrix problem.
\begin{lem}
\label{lem:OmegaSoS} The degree $d$ polynomial
$f-\lambda$  is a sum of squares with $r$
squares of polynomials each having degree $\leq n$ if and only if
there exists a positive semi-definite matrix
$\Omega \in \Rnn$
of rank $r$ such that
\beq
\label{eq:Omegaf} {\rm tr} ( B_\alpha \Omega) =
 \ f_\alpha - \lambda \delta_{\alpha_0} \ \
 \ for \ all \ |\alpha | \leq d.
\eeq
\end{lem}

\begin{proof}.\  Write the symmetric positive semi-definite $n \times n$
matrix
 $ \Omega = \sum_j^r q_j q_j^T$
were $r$ is the rank of $\Omega$ and $q_j \in \mathbb R^n$. Then
\begin{equation}
 {\rm tr} ( \sum_\alpha B_\alpha x^\alpha  \Omega)=
 \ \sum_\alpha f_\alpha x^\alpha - \lambda
\end{equation}
which gives
$$
{\rm tr} (\  V(x) V(x)^T \sum_j^r q_j q_j^T  )
  = f(x) - \lambda
$$
and
$$
 \sum_j  q_j^T V(x) V(x)^T q_j=
f(x) - \lambda
$$
so we obtain
$$  \sum_j  Q_j^T(x)  Q_j(x) = f(x) - \lambda $$
where $  Q_j(x):= V(x)^T q_j$. The argument reverses, so $f$ a sum
of squares implies $\Omega = \sum_j^r q_j q_j^T$ which makes
$\Omega$ positive semi-definite.\qed  \end{proof}

Clearly, there are many matrices $\Omega \geq 0$ satisfying
(\ref{eq:Omegaf}). A canonical choice, appearing in the next
Lemma, was proposed by Nesterov.

\begin{lem}
Suppose there is a positive definite solution to
(\ref{eq:Omegaf}), then one of them ${\breve \Omega}$ has inverse,
${\breve \Omega}^{-1}$ which is a Hankel matrix.
\end{lem}

\begin{proof}. \ \ We show that ${\breve \Omega}$ is
 the  ``maximum entropy solution" to (\ref{eq:Omegaf}),
namely,
 $\breve \Omega$ is the (unique) solution
to
$$
\max_\Omega \ \ln \; \det \  \Omega \ \ \  \ \ \ \  {\rm  subject \ to} \
(\ref{eq:Omegaf}).
$$
It is standard that if a positive definite $\Omega$ exists
maximizing entropy, then it ``keeps its eigenvalues positive", so
$\breve \Omega$ is positive definite. We use the standard formula
$$
\frac{d \ \ln \; \det (\breve \Omega + t \Delta)}{dt}|_{t=0} = -
{\rm tr} [{\breve \Omega}^{-1} \Delta]
$$
which is 0 for all $\Delta$  satisfying $
 {\rm tr}  [B_\alpha \Delta] =0 $ for all  $\alpha$.
Now ${\rm tr} [{\breve \Omega}^{-1} \Delta] =0$ says that
 ${\breve \Omega}^{-1}$ is in the orthogonal
complement of the orthogonal complement of $span \ \{B_\alpha \}$,
thus $ {\breve \Omega}^{-1} \in span \ \{B_\alpha \}$, in other
words, it is a Hankel matrix. \qed \end{proof}

Our next step is the minimization problem:
$$
{\hat L}_f:= \min_L L(f) \ \ \ \ {\rm subject \ to} \ L(\Sigma^2)
\geq 0 \ {\rm and} \ L(1)=1.
$$
Clearly,
$${\hat \lambda}_f \leq \min_\mu \int f d\mu =\min_x f(x)$$
as $\mu$ ranges over all probability measure. If the minimum is
attained, Dirac's measure at the optimal $x$ yields ${\hat
\lambda}_f = \min_x f(x).$

Since $L = L_y$ corresponding to the moment sequence $y$ satisfies
 $ L_y(f) = \sum_\alpha \ f_\alpha y_\alpha$
our basic "primal" problem is:\\
$$(PRIMAL) \qquad \qquad
\min_y  \sum_\alpha f_\alpha y_\alpha \ \ \ \ {\rm subject \ to} \
\ \ M_y \geq 0 \ \ y_{\alpha_0}=1.
$$

The "dual"  problem is
$$(DUAL) \qquad
\max \lambda \ \ \ {\rm  \ subject \  to} \ \ \ \Omega \geq 0 \ \
\ {\rm and \ \ tr} ( B_\alpha \Omega ) \ =  \ f_\alpha - \lambda
\delta_{\alpha_0}
$$
for all $\alpha$ which we saw in Lemma \ref{lem:OmegaSoS} as
solving the sum of squares problem is the same as
$$
\hat \lambda_\Sigma:= \max  \lambda  \ \ \ {\rm  subject \
to} \ \ f-\lambda \ {\rm is \ in} \ \Sigma.
$$
To derive that these problems are indeed  dual to each other
 define a "Lagrangian" by
$$
\mathcal L(y, \Omega, \lambda):=
 \sum_\alpha \ f_\alpha y_\alpha -
{\rm tr} ( \sum_\alpha \; y_\alpha B_\alpha \Omega ) - \ (
y_{\alpha_0}- 1) \lambda
$$
on the set $ M_y \geq 0, \Omega \geq 0, \lambda \in \bbR$.
For convenience write
\beq \mathcal L(y, \Omega, \lambda):=
 \sum_\alpha \ [ ( f_\alpha -
{\rm tr} ( \sum_\alpha \;  B_\alpha \Omega ) - \delta_{\alpha_0}
\lambda ]
 y_\alpha - \  y_{\alpha_0} \lambda.
\eeq Then
$$
\min_{M_y \geq 0} \  \max_{\Omega \geq 0, \lambda} \; \mathcal
L(y, \Omega, \lambda) = \min_{M_y \geq 0}  \sum_\alpha \ f_\alpha
y_\alpha \ \  {\rm if} \ \ \  y_{\alpha_0} =1 \ \ \ {\rm ( is}
+\infty \ {\rm if} \ y_{\alpha_0} \neq 1).
$$
which is the  primal problem. Next
$$
\max_{\Omega \geq 0, \lambda} \ \min_{M_y \geq 0} \;
 \mathcal L(y, \Omega , \lambda)
 = \max  \lambda  \ \ {\rm if} \ \
 f_\alpha  - {\rm tr} (  B_\alpha \Omega )
- \lambda \delta_{\alpha_0} = 0 \ \ \  ({\rm is}  \  \ -\infty \
{\rm otherwise}).
 $$
 which is the dual problem.
 We summarize with

 \begin{lem} For $\Omega, M_y $ in $\Rnn$,
define sets $$\cN:= \{\Omega \geq0 , \lambda \in \bbR :
 \ f_\alpha  - {\rm tr} (  B_\alpha \Omega )
- \lambda \delta_{\alpha_0} = 0 \ \ all \ |\alpha| \leq d \}$$
$$\cM:= \{ y: \ M_y \geq 0,\; y_{\alpha_0} =1 \}.$$
Then $\cN$ is not empty for some dimension $n$ if and only if
$f - \lambda$ is a sum of squares and we have
$$
 {\hat \lambda_\Sigma}:= \max_{ \Omega, \lambda \in \cN} \lambda =
 \max_{ \Omega, \lambda \in \cN}
\min_{y \in \cM} \;
 \mathcal L(y, \Omega , \lambda)
 \qquad \qquad\qquad\qquad
 $$
 $$
 \qquad  \qquad\qquad\qquad \leq
\min_{y \in \cM} \  \max_{ \Omega, \lambda \in \cN} \; \mathcal
L(y, \Omega, \lambda)
= \min_{y \in \cM}  \sum_\alpha \ f_\alpha
y_\alpha =:
 {\hat \lambda_f}
$$
 \end{lem}

A saddle point $\hat \Omega, \hat
\lambda, \hat y$ is defined as one which satisfies
$$
 {\hat \lambda_\Sigma}
= \mathcal L(\hat y, \hat \Omega,\hat \lambda)
= {\hat \lambda_f}.
$$
The lemma verifies our claim that
 the problems PRIMAL and DUAL are dual with respect to each
other. Also if they have a saddle point they have the same optimal
value $\hat \lambda:=
 {\hat \lambda_\Sigma}=
 {\hat \lambda_f}$.
Existence of a saddle point,
 because of the bilinearity of $\mathcal L$,
 is in the perview of the
 von Neumann Minmax Theorem, but we do not discuss this here.
We refer the reader to \cite{Lasserre0, Lasserre1,
Henrion3} for further details. Although their approach
is a bit different.

Now we make a few remarks.
 Firstly, the saddle point condition
$$
\hat \lambda = \mathcal L(\hat y, \hat \Omega,\hat \lambda)=
 \sum_\alpha \ f_\alpha \hat y_\alpha -
{\rm tr} ( \sum_\alpha \; \hat y_\alpha B_\alpha \hat \Omega ) - \
( \hat y_{\alpha_0}- 1) \hat \lambda
$$ reduces to
$ 0= {\rm tr} ( \sum_\alpha \; \hat y_\alpha B_\alpha \hat \Omega
) =tr ( M_{\hat y} \hat \Omega ) $. Since $ M_{\hat y}, \hat
\Omega$ are both positive semi-definite, this forces the
``complementarity" of an optimal moment matrix and an optimal sum
of squares representor $\hat \Omega$. \beq
 M_{\hat y} \hat \Omega=0.\eeq
 Secondly, an  equation ``balanced" between
 primal and dual   is
 $$
tr ( M_y \Omega ) = L_y (f) - \lambda.
 $$

We will continue this line of thought in a separate article.

\subsection{Engineering}

For nonlinear systems,
sum of squares techniques can be
used to find Lyapunov functions by direct computation.
Here the problem is to check if a differential equation
$$ \frac{d x}{dt}= a(x)$$
on $\RR^g$ is stable.
The most common technique is to seek a function $V >0$
except $V(0)=0$ satisfying the differential inequality
 $$\nabla V(x) \cdot a(x) \leq 0 \ \ \ {\rm for  \ all} \ x \ in \ \RR^g;$$
such $V$ are called {\it Lyapunov functions} \index{Lyapunov
functions}. If $a$ is a vector field with  polynomial entries, it
is natural to seek $V$ which is a sum of squares of polynomials,
and this reduces to a semi-definite program. Solution can be
attempted numerically and if successful produces $V$ a Lyapunov
function; if not one can modify the sum of squares from
polynomials to   some rational sum of squares and try again, see
\cite{parThesis}. Such  techniques lay out numerical solutions to
midsized (dimension 8 or so) nonlinear control problems.

More generally for control problems one seeks
to find a feedback law $u= k(x)$ which stabilizes
$$ \frac{d x}{dt}= a(x)+ b(x) u.$$
Here there is a function $V$, beautifully tamed by E. Sontag,
called a ``control Lyapunov function" generalizing the classical
Lyapunov function. Unfortunately, no Positivstellens\"atz
technique is known for for finding $V$. However, A. Rantzer
cleverly introduced a
 ``dual control Lyapunov function"
and in \cite{PPR04} showed that it
is quite amenable to sum of squares techniques.

Another direction
gives a generalization of the
classical S-procedure which finds performance
bounds on broad classes of problems,
see \cite{parThesis}.
There are also applications to
combinatorial problems described there.

Recently, \cite{Henrion3} have given a technique for converting
system engineering problems to polynomial minimization. The wide
scope of the technique is very appealing.


\section{Linear matrix inequalities and computation of sums of squares}
\label{sec:LMIs and SoS}

Numerical computation of a sum of squares and a Positivstellensatz
is based on a revolution which started about 20 years ago in
optimization; the rise of interior point methods. We avoid delving
into yet another topic but mention the special aspects concerning
us. Thanks to the work of Nesterov and Nemirovskii in the early
1990s one can solve Linear Matrix Inequalities (LMIs in short)
numerically using interior point optimization methods, called {\it
semi-definite programming} \index{semi-definite programming}. An
LMI is an inequality of the form
\begin{equation}
A_0 + A_1 x_1 + \cdots A_g x_g  \geq 0
\end{equation}
\index{LMIs, Linear Matrix Inequalities}where the $A_j$ are
symmetric matrices and the numerical goal is to compute $ x \in
\RR^g $ satisfying this. The sizes of matrix unknowns treatable by
year 2006 solvers exceed 100 $\times$ 100; with special structure
dimensions can go much higher. This is remarkable because our LMI
above has about $5000 g$ unknowns.

\subsection{SOS and LMIs}

Sum of squares and Positivstellens\"atze problems convert readily
to LMIs and these provide an effective solution for polynomials
having modest number of terms. These applications  make
efficiencies in numerics a high priority. This involves shrewd use
of semi-algebraic theory and computational ideas to produce a
semi-definite programming package, for a recent paper see
\cite{Parr05}; also there is recent work of L. Vandenberghe.
Semi-algebraic geometry  packages are: SOS tools \cite{PPSP04} and
GloptiPoly  \cite{Henrion1}.

A lament is that  all current computational semi-algebraic
geometry projects use a packaged semi-definite solver, none write
their own. This limits efficiencies for sum of squares
computation.

Special structure leads to great computational improvement as well
as elegant mathematics. For example, polynomials which are
invariant under a group action, the delight of classical invariant
theory, succumb to rapid computation, see \cite{GP04}
\cite{CKSprpept}.

\subsection{LMIs and the world}
LMIs have a life extending far beyond
computational sum of squares and are being found in many
areas of science.
Later in this paper \S \ref{sec:engLMI} we shall
glimpse at their use in systems engineering,
a use preceding sum of squares applications by 10 years.
The list of other areas includes statistics,
chemistry, quantum computation together with
more; all to vast for us to attempt description.

A paradigm mathematical question here is:\\

 {\it Which convex
sets $\cC$ in $\RR^g$ with algebraic boundary can be represented
with some monic LMI?}\\

 That is,
$$ \cC = \{ x \in \RR^g : \ I + A_1 x_1
 + \cdots A_g x_g  \geq 0 \},$$
where $A_j$ are symmetric matrices. Here we have assumed the
normalization $0 \in \cC$. This question was raised by Parrilo and
Sturmfels \cite{PS}. The paper \cite{HVprept} gives an obvious
necessary condition \footnote{This is in contrast to the free
algebra case where all evidence (like that in this paper)
indicates that convexity is the only condition required.} on $\cC$
for an LMI representation to exist and proves sufficiency when
$g=2$.

The main issue is that of determinantal representations
\index{determinantal representation} of a polynomial $p(x)$ on
$\RR^g$, namely, given $p$ express it in the form
\begin{equation}
p(x) =det(A_0 + A_1 x_1 + \cdots A_g x_g ).
\end{equation}
That this is possible for some matrices is due to the computer
scientist Leslie Valiant \cite{Val}. That the matrices can be
taken real and symmetric is in \cite{HMVprept} as is the fact the
a representation of $\det \ p(X)$ always holds for polynomials in
non-commuting (free) variables, as later appear in \S
\ref{sec:nonCom}. A symbolic computer algorithm due to N.
Slinglend and implemented by J. Shopple runs under the Mathematica
package NCAlgebra.

The open question is which polynomials can we represent monicaly;
that is with $A_0=I$. Obviously, necessary is the {\it real zero
condition} \index{real zero condition}, namely,

\begin{center}
{\it the polynomial
$f(t):= p(tx)$  in one complex variable $t$ \\ has only real
zeroes,}
\end{center}
but what about the converse? When $g=2$ the real zero
condition on $p$ insures that it has a monic representation; this
is the core of \cite{HVprept}.

What about higher dimensions? Lewis, Parrilo and Ramana
\cite{LPR05} showed that this $g=2$ result (together with a
counterexample they concocted) settles a 1958 conjecture of Peter
Lax, which leads to the surmise that sorting out the $g >2$
situation may not happen soon. Leonid Gurvitz pointed out the
Valient connection to functional analysts and evangelizes that
monic representations have strong implications for lowering the
complexity of certain polynomial computations.

\section{Non-commutative algebras}
\label{sec:nonCom}

A direction in semi-algebraic geometry, recently blossoming  still
with many
 avenues to explore,  concerns variables
 which do not commute.
As of today  versions of the strict Positivstellens\"atze we saw
in \S \ref{sec:PosSS} are proved for  a free $*$- algebra and for
the enveloping algebra of a Lie algebra; here the structure is
cleaner or the same as in the classical commutative theory. The
verdict so far on noncommutative Nullstellens\"atze is mixed. In a
free algebra it goes  through so smoothly that no radical ideal is
required. This leaves us short of the remarkable perfection we see
in the Stengle -Tarski - Seidenberg commutative landscape. Readers
will be overjoyed to hear that the proofs needed above are mostly
known to them already: just as in earlier sections, non-negative
functionals on the sums of squares cone in a $\ast$-algebra can be
put in correspondence with tuples of non-commuting operators, and
this carries most of the day.

This noncommutative semi-algebraic foundation underlies  a rigid
structure (at least) for free $*$-algebras which has recently
become visible. A noncommutative polynomial $p$ has second
derivative $p''$ which is again a polynomial and if $p''$  is
positive, then our forthcoming free $*$-algebra Positivstellensatz
tells us that $p''$ is  a sum of squares. It is a bizarre twist
that this and the derivative structure are incompatible, so
together imply that a ``convex polynomial" in a free $*$- algebra
has degree 2 or less; see \S \ref{sec:convexity}. The authors
suspect that this is a harbinger of a very rigid structure in a
free $*$-algebra for ``irreducible varieties" whose curvature is
either nearly positive or nearly negative; but this is a tale for
another (likely distant) day. Some of the material in
this section on higher derivatives and the next
is new.

A final topic on semi-algebraic geometry in a free $*$- algebra is
applications to engineering, \S \ref{sec:engLMI}. Arguably the
main practical development in systems and control through the
1990's was the reduction of linear systems problems to  Linear
Matrix Inequalities, LMIs.
For theory and numerics to be highly
successful something called
 ``Convex Matrix Inequalities", henceforth denoted in short CMIs,
\index{CMIs Convex Matrix Inequalities} will do nicely. Most
experts would guess that the class of problems treatable with CMIs
is much broader than with LMIs. But no, as we soon see, our
draconian free $*$ convexity theorems suggest that for systems
problems fully characterized by performance criteria based on
$L^2$ and signal flow diagrams (as are most textbook classics),
convex matrix inequalities give no greater generality than LMIs.

These systems problems have the key feature that their statement
does not depend on the dimension of the systems involved. Thus we
summarize our main engineering contention:

\begin{center}
{\it Dimension free convex problems are equivalent to an LMI}\\
\end{center}

This and the next sections tells the story we just described
but there is a lot it does not do.
Our focus in this paper has been on inequalities,
where  various noncommutative equalities are of course a special
and often well developed case.
For example,
algebraic geometry based on the Weyl algebra and corresponding
computer algebra implementations, for example, Gr\"obner basis
generators for the Weyl algebra are in the standard computer
algebra packages such as Plural/Singular.

A very different and elegant  area is that of rings with a
polynomial identity, in short PI rings \index{ PI rings}, e.g. $N
\times N$ matrices for fixed $N$. While most PI research concerns
identities, there is one line of work on polynomial inequalities,
indeed sums of squares, by Procesi-Schacher \cite{Procesi}.
A Nullstellensatz for PI rings is discussed in \cite{Amitsur}.

\subsection{Sums of squares in a free $\ast$-algebra}
\label{sec:SoSfree}

\def\hh{q}

Let $ \FF$
\index{${\mathbb{R}} \langle x, x^\ast \rangle$} denote the
polynomials with real numbers as coefficients in variables
$x_1,...,x_g, x_1^\ast,...,x_g^\ast$.
These variables do not
commute, indeed they are free of constraints other than
 $^*$ being an anti-linear involution:
$$ (f \hh)^\ast = \hh^\ast f^\ast, \ \ \  \ \ (x_j)^\ast = x_j^\ast.$$
Thus $\FF$ is called the {\it real free $*-$ algebra}
\index{free $*-$ algebra} on generators $x, x^*$.

Folklore has it that analysis in  a free $*$-algebra gives results
like ordinary commutative analysis in one variable.
The SoS
phenomenon we describe in this section  is consistent with this
picture, but convexity properties in the next section do not.
Convexity in a free algebra is much more rigid.

We invite those who work in a free algebra (or their students) to
try NCAlgebra, the free free-$\ast$ algebra computer package
\cite{HSM05}. Calculations with it had  a profound impact on the
results in \S \ref{sec:nonCom} and \ref{sec:convexity}; it is a
very powerful tool.

The cone of sums of squares is the convex hull:
$$ \Sigma^2 = {\rm co} \{ f^\ast f; \  f \in \FF \}.$$
A linear functional $L \in \FF'$ satisfying $L|_{\Sigma^2} \geq 0$
produces a positive semidefinite bilinear form
$$ \langle f, \hh \rangle = L(\hh^\ast f)$$ on $\FF$.
We use the same construction
introduced in section \ref{sec:moments},
namely, mod out the null
space of $\langle f, f \rangle$ and denote the
 Hilbert space completion by $H$, with
$\mathcal D$  the dense subspace of $H$ generated by $\FF$.
The  separable Hilbert space  $H$ carries the
multiplication operators
$M_j : \mathcal D \longrightarrow \mathcal D$:
$$ M_j f = x_j f, \ \ f \in \mathcal D, \ 1 \leq j \leq n.$$
One verifies from the definition that each $M_j$ is well defined and
$$ \langle M_j f , \hh \rangle = \langle x_j f,  \hh  \rangle = \langle
f, x_j^\ast  \hh  \rangle, \ \ f, \hh  \in \mathcal D.$$ Thus $M_j^\ast =
M_{x_j^\ast}$. The vector $1$ is still $\ast$-cyclic, in the sense
that the linear span $\vee_{p \in \FF} p(M,M^\ast)1$ is dense in $H$.
Thus, {\it mutatis mutandis}, we have obtained the following
result.

\begin{lem}
There exists a bijective correspondence between positive
linear functionals, namely
$$ L \in \FF'\ \  \ \ and \ \  \ \ L_{|{\Sigma^2} } \geq 0,$$
and $g$-tuples of unbounded linear operators $T$ with a star
cyclic vector $\xi$, established by the formula
$$ L(f) = \langle f(T,T^\ast) \xi, \xi \rangle,
\ \  \ \ f \in \FF.$$
\end{lem}

We stress that the above operators do not commute, and might be
unbounded. The calculus $f(T,T^\ast)$ is the non-commutative
functional calculus: $x_j(T) = T_j, \
   x_j^\ast (T) = T_j^\ast.$

An important feature of the above correspondence is that it can be
restricted by the degree filtration. Specifically, let
$\FF_k = \{ f ; \ {\rm deg} f \leq k\}$,
and similarly, for a quadratic form $L$ as in the
lemma, let $\mathcal D_k$ denote the finite dimensional subspace
of $H$ generated by the elements of $\FF_k$.
Define also
$$ \Sigma^2_k = \Sigma^2 \cap \FF_k.$$
\index{$ \Sigma^2_k $}

Start with a functional $L \in \FF_{2k}'$ satisfying
$L|_{\Sigma^2_{2k}} \geq 0.$ One can still construct a finite
dimensional Hilbert space $H$, as the completion of $\FF_k$ with
respect to the inner product
$\langle f,  \hh  \rangle = L( \hh^\ast f), \ f, \hh  \in \FF_k.$
The multipliers
  $$ M_j : \mathcal D_{k-1} \longrightarrow H, \ M_j f = x_j
  f,$$
are well defined and can be extended by zero to the whole $H$.
Let
$$N(k) = \dim \FF_k = 1 + (2g) + (2g)^2 +...+ (2g)^k =
\frac{(2g)^{k+1}-1}{2g-1}.$$
In short, we have proved the
following specialization of the main Lemma.

\begin{lem}\label{fd lemma}
Let
$L \in \FF_{2k}'$
satisfy
$L|_{\Sigma^2_{2k}} \geq 0$.
 There exists a Hilbert space of
dimension $N(k)$ and an $g$-tuple of linear operators $M$ on $H$,
with a distinguished vector $\xi \in H$, such that
\begin{equation}\label{fd representation}
L(p) = \langle p(M,M^\ast) \xi, \xi \rangle, \ \ p \in
\FF_{2k-2}.
\end{equation}
\end{lem}

Following the pattern of the preceding section, we will derive now
a Nichtnegativstellensatz.

\begin{thm}[\cite{H}]
\label{free SOS}
\index{Sum of Squares Theorem, Free}
Let $p \in \FF_d$ be a non-commutative
polynomial satisfying $p(M,M^\ast) \geq 0$ for all $g$-tuples of
linear operators $M$ acting on a Hilbert space of dimension at most $N(k), \
2k \geq d+2$. Then $ p \in \Sigma^2$.
\end{thm}

\begin{proof}.
The only necessary technical result we need is the closedness of
the cone $\Sigma^2_k$ in the Euclidean topology of the finite
dimensional space $\FF_k$. This is done as in the commutative
case, using Carath\'edodory's convex hull theorem. More exactly,
every element of $\Sigma^2_k$ is a convex combination of at most
$\dim \FF_k +1$ elements, and on the other hand there are finitely
many positive functionals on $\Sigma^2_k$ which separate the
points of $\FF_k$. See for details \cite{HMP1}.

Assume that $p \notin \Sigma^2$ and let $k \geq (d+2)/2$, so that
$p \in \FF_{2k-2}$. Once we know that $\Sigma_{2k}^2$ is a closed
cone, we can invoke Minkowski separation theorem and find a
functional $L \in \FF_{2k}'$
providing the strict separation:
$$
L(p)<0 \leq L(f), \ \ f \in \Sigma_{2k}^2.$$
According to Lemma \ref{fd lemma} there exists a tuple $M$ of
operators acting on a Hilbert space $H$ of dimension $N(k)$ and a
vector $\xi \in H$, such that
$$ 0 \leq  \langle p(M,M^\ast) \xi, \xi \rangle = L(p) <0,$$
a contradiction.  \ \qed
\end{proof}

When compared to the commutative framework, this theorem is
stronger in the sense that it does not assume a strict positivity
of $p$ on a well chosen "spectrum". Variants with supports (for
instance for spherical tuples $M: \ M_1^\ast M_1 + ...+ M_g^\ast
M_g \leq I$) of the above result are discussed in \cite{HMP1}.

We state below an illustrative and generic result,
from \cite{HM04pos}, for sums of
squares decompositions in a free $\ast$-algebra.

\begin{thm}
\index{Free Positivstellensatz, with supports}
Let $p \in \FF$  and let $q = \{ q_1,...,q_k\} \subset \FF$
be a set of polynomials, so that the non-commutative quadratic
module
$$QM(q) = {\rm co} \{ f^\ast q_k f; \ f \in \FF, \ 0 \leq i \leq
k\}, \ q_0 =1,$$ contains $1-x_1^\ast x_1 -...-x_g^\ast x_g$ . If
for all tuples of linear bounded Hilbert space operators $X =
(X_1,...,X_g)$ subject to the conditions
$$ q_i(X,X^\ast) \geq 0, \ 1 \leq i \leq k,$$
we have
$$ p(X,X^\ast) > 0,$$
then $p \in QM(q)$.
\end{thm}

Notice that the above theorem covers relations of the form
$r(X,X^\ast) =0$, the latter being assured by $\pm r \in QM(q)$.
For instance we can assume that we evaluate only on commuting
tuples of operators, in which situation all commutators $[x_i,
x_j]$ are included among the (possibly other) generators of
$QM(q).$

Some interpretation is needed in degenerate cases, such as those
where no bounded operators satisfy the relations $q_i(X,X^\ast)
\geq 0$, for example, if some of $q_i$ are the defining relations
for the Weyl algebra; in this case, we would say $ p(X,X^\ast) >
0$, since there are no $X$. Indeed $p \in QM(q)$ as  the theorem
says.

\bigskip

\begin{proof} Assume that $p$ does not belong to
the convex cone $QM(q)$. Since the latter is archimedean, by the
same Minkovski principle there exists a linear functional $L \in
\FF'$, such that
$$ L(p) \leq 0 \leq L(f), \ \ f \in QM(q).$$
Define the Hilbert space $H$ associated to $L$, and remark that
the left multipliers $M_{x_i}$ on $\FF$ give rise to linear
bounded operators (denoted by the same symbols) on $H$. Then
$$ q_i(M, M^\ast) \geq 0, \ \ 1 \leq i \leq k,$$
by construction, and
$$ \langle p(M,M^\ast)1, 1\rangle = L(p) \leq 0,$$
a contradiction.
\end{proof}
\bigskip

The above statement allows a variety of specialization to quotient
algebras. Specifically, if $I$ denotes a bilateral ideal of $\FF$,
then one can replace the quadratic module in the statement with
$QM(q) + I$, and separate the latter convex cone from the
potential positive element on the set of tuples of matrices $X$
satisfying simultaneously
$$ q_i(X,X^\ast) \geq 0, \ 0 \leq i \leq k, \ \ f(X) = 0, \ \ f
\in I.$$

For instance, the next simple observation can also be deduced from
the preceding theorem.

\begin{cor} Let $J$ be the bilateral ideal of $\FF$ generated by the
commutator polynomial $[x_1+x_1^\ast, x_2+x_2^\ast] - 1$. Then $J
+ QM(1-x_1^\ast x_1 -...- x_g^\ast x_g) = \FF$.
\end{cor}

\begin{proof} Assume by contradiction that $ J +
QM(1-x_1^\ast x_1 -...- x_g^\ast x_g) \neq \FF$. By our basic
separation lemma, there exists a linear functional $L \in \FF'$
with the properties:
$$ L_{J
+ QM(1-x_1^\ast x_1 -...- x_g^\ast x_g)} \geq 0, \ \ {\rm and}\ \
L(1)>0.$$

Then the GNS construction will produce a tuple of linear bounded
operators $X$, acting on the associated non-zero Hilbert space
$H$, satisfying $ X_1^\ast X_1 +...+X_g^\ast X_g \leq I$ and
$$ [X_1^\ast + X_1, X_2^\ast + X_2] = I.$$
The latter equation is however impossible, because the left hand
side is anti-symmetric while the right hand side is symmetric and
non-zero.
\end{proof} \qed
\bigskip

Similarly, we can derive following the same scheme the next
result.

\begin{cor} Assume, in the condition of the above Theorem, that
$p(X,X^\ast) >0$ for all {\sc commuting} tuples $X$ of matrices
subject to the positivity constraints $q_i(X,X^\ast) \geq 0, 0
\leq i \leq k.$ Then
$$ p \in QM(q) + I,$$
where $I$ is the bilateral ideal generated by all commutators
$[x_i, x_j], [x_i, x_j]^\ast, \ 1 \leq i,j \leq g.$
\end{cor}
\bigskip

With similar techniques (well chosen, separating,
$\ast$-representations of the free algebra)
one can prove a series of Nullstellens\"atze.
We state for information one of them, see
for an early version \cite{HMP2}.

\begin{thm}
\label{thm:NullSS}
  Let $p_1(x),...,p_m(x) \in \RRx$
  be polynomials not depending on the $x_j^\ast$ variables and let
  $q(x,x^\ast) \in \FF$. Assume that for every $g$ tuple $X$ of linear
  operators acting on a finite dimensional Hilbert space $H$, and
every vector $v \in
  H$, we have:
  $$ (p_j(X)v = 0, \ 1 \leq j \leq m) \ \Rightarrow \ (q(X,X^\ast)v =
  0).$$
  Then $q$ belongs to the left ideal $\FF p_1 + ... + \FF p_m$.
\end{thm}

Again, this proposition is stronger than its commutative
counterpart. For instance there is no need of taking higher powers
of $q$, or of adding a sum of squares to $q$.

We refer the reader to \cite{HMP3} for the proof of Proposition
\ref{thm:NullSS}. However, we say a few words about the intuition
behind it. We are assuming
$$ p_j(X)v =0, \forall j \ \ \ \implies \ \ \ q(X,X^\ast)v =0.$$
On  a very large vector space if $X$ is determined on a small
number of vectors, then $X^\ast$ is not heavily constrained; it is
almost like being able to take $X^\ast$ to be a completely
independent tuple $Y$.
If it were independent, we would have
$$ p_j(X)v =0, \forall j \ \ \ \implies \ \ \ q(X,Y)v =0.$$

Now, in the free algebra $\mathbb R \langle x,y \rangle$, it is
much simpler to prove that this implies
$q \in \ \sum_j^m  \mathbb R \langle x,y \rangle \; p_j $,
as required. We isolate this fact in a separate lemma.

\begin{lem}
Fix a finite collection  $p_1,...,p_m$ of polynomials in
non-commuting variables $\{x_1,\dots,x_g \}$ and let  $q$ be a
given polynomial in $\{x_1,\dots,x_g\}$.
Let $d$ denote the
maximum of the $\mbox{\rm deg}(q)$ and $\{\mbox{\rm deg}(p_j): 1 \leq j
\leq m \}$.

There exists a real Hilbert space $\mathcal H$ of dimension
$\sum_{j=0}^d g^j$, such that, if
\begin{equation*}
q(X)v=0
\end{equation*}
whenever $X=(X_1,\dots,X_g)$ is a tuple of operators on $\mathcal
H$, $v\in\mathcal H$, and
\begin{equation*}
p_j(X)v=0 \mbox{ for all } j,
\end{equation*}
then $q$ is in the left ideal generated by $p_1,...,p_m$.
\end{lem}

\begin{proof} (of Lemma). We sketch a proof based on an idea
of G. Bergman,
see
\cite{HM04pos}.

Let $\cI$ be the left ideal generated by $p_1,...,p_m$ in $F =
\mathbb R \langle x_1,...,x_g \rangle.$ Define $\cV$ to be the
vector space $F/\cI$ and denote by $[f]$ the equivalence class of
$f \in F$ in the quotient $F/\cI$.

Define $X_j$ on the vector space $F/\cI$ by $X_j [f]=[x_j f]$ for
$f \in F$, so that $x_j \mapsto X_j$ implements a quotient of the
left regular representation of the free algebra $F$.

  If $ \cV:= F/\cI$ is finite dimensional, then the linear
operators $X=(X_1,\dots,X_g)$ acting on it can be viewed as a
tuple of matrices and we have, for $f\in F$,
\begin{equation*}
f(X)[1]=[f].
\end{equation*}
In particular,  $p_j(X)[1]=0$ for all $j$. If we do not worry
about the dimension counts, by assumption, $0=q(X)[1]$, so $ 0
=[q]$ and therefore $q \in \cI$. Minus the precise statement about
the dimension of $\mathcal H$ this establishes the result when
$F/\cI$ is finite dimensional.

  Now we treat the general case where we do not assume finite
dimensionality of the quotient. \  \ Let $\cV$ and $\cW$ denote
the vector spaces
$$
\cV:= \{[f]: \ f \in F, \ \deg (f) \le d\},
$$ $$
  \cW:= \{[f]: \ f \in F, \ \deg (f) \le d-1 \}.$$
Note that the
dimension of $\cV$ is at most $\sum_{j=0}^d g^j$. We define $X_j$
on $\cW$ to be multiplication by $x_j$. It maps $\cW$ into $\cV$.
Any linear extension of $X_j$ to the whole $\cV$ will satisfy: if
$f$ has degree at most $d$, then $f(X)[1]=[f]$. The proof now
proceeds just as in the part 1 of the proof above. \ \qed
\end{proof}

With this observation we can return and finish the proof of
Theorem \ref{thm:NullSS} Since $X^\ast$ is dependent on $X$, an
operator extension with properties stated in the lemma below
gives just enough structure to make the above free algebra
Nullstellensatz apply; and we prevail.

\begin{lem}
\label{prop:lift}
  Let
  $x = \{x_1, \ldots, x_m \},\  y = \{y_1, \ldots, y_m\}$ be free,
  non-commuting variables.
  Let $H$ be a finite dimensional Hilbert space,
        and let
        $X,Y$ be two $m$-tuples of linear operators acting on
        $H$. Fix a degree $d \geq 1$.

  Then there exists a larger Hilbert space $K \supset H$,
        an $m$-tuple of linear transformations $\tX$ acting on $K$,
        such that
$$
\tX_j|_H = X_j, \ \  1 \leq j \leq g, $$
        and for every polynomial
        $q \in \mathbb R \langle x, x^\ast \rangle$
        of degree at most $d$  and vector $v \in H$,
        $$
        q(\tX,\tX^\ast)v  = 0 \ \Rightarrow  q(X,Y)v = 0.
$$
\end{lem}
For the matrical construction in the proof see \cite{HMP3}.
\bigskip

We end this subsection with an example, see \cite{HM04pos}.

\begin{example}
  \label{ex:unreal}
Let $p=(x^\ast x+xx^\ast)^2$ and $q=x+x^\ast$ where $x$ is a
single variable. Then, for every matrix $X$ and vector $v$
(belonging to the space where $X$ acts), $p(X)v=0$ implies
$q(X)v=0$; however,
  there does not exist
  a positive integer $m$ and $r,r_j \in \mathbb R \langle x,
x^\ast\rangle,$ so that
  \begin{equation}
   \label{eq:unrealrep}
     q^{2m}+\sum r_j^\ast r_j =  pr + r^\ast p .
  \end{equation}
Moreover,  we can modify the example to add the condition $p(X)$
is  positive semi-definite  implies $q(X)$ is  positive
semi-definite and still not obtain this representation.
\end{example}

\proof \ \ Since $A:=XX^\ast+X^\ast X$ is self-adjoint, $A^2 v=0$
if and only if $Av=0$. It now follows that if $p(X)v=0$, then
$Xv=0=X^\ast v$ and therefore $q(X)v=0$.

For $\lambda \in \mathbb R$, let
\begin{equation*}
  X=X(\lambda)=\begin{pmatrix} 0 & \lambda & 0 \\ 0 & 0 & 1 \\ 0&0&0
\end{pmatrix}
\end{equation*}
viewed as an operator on $\mathbb R^3$ and let $v=e_1$, where
$\{e_1,e_2,e_3\}$ is the standard basis for $\mathbb R^3$.

We begin by calculating the first component of even  powers of the
matrix $q(X)$. Let $Q=q(X)^2$ and verify,
\begin{equation}
  \label{eq:unreal1}
  Q=\begin{pmatrix} \lambda^2 &0&\lambda \\ 0& 1+\lambda^2 & 0\\
\lambda & 0 & 1\end{pmatrix}.
\end{equation}
For each positive integer $m$ there exist a  polynomial $q_m$ so
that
\begin{equation}
  \label{eq:unreal2}
    Q^m e_1 =\begin{pmatrix} \lambda^2 (1+\lambda q_m(\lambda)) \\ 0
\\ \lambda(1+\lambda q_m(\lambda)) \end{pmatrix}
\end{equation}
which we now  establish by an induction argument. In the case
$m=1$, from equation (\ref{eq:unreal1}), it is evident that
$q_1=0$. Now suppose equation (\ref{eq:unreal2}) holds for $m$.
Then,  a computation of $QQ^me_1$ shows that equation
(\ref{eq:unreal2}) holds for $m+1$ with $q_{m+1}=\lambda (q_m
+\lambda +\lambda q_m)$. Thus, for any $m$,
\begin{equation}
  \label{eq:unreal3}
  \lim_{\lambda \to 0} \frac{1}{\lambda^2} <Q^m e_1,e_1> =
  \lim_{\lambda \to 0} ( 1 +\lambda q_m( \lambda) ) =1.
\end{equation}

Now we look at $p$ and get
\begin{equation*}
  p(X)=\begin{pmatrix} \lambda^4 & 0&0 \\ 0 & (1+\lambda^2)^2 &
0\\0&0& 1 \end{pmatrix}.
\end{equation*}
Thus
\begin{equation*}
  \lim_{\lambda \to 0} \frac{1}{\lambda^2} \left (< r(X)^\ast
p(X)e_1,e_1>+<p(X)r(X) e_1,e_1>\right)  =0.
\end{equation*}

If the representation of equation (\ref{eq:unrealrep}) holds, then
apply $< \cdot \ e_1, e_1 >$ to both sides and take $\lambda$ to
0. We just saw that the right side is 0, so the left side is 0,
which because
$$<\sum r_j(X)^\ast r_j(X)e_1,e_1> \; \ge 0$$
forces
\begin{equation*}
  \lim_{\lambda \to 0} \frac{1}{\lambda^2} <Q^m e_1,e_1> \;  \leq 0
\end{equation*}
a contradiction to equation ( \ref{eq:unreal3} ). Hence the
representation of  equation (\ref{eq:unrealrep}) does not hold.

The last sentence claimed in the example is true when we use the
same polynomial $p$ and replace $q$ with $q^2$.  \ \qed

There are more Positivstellens\"atze in a free *-algebra which
fill in more of the picture. The techniques proving them are not
vastly beyond what we illustrated here. For example,
Klep-Schweighofer \cite{KSprept} \index{Positivstellensatz, Free
(a)} do an analog of Stengle's Theorem
\ref{thm:stengle}(a), while Theorem \ref{quadratic module} is
faithfully made free in \cite{HM04pos}. In spite of the above
results we are still far from having a full understanding (\`a la
Stengle's Theorem) of the Null- and Positiv-stellens\"atze
phenomena in the free algebra.

\subsection{The Weyl algebra}
Weyl's algebra, that is the
enveloping algebra of the Heisenberg group is interesting because,
by a deep result of Stone- von Neumann, it has a single
irreducible representation; and that is infinite dimensional.
Thus, to check on the spectrum the positivity of an element, one
has to do it at a single point. The details were revealed by
Schm\"udgen in a very recent article \cite{Schmudgen2}. We
reproduce from his work the main result.

Fix a positive integer $g$ and consider the unital $\ast$-algebra
$W(g)$ generated by $2g$ self-adjoint elements $p_1,...,p_g,
q_1,...,q_g$, subject to the commutation relations:
$$ [p_k, q_j] = - \delta_{kj} (i \cdot 1),\ \ [p_k,p_j] = [q_j,q_k]
=0, \ \ 1 \leq j,k \leq g.$$ The unique irreducible representation
$\Phi$ of this algebra is given by the partial differential
operators
$$ \Phi(p_k)f = -i \frac{\partial f}{\partial x_k}, \ \ \Phi(q_k) f =
x_k f,$$ acting on Schwartz space $\mathcal S (\mathbb R^g)$. Via
this representation, the elements of $W(g)$ are identified with
linear partial differential operators with polynomial coefficients
(in $g$ variables). These operators can be regarded as densely
defined, closed graph operators from $\mathcal S(\mathbb R^g)$ to
$L^2(\mathbb R^g)$. The set
$$W(g)_+ = \{ f \in W(g); \ \ \langle \Phi(f)\xi, \xi \rangle \geq
0, \ \xi \in \mathcal S(\mathbb R^g)\}$$ consists of all symmetric,
non-negative elements, with respect to the representation $\Phi$.

Define
$$ a_k = \frac{q_k + ip_k}{\sqrt{2}}, \ \ a_{-k} = \frac{q_k -
ip_k}{\sqrt{2}},$$ so that $a_k^\ast = a_{-k}$. Fix a positive
number $\alpha$ which is not an integer, and let
$$ N = a_1^\ast a_1+...+a_g^\ast a_g;$$
denote by $\mathcal N$ the set of all finite products of elements
$N + (\alpha+n)1$, with $n \in \mathbb Z$.

The algebra $W(g)$ carries a natural degree, defined on generators
as
$$ \deg (a_k) = \deg(a_{-k}) = 1.$$
Every element $f \in W(g)$ can be decomposed into homogeneous
parts $f_s$ of degree $s$:
$$ f = f_m + f_{m-1} +...+f_0.$$
We can regard $f_k$ as a homogeneous polynomial of degree $k$, in
the variables $a_{\pm 1},...,a_{\pm g}.$ The principal symbol of
$f$ is the polynomial\\
$f_m(z_1,...,z_g,\overline{z}_1,...,\overline{z}_g)$, where $a_j$
was substituted by $z_k$ and $a_{-k}$ by $\overline{z_k}$.

\begin{thm}
\cite{Schmudgen2}
\index{Sum of Squares Theorem, Weyl Algebra}
Let $f \in W(g)$ be a self-adjoint element
of even degree $2m$, and let $P(z,\overline{z})$ be its principal
symbol. If

a). There exists $\eps >0$ such that $f - \eps \cdot 1 \in
W(g)_+$,

b). $P(z,\overline{z})>0$ for $z \neq 0$, \\
then, if $m$ is even
there exists $b \in \mathcal N$ such that $bfb \in \Sigma^2 W(g)$;
if $m$ is odd, there exists $b \in \mathcal N$ such that
$\sum_{j=1}^g b a_j f a_{-j} b \in \Sigma^2 W(g).$
\end{thm}

For examples and details see \cite{Schmudgen2}.

Already mentioned and annotated  was our serious omission of any
description of the Nullstellensatz in a Weyl Algebra.

\subsection{Sums of squares modulo cyclic equivalence} A still
open, important conjecture in the classification theory of von
Neumann algebras was recently reduced by F. Radulescu to an
asymptotic Positivstellensatz in the free algebra. We reproduce
from his preprint \cite{Radulescu} the main result. We do not
explain below the standard terminology related to von Neumann
algebras, see for instance \cite{Takesaki}.

The following conjecture was proposed thirty years
ago in \cite{Connes}:\index{Connes conjecture}\\

{\it Every type $II_1$ factor can be embedded into an ultraproduct
of the hyperfinite factor.}\\

There are presently quite a few reformulations or reductions of
this conjecture. The one of interest for this survey can be
formulated as follows.

Let $F = \mathbb C \langle x_1,..., x_g \rangle$ be the free
algebra with anti-linear involution $x_j^\ast = x_j, \ 1 \leq j
\leq g$. We complete $F$ to the algebra of convergent series
$$ \hat{F} = \{ \sum_w a_w w; \ \sum_w |a_w| r^{|w|} <\infty, \
\forall r >0 \},$$ where $w$ runs over all words in $F$ and $a_w
\in \mathbb C$. The resulting Fr\'echet space $\hat{F}$ carries a
natural weak topology denoted $\sigma(\hat{F}, \hat{F}^\ast)$.

A trace $\tau$ in a von-Neumann algebra $M$
is a linear functional which has by definition the
cyclic invariant property $\tau (a_1... a_n) = \tau(a_2 a_3...a_n
a_1)$. Two series $f_1, f_2 \in \hat{F}$ are called {\it
cyclically equivalent} if $f_1-f_2$ is the weak limit of a linear
combination of elements $w-w'$, where $w \in F$ is a word and $w'$
is a cyclic permutation of it.

The following asymptotic Positivstellensatz holds.

\begin{thm}
\cite{Radulescu} \index{Radulescu Theorem} Let $f \in \hat{F}$ be
a symmetric series with the property that for every separable,
type $II_1$ von Neumann algebra $(M,\tau)$ and every $g$-tuple of
self-adjoint elements $X$ of $M$ we have $\tau(f(X)) \geq 0$. Then
$f$ is cyclically equivalent to a weak limit of sums of squares
$s_n$,  $s_n \in \Sigma^2F.$
\end{thm}

It is not known whether one can replace the test $II_1$ algebras
by finite dimensional algebras, but an answer to this querry would
solve Connes conjecture.

\begin{cor} Connes embedding conjecture holds if and only if for
every symmetric element $f \in \hat{F}$ the following assertion
holds:

f is cyclically equivalent to a weak limit of sums of squares
$s_n$, $s_n \in \Sigma^2F,$ if and only if for any positive integer
$d$ and $g$-tuple of self-adjoint $d \times d$ matrices $X$ one
has ${\rm trace} f(X) \geq 0$.
\end{cor}

The proofs of Radulescu's theorem and the corollary follow the
same pattern we are by now familiar with: a convex separation
argument followed by a GNS construction. See for details
\cite{Radulescu}, and for a last minute refinement \cite{KS06}.

\section{Convexity in a free algebra}
\label{sec:convexity}

Convexity of functions,  domains
and their close relative, positive curvature of varieties, are
very natural notions  in a $\ast$-free algebra.
A shocking thing happens: these convex functions are
so rare as to be almost trivial.
This section illustrates a simple  case,
that of convex polynomials, and we see how
in a free algebra the  Nichtnegativtellens\"atze
have extremely strong consequences for
inequalities on  derivatives.
The phenomenon
 has direct qualitative consequences
for systems engineering as we see in \S \ref{sec:engLMI}.
The results of this section can be read independently
of all but a few definitions in \S \ref{sec:nonCom},
and the proofs require only a light reading of it.

This time $ \mathbb R \langle x \rangle$ denotes the free
$\ast$-algebra in indeterminates $x = (x_1,...,x_g)$, over the
real field. There is an involution $x_j^\ast = x_j$ which reverses
the order of multiplication $ (f p)^\ast = p^\ast f ^\ast.$
In this
exposition we take symmetric variables $x_j=x_j^*$, but in the
literature we are summarizing typically $x_j$ can be taken either
free or symmetric with no change in the conclusion, for example,
the results also hold for symmetric polynomials in $\RRxxS$.

A symmetric polynomial $p, p^\ast = p,$ is {\it matrix convex}
\index{convex polynomial} if for each positive integer $n$, each
pair of tuples $X=(X_1,\dots,X_g)$ and $Y=(Y_1,\dots,Y_g)$ of
symmetric $n\times n$ matrices, and each $0\le t \le 1$,
\begin{equation}
\label{eqn:matrixconvex} p(tX+(1-t)Y)\leq tp(X)+(1-t)p(Y).
\end{equation}
Even in one-variable, convexity in the noncommutative setting
differs from convexity in the commuting case because here $Y$ need
not commute with $X$. For example, to see  that the polynomial
$p=x^4$ is not matrix convex, let
$$
X=\left ( \begin{array}{cc} 4&2\\2&2\end{array}\right ) \mbox{ and
} Y=\left ( \begin{array}{cc} 2&0\\0&0\end{array}\right )
$$
and compute
$$
\frac12 X^4+\frac12 Y^4 - (\frac12 X+\frac12 Y)^4 = \left (
\begin{array}{cc}164&120\\120&84\end{array}\right )
$$
which is not positive semi-definite. On the other hand,  to verify
that $x^2$ is a matrix convex polynomial, observe that
\begin{eqnarray}
\nonumber
tX^2+(1-t)Y^2&-&(tX+(1-t)Y)^2\\
\nonumber &=&t(1-t)(X^2-XY-YX+Y^2)=t(1-t)(X-Y)^2 \geq 0.
\end{eqnarray}

\index{convex degree 2 Theorem}
\begin{thm}
\label{thm:matrixconvexI}\cite{HM04conv}
   Every convex symmetric polynomial
   in the free algebra $\RRx$ or $\RRxxS$
    has degree two or less.
\end{thm}

As we shall see convexity of $p$ is equivalent to its
``second directional derivative" being a positive polynomial.
As a matter of fact, the  phenomenon has nothing to do
with  order two derivatives and
the extension of this to  polynomials with $k^{th}$ derivative
nonnegative is given later in Theorem \ref{thm:matrixderivI}.

Yet stronger about convexity  is the next local implies global theorem.

Let $\cP$ denote a collection of  symmetric  polynomials in
non-commutative variables $x=\{x_1, \cdots , x_g\}$.  Define the
matrix nonnegativity domain $\POD{\cP}$
associated to $\mathcal P$
\index{ $\POD{\cP}$ matrix nonnegativity domain}
to be the set of tuples $X=(X_1, \cdots , X_g)$
\index{ $X = (X_1, \cdots , X_g)$} of
finite dimensional {\it real} matrices of all sizes, except 0
dimensions,
%
%
making $p(X_1, \cdots ,X_g)$
a positive semi-definite matrix.

\index{convex polynomial, local degree 2 Theorem}
\begin{thm}
\label{thm:matrixconvexII}\cite{HM04conv}
Suppose there is a set $\cP$
of symmetric polynomials, whose matrix nonnegativity domain
$\POD{\cP}$  contains open sets in all large enough dimensions.
Then
   every symmetric polynomial $p$ in $\RRx$ or in
   $\mathbb R \langle x, x^\ast
\rangle$
   which is matrix convex on $\cD_\cP$ has degree two or less.
\end{thm}

The first convexity theorem follows from Theorem \ref{free SOS},
and we outline below the main ideas in its proof. The proof of the
more general, order $k$ derivative,
is similar and we will return to it later in this section.
The
proof of Theorem \ref{thm:matrixconvexII} requires different
machinery (like that behind representation (\ref{eq:repM} ))
and is not presented here.

At this point we  describe a bit of history.
In the beginning was
Karl L\"owner who studied a class of real analytic
functions in one real variable called matrix monotone,
which we shall not  define  here.
L\"owner gave  integral representations
and these have developed beautifully over the years.
The impact on our story comes
a few years later when L\"owner's student
Klaus \cite{K36} introduced matrix convex functions
$f$ in one variable.
Such a function $f$ on $[0,\infty] \subset \bbR$
 can be represented as $f(t) = tg(t)$  with $g$ matrix monotone,
 so the representations for $g$ produce representations  for $f$.
Modern references are \cite{OSTprept}, \cite{U02}.
Frank Hansen has extensive deep work
on matrix convex an monotone functions
whose  definition  in several variables is different than the
one we use here, see\cite{HanT06};
for a recent reference see  \cite{Han97}.

For a polynomial $p \in \RRx$ define the {\it directional
derivative}\index{directional derivative}:
$$p'(x)[h]= \frac{d}{dt} p(x+th)_{|_{t = 0}}. $$
It is a linear form in $h$.
Similarly, the $k^{th}$ derivative
$$p^{(k)}(x)[h]= \frac{d^k}{dt^k} p(x+th)_{|_{t = 0}} $$
is homogeneous of degree $k$ in $h$.

More formally, we regard
  the directional derivative $p^\prime(x)[h] \in \mathbb R\langle x,h
\rangle$ as a
       polynomial in $2g$ free  symmetric
       (i.e. invariant  under $^*$) variables
      $(x_1,\dots,x_g,h_1,\dots,h_g)$;
      In the case of
     a word
     $w=x_{j_1}x_{j_2}\cdots x_{j_n}$ the derivative is:
    \begin{equation*}
      w'[h] = h_{j_1}x_{j_2}\cdots x_{j_n}
       + x_{j_1}h_{j_2}x_{j_3}\cdots x_{j_n} + \ \dots \
       + x_{j_1}\cdots x_{j_{n-1}}h_{j_n}
    \end{equation*}
      and for a polynomial
      $p=       p^\prime(x)[h]=\sum p_w w$
      the derivative is
    \begin{equation*}
      p^\prime(x)[h]=\sum p_w w'[h].
    \end{equation*}
      If $p$ is symmetric, then so is $p^\prime$.

     For $g$-tuples of symmetric matrices of a fixed size
     $X,H ,$ observe that the evaluation formula
   \begin{equation*}
     p^\prime(X)[H]=\lim_{t\to 0} \frac{p(X+tH)-p(X)}{t}
   \end{equation*}
     holds. Alternately, with $q(t)=p(X+tH)$, we find.
   \begin{equation*}
     p^\prime(X)[H]=q^\prime(0).
   \end{equation*}

Likewise for a polynomial $p \in \RRx$, the {\it Hessian}
$p^{\prime\prime}(x)[h]$ of $p(x)$ can be thought of as the formal
second directional derivative of $p$ in the ``direction'' $h$.
Equivalently, the Hessian of $p(x)$ can also be defined as the
part of the polynomial
$$r(x)[h]:=p(x+h)-p(x)$$
in the free algebra in the symmetric variables that is homogeneous
of degree two in $h$.

If $p^{\prime\prime} \neq 0$, that  is,
  if $\textrm{degree}\, p \ge 2$, then
the degree of $p^{\prime\prime}(x)[h]$ as a polynomial in the $2g$
variables $x_1,\ldots,x_g, h_1\ldots,h_g$ is equal to the degree
of $p(x)$ as a polynomial in $x_1,\ldots,x_g$.

Likewise for $k^{th}$ derivatives.

\begin{example}

1.\ \ $p(x)= x_2 x_1 x_2$
$$p'(x)[h]= \frac{d}{dt}
[ (x_2+th_2) (x_1 +th_1) (x_2 +h_2) ]_{|_{t = 0}}= h_2 x_1 x_2 +
x_2 h_1 x_2 + x_2 x_1 h_2. $$

2. \ \  One variable $p(x)=x^4$. Then
$$p'(x)[h]= hxxx + xhxx  + xxhx + xxxh$$
Note each term is linear in $h$ and $h$ replaces each occurrence
of $x$ once and only once:
$$p''(x)[h]= $$
$$hhxx + hhxx  + hxhx + hxxh+$$
$$hxhx + xhhx  + xhhx + xhxh+$$
$$hxxh + xhxh  + xxhh + xxhh,$$
which yields
$$ p''(x)[h]=
2hhxx   + 2hxhx + 2 hxxh
  + 2xhhx   + 2 xhxh  + 2 xxhh. $$
Note each term is degree two  in $h$ and $h$ replaces each pair of
$x$'s exactly once.
Likewise
$$
p^{(3)}(x)[h]=
 6 (hhhx + hhxh + hxhh +xhhh)
$$
and
$
p^{(4)}(x)[h]=
 24 hhhh
$
and $p^{(5)}(x)[h]=0$.

\bigskip

3. \  \  $p=x_1^2 x_2$
$$
p^{\prime\prime}(x)[h]= h_1^2 x_2 + h_1 x_1 h_2 + x_1 h_1 h_2.
$$
\end{example}

The definition of a convex polynomial can be easily adapted to
domains. Then one remarks without difficulty that,
in exact analogy with
the commutative case, a polynomial $p$ is convex (in a domain) if
and only if the Hessian evaluated at the respective points is
non-negative definite.
Because of this
Theorem \ref{thm:matrixconvexI}
is an immediate consequence of the  next theorem
restricted to $k=2$.

\begin{thm}
\label{thm:matrixderivI}
  Every  symmetric polynomial $p$
   in the free algebra $\RRx$ or $\RRxxS$
   whose $k^{th}$ derivative is a matrix positive polynomial
    has degree $k$ or less.
\end{thm}

\proof \ (when the variables $x_j$ are symmetric).

Assume $p^{(k)} (x)[h]$ is a matrix positive polynomial, so that,
in view of Theorem \ref{free SOS} we can write
it
as a sum of squares:
$$p^{(k)} (x)[h] = \sum f_j^\ast f_j; $$
here each $f_j(x,h)$ is a polynomial in the free algebra
$\mathbb R \langle x,h \rangle$.\\

If $p^{(k)} (x)[h]$ is identically equal to zero, then the statement follows.
Assume the contrary, so that $p^{(k)} (x)[h]$ is homogeneous
of degree $k$ in $h$, and there are tuples of matrices $X,H$
and a vector $\xi$ in the underlying finite dimensional Hilbert space,
so that
$$\langle  p^{(k)} (X)[H]\xi, \xi \rangle >0.$$
By multiplying $H$ by a real scalar $t$ we find
$$ t^k \langle  p^{(k)} (X)[H]\xi, \xi \rangle = \langle  p^{(k)} (X)[tH]\xi, \xi \rangle >0,$$
whence $k = 2\mu$ is an even integer.

Since in a sum of squares the highest degree terms cannot cancel,
the degree of each $f_j$ is at most $\nu$ in $x$ and $\mu$ in $h$,
where $2\nu$ is the degree of $\pk$ in $x$.

Since $\pk$ is a directional derivative,
it must have a highest degree
term of the form $h_{i_1} \cdots h_{i_k} m(x)$
where the monomial $m(x)$ has degree
equal to degree $\pk - k$; also
 $h_{i_j}$ is allowed to equal $h_{i_\ell}$.
Thus some
product, denote it $f_J^\ast f_J$, must contain such a term.
(Note the the order of the $h's$ vs. the $x's$ matters.)
This
forces $f_J$ to have the form
$$ f_J= c_1 ( h_{i_{\mu \; +1}} \cdots h_{i_k} ) m(x)
\ + \ c_2 (h_{i_1} \cdots h_{i_{\mu} }) \ + \ ...\  ,$$
the $c_j$ being scalars.

To finish the proof
use
that $f_J^\ast f_J$ contains
$$c^2 \; m(x)^\ast
(h_{i_{\mu \; +1}} \cdots h_{i_k})^*
(h_{i_{\mu \; +1}} \cdots h_{i_k} )  m(x)$$
and
this can not be cancelled  out,
  so
  $$\deg  \pk \
  = \
    k + 2 (\deg \pk -k)
  =   2 \deg \pk - k .$$
Solve this to find    $ \deg \pk =k$.
  Thus $p$ has degree k. \qed \\

We use a previous example in order to illustrate this proof
when $k=2$.
\begin{example}
{\it Example $p=x^4$ is not matrix convex; here $x=x^\ast$.}
\\
Calculate that
$$ p''(x)[h]= 2hhxx + 2 hxhx + 2 hxxh
            + 2 xhhx + 2 xhxh + 2 xxhh.
$$
Up to positive constants
some polynomial $f_J^\ast f_J$ contains a term $ hhxx$, so
  $f_J=   hxx + h +  \ldots $.

So $f_J^\ast f_J$ contains $xxhhxx$. This is a highest order
perfect square so can be cancelled out. Thus is appears in $p''$,
which as a consequence has degree 6. This a contradiction.
\end{example}

We call the readers attention to work which goes beyond what we
have done in several directions. One \cite{HMVprept} concerns a
noncommutative rational function $r$ and  characterizes those
which are convex near 0.
\index{convex rational classification Theorem}
\ It is an extremely small and rigidly
behaved class, for example, $r$ is convex on the entire component
of the "domain of $r$" which contains 0. This rigidity is in
analogy to convex polynomials on some "open set" having degree 2
or less and this implying they are convex everywhere. Another direction
is the classification of noncommutative polynomials whose Hessian
$p''(x)[h]$ at most $k$ "negative noncommutative eigenvalues" In
\cite{DHMprept} it is shown that this implies
$$ \deg  p \ \ \leq 2k +2.$$
Of course the special case we studied in this section is exactly
that of polynomials with $k=0$.

\section{Dimension free engineering: LMIs vs. CMIs}
\label{sec:engLMI}

This section demonstrates the need for real algebraic geometry (in
the broad sense) aimed at convexity over a free or nearly free
$*$- algebra. From this viewpoint the theory in this survey goes
in an essential direction but much more is called for in order to
do general engineering problems. Hopefully the brief description
in this section will give a little feel for where we now stand. We
are aiming at one of the major issues in linear systems theory:\\

{\it Which problems convert to a convex matrix inequality, CMI?
How does one do the  conversion?}\\

To be in line with the engineering literature, we use below a
slightly different notation than the rest of the article. For
instance $A^T$ denotes the transpose of a (real entries) matrix,
and $A^T$ replaces in this case the involution $A^\ast$ we have
encountered in the previous sections. The inner product of vectors
in a real Hilbert space will be denoted $u \cdot v$.

\subsection{Linear systems}
A {\it linear system}\index{linear system} $\mathfrak F$ is given
by the linear differential equations
$$\frac{dx}{dt}=Ax + Bu$$
$$y=Cx$$
with the vector
\begin{itemize}
\item
  $x(t)$ at each time $t$
being in the  vector space $\cX$ called the {\it state space},
\index{state space}
\item $u(t)$ at each time $t$ being in the
vector space $\cU$ called the {\it input space}, \index{input
space}
\item $y(t)$ at each time $t$ being in the vector space
$\cY$ called the {\it output space}\index{output space},
\end{itemize}
and $A,B,C$ being linear maps on the corresponding vector spaces.

\subsection{Connecting linear systems}
Systems can be connected in incredibly complicated configurations.
We  describe a simple connection and this goes along way toward
illustrating the  general idea. Given two linear systems
$\mathfrak F$, $\mathfrak G$, we describe the formulas for
connecting them as follows.
{\flushleft
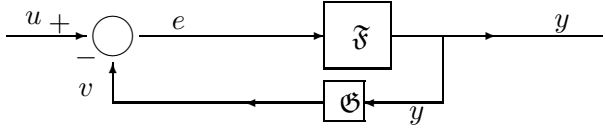
{\begin{picture}(50,80)(10,20)
\put(50,70){\vector(1,0){30}} \put(60,77){\makebox(0,0){{$u$}}}
\put(90,70){\circle{15}} \put(70,75){\makebox(0,0){{$+$}}}
\put(80,62){\makebox(0,0){{$-$}}} \put(100,70){\vector(1,0){70}}
\put(80,50){\makebox(0,0){{$v$}}}
\put(115,75){\makebox(0,0){{$e$}}}
\put(170,58){\framebox(25,25){}}
\put(184,70){\makebox(0,0){{$\mathfrak F$}}}
\put(195,70){\vector(1,0){40}} \put(235,70){\line(1,0){40}}
\put(260,75){\makebox(0,0){{$y$}}} \put(90,45){\vector(0,1){15}}
\put(90,45){\line(1,0){50}} \put(170,45){\vector(-1,0){30}}
\put(170,38){\framebox(15,15){}}
\put(180,45){\makebox(0,0){{$\mathfrak G$}}}
\put(215,45){\vector(-1,0){30}} \put(215,45){\line(0,1){25}}
\put(205,40){\makebox(0,0){{$y$}}}
\end{picture}}}

Systems $\mathfrak F$ and $\mathfrak G$ are respectively given by
the linear differential equations
$$\frac{dx}{dt}=Ax+Be,\qquad  \qquad   \frac{d\xi}{dt} =a\xi+bw,$$
$$y=Cx,\qquad  \qquad  v=c\xi.$$
The connection diagram is equivalent to the algebraic statements
$$w=y \ \ \ \ {\rm and} \ \ \ \ \ e=u-v.$$
The {\it closed loop system} \index{closed loop system} is a new
system whose differential equations are
$$\frac{dx}{dt}=Ax  - Bc\xi + Bu,$$
$$\frac{d\xi}{dt}=a\xi+by=a\xi + bCx,$$
$$y=Cx.$$
In matrix form this is
\begin{equation}
\label{eq:cloopsys} \frac{d}{dt}{\left(\begin{array}{c} x\\ \xi
\end{array}\right)} = \left(\begin{array}{cc} A& -Bc\\bC & a
\end{array}\right)\left(\begin{array}{c} x\\ \xi
\end{array}\right)+\left(\begin{array}{c} B\\0
\end{array}\right)u,
\end{equation}
$$
y=(C\   0)\left(\begin{array}{c} x\\ \xi   \end{array}\right),
$$
where the state space of the closed loop systems is the direct sum
`$\cX \oplus \cY$' of the state spaces $\cX$ of $\mathfrak F$ and
$\cY$ of $\mathfrak G$.
The moral of the story is:\\

{\it System connections produce a new system whose coefficients
are matrices with entries which are polynomials in the
coefficients of the component systems. }\\

 Complicated signal
flow diagrams give complicated matrices of polynomials. Note in
what was said the dimensions of vector spaces and matrices never
entered explicitly; the algebraic form of (\ref{eq:cloopsys} ) is
completely determined by the flow diagram. We have coined the term
{\it dimension free}\index{dimension free} for such problems.

\subsection{Energy dissipation}
We have a system $\mathfrak F$ and want a condition which checks
whether
$$
\int_{0}^{\infty}{|u|}^{2}dt \geq \int_{0}^{\infty}{|\mathfrak
Fu|}^{2}dt,
    \qquad x(0)=0,
$$
holds for all input functions $u$, where $\mathfrak F u =y$ in the
above notation . If this holds $\mathfrak F$ is called a { \it
dissipative system} \index{dissipative system} {\flushleft
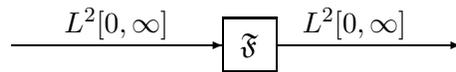
\hspace{1in}

\begin{picture}(100,70)(10,20)
\put(130,55){$L^{2}[0, \infty$]} \put(110,50){\vector(1,0){80}}
\put(190,40){\framebox(20,20){}} \put(220,55){$L^{2}[0, \infty$]}
\put(210,50){\vector(1,-0){70}}
\put(200,50){\makebox(0,0){{$\mathfrak F$}}}
\end{picture}}


This is  analysis but it converts to algebra because of the
following construction. Hope there is a "potential energy" like
function $V \geq 0, \ V(0)=0,$ on the state space;  it should
satisfy:\\

{\it
  potential energy now + energy in $\geq$
potential energy then
+  energy out.}\\

\noindent
  In mathematical notation this is
  $$
  V(x(t_{1})) + \int_{t_{1}}^{t_{2}} |u(t)|^{2} \ \ \geq \ \
   V(x(t_{2}))  + \int_{t_{1}}^{t_{2}}{|y(t)|^{2}}
   $$
and a $V \geq 0, \ V(0)=0,$
  which satisfies this for all input functions $u$ and
  initial states $x(t_1)$ is called
  a {\it storage function}\index{storage function}.
We can manipulate this integral condition  to obtain first a
differential inequality and then an algebraic inequality, as
follows:
$$
0\geq \frac{V(x(t_{2}))-V(x(t_{1}))}{t_{2}-t_{1}} \ \ + \ \
\frac{1}{t_{2}-t_{1}}\int_{t_{1}}^{t_{2}}|y(t)|^{2}-|u(t)|^{2},
$$
  $$
  0 \geq \nabla V(x(t_{1}) ) \cdot \frac{dx}{dt}(t_{1})
  \ \  + \ \ |y(t_{1})|^{2}  - |u(t_{1})|^{2}.
  $$
  Use $\frac{dx}{dt}(t_{1}) = \ Ax(t_{1})+Bu(t_{1})$
to get
  $$
0 \ \geq \ \ \nabla  V(x(t_{1}))\cdot (Ax(t_{1})+Bu(t_{1})) \ \ +
\ \ |Cx((t_{1}))|^{2} -  |u(t_{1})|^2.
$$
The system is dissipative if and only if this holds for all
$u(t_{1})$, $x(t_{1})$ which can  occur when it runs (starting at
$x(0)=0$). All vectors  $u(t_{1})$ in $\cU$ can certainly occur as
an input and if all  $x(t_{1})$ can occur we call the system {\it
reachable} \index{reachable system}. Denote $x(t_{1})$ by $x$ and
$u(t_{1})$ by $u$
  \begin{equation}
    \label{eq:dissipux}
0 \ \geq \ \ \nabla  V(x) \cdot (Ax  + Bu) \ \ + \ \ |Cx|^{2} -
|u|^2,
  \end{equation}
and conclude:
\begin{thm}
   The system $A,B,C$
    is dissipative if inequality (\ref{eq:dissipux})
    holds for all $u \in \cU, x \in \cX$.
    Conversely, if $A,B,C$ is reachable,
    then dissipativity implies inequality (\ref{eq:dissipux})
    holds for all $u \in \cU, x \in \cX$.
\end{thm}
For a linear system we try  $V$ which is quadratic, so $V(x)=Px
\cdot x$ with $P \geq 0$ and $\nabla V(x)= 2 Px$. At this point
there are two commonly pursued paths which constitute the next two
subsections.

\subsubsection{Riccati inequalities}
Use that $\nabla V(x)= 2 Wx$ in (\ref{eq:dissipux}) to get
$$0\geq 2 Wx\cdot (Ax+Bu) + |Cx|^{2} - |u|^{2}, \ \
\ \ \ \ {\rm for \ all }\ u,x,$$ so
\begin{equation}
\label{eq:lindisux} 0\geq \max_u \big( [WA+A^{T}W+C^{T}C]x \cdot x
+ 2B^{T}Wx \cdot u- |u|^{2} \big).
\end{equation}
The maximizer in $u$ is $u = B^T Wx$, hence
$$
0 \geq 2Wx \cdot Ax + 2{|B^T W x |}^{2} + {|Cx|}^{2} - {|B^T W x
|}^{2}.
$$
Which in turn is
$$0 \geq [WA + A^{T}W + W BB^{T} W + C^{T}C]x\cdot x .$$
This is  the classical {\it Riccati matrix inequality}
\index{Riccati matrix inequality}
$$0 \geq \ WA + A^{T}W + WBB^{T}W + C^{T}C $$
which together with $W \geq 0$ insures dissipativity and is also
necessary for it when the system is reachable.

\subsubsection{Linear Matrix Inequalities (LMI)}
Alternatively we do not need to compute $\max_{u}$ but can express
(\ref{eq:lindisux}) as  the inequality:
$$L(W):
=\left(\begin{array}{cc} WA+A^{T}W+C^{T}C & WB\\ B^{T}W & -I
\end{array}\right)\left(\begin{array}{c} x\\ u
\end{array}\right)
\cdot
  \left(\begin{array}{c} x\\ u \end{array}\right)
  \leq 0
$$
for all $u \in \cU, x \in \cX$. That is the following matrix
inequality holds:
$$\left(\begin{array}{cc} WA+A^{T}W+C^{T}C & WB\\
B^{T}W & -I   \end{array}\right) \leq 0.$$ Here $A$, $B$, $C$
describe the system and $W$ is an unknown matrix. If the system is
reachable, then $A$, $B$, $C$ is dissipative if and only if $L(W)
\leq 0$ and $W \geq 0$.

Recall that the {\it Schur complement}\index{Schur complement} of
a matrix is defined by
$$SchurComp \; \left(\begin{array}{cc}
\alpha & \beta \\ \beta^{T} & \gamma   \end{array}\right): =
\alpha -\beta \gamma^{-1}\beta^{T}.$$ Suppose $\gamma $ is
invertible. The matrix is positive semi-definite if and only if
$\gamma > 0$ and its Schur complement is positive semi-definite.
Note that
$$
SchurComp \; L(W) \ = \ WA + A^{T}W + WBB^{T}W + C^{T}C
$$
featuring the Riccati inequality we saw before. Indeed, $L(W) \leq
0$ if and only if $WA + A^{T}W + WBB^{T}W + C^{T}C \leq 0$, since
this implies $WA+A^{T}W+C^{T}C \leq 0 $.
Thus the Riccati approach
and the LMI approach give equivalent answers.
\subsection{Example: An $H^\infty$ control problem}
\label{sec:HinftyProb}
 Here is a basic engineering problem:\\

{\it Make a given system dissipative by designing a feedback
law.}\\

To be more specific, we are given a signal flow diagram:
{\flushleft
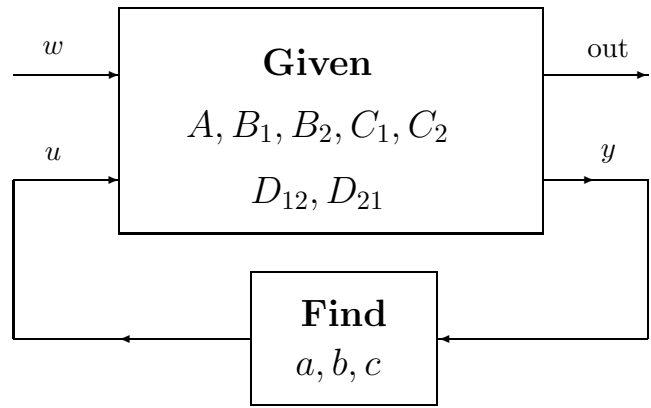
\begin{picture}(100,200)(10,20)
\put(150,100){\framebox(160,85){}} \put( 110,160){\vector( 1, 0){
40}} \put( 110,120){\vector( 1, 0){ 40}}
\put(225,165){\makebox(0,0){{\Large {\bf Given}}}}
\put(225,140){\makebox(0,0){{\Large${A, B_1, B_2, C_1, C_2}$}}}
\put(225,115){\makebox(0,0){{\Large$ {D_{12}, D_{21}}$}}}
\put(310,160){\vector( 1, 0){ 40}} \put(310,120){\vector( 1, 0){
20}} \put(330,120){\line( 1, 0){ 20}}
\put(125,170){\makebox(0,0){{$w$}}}
\put(335,170){\makebox(0,0){{out}}}
\put(125,130){\makebox(0,0){{$u$}}} \put(
335,130){\makebox(0,0){{$y$}}} \put(200,35){\framebox(70,50){}}
\put(110,60){\line(0,1){60}} \put(350,60){\line(0,1){60}}
\put(200,60){\vector(-1, 0){50}} \put(110,60){\line(1, 0){40}}
\put(350,60){\vector(-1,0){79}} \put(
235,70){\makebox(0,0){{\Large {\bf Find}}}} \put(
235,50){\makebox(0,0){{\Large  $a,b,c$ }}}
\end{picture}}

where the given system is {\begin{eqnarray*}
\frac{ds}{dt}&=& As+B_1w +B_2u,\\
{\rm out}&=&C_1s+D_{12}u+D_{11}w,\phantom{l^l}\\
y&=& C_2s+D_{21}w,
\end{eqnarray*}}
$$
  D_{21}=I,\qquad D_{12}D^\prime_{12} = I, \qquad
D^\prime_{12}D_{12}=I, \qquad D_{11}=0.$$
  We want to find an
unknown  system
$$\frac{d\xi}{dt} =  a \xi +  b,
\qquad\qquad u= c \xi,$$ called the {\it
controller}\index{controller},
  which makes
the system dissipative over every finite horizon. Namely:
$$
\int^T_0|w(t)|^2dt \geq \int\limits^T_0 |out(t)|^2 dt,\ \ \  \
s(0)= 0.
$$
So $ { a, b , c}$ are the critical unknowns.

\subsubsection{Conversion to algebra}
The dynamics of the ``closed loop" system has the form
  $$\frac{ d}{dt}
  \left(%
\begin{array}{c}
   s \\
   \xi \\
\end{array}%
\right)
  =  \cA \left(%
\begin{array}{c}
   s \\
   \xi \\
\end{array}%
\right) + \cB w
$$
$$
  out = \cC  \left(%
\begin{array}{c}
   s \\
   \xi \\
\end{array}%
\right) + \cD w
$$
where $\cA, \cB, \cC,\cD$ are "$ 2 \times 2$ block matrices" whose
entries are polynomials in the  $A's, B's, \cdots,  a, b,c$ etc.
The storage function inequality which corresponds to energy
dissipation has the form
\begin{equation}
   \label{ineq:Hleq0}
  H:
={\cA}^T E + E {\cA} + E  {\cB}  {\cB}^T E +   {\cC}^T {\cC} \leq
0
\end{equation}
where $E$ has the form
$$
E = \left(\begin{array}{cc} E_{11}& E_{12}\\ E_{21}&
E_{22}\end{array}\right),\ \ \ \ \ \ \ \ {\ } E_{12}=E_{21}^T.
$$

The algebra problem above in more detail is to solve inequality
(\ref{ineq:Hleq0})
$$
H = \left(\begin{array}{cc} H_{ss}& H_{sy}\\ H_{ys}&
E_{yy}\end{array}\right) \leq 0, \ \ \ \ \ \ \ \ \ \ \ \ \ \ {\ }
H_{sy}=H^T_{ys},
$$
where the entries of $H$ are the polynomials:\bigskip

$$H_{ss} = E_{11} \, A + A^T \, E_{11} + C_1^T \, C_1 + E_{12}^T \,
   b \, C_2 + C_2^T \, b^T \, E_{12}^T
    + E_{11} \, B_1 \, b^T \, E_{12}^T +$$ $$
    E_{11} \, B_1 \, B_1^T \, E_{11} +
E_{12} \, b \, b^T \, E_{12}^T +
    E_{12} \, b \, B_1^T \, E_{11},$$
$$H_{sz}=E_{21} \, A + \frac{a^T \, (E_{21} + E_{12}^T)}{2} +
   c^T \, C_1 + E_{22} \, b \, C_2 +
    c^T \, B_2^T \,
   E_{11}^T +$$ $$
\frac{E_{21} \, B_1 \, b^T\, (E_{21} + E_{12}^T)}{2} +
    E_{21} \,
   B_1 \, B_1^T \, E_{11}^T +
\frac{E_{22} \, b \, b^T \, (E_{21} + E_{12}^T)}{2} +
    E_{22} \, b \, B_1^T \, E_{11}^T, $$
$$H_{zs}=A^T \, E_{21}^T + C_1^T \, c + \frac{(E_{12} + E_{21}^T)
\, a}{2} +
   E_{11} \,
B_2 \, c +
    C_2^T \, b^T \, E_{22}^T +
   E_{11} \, B_1 \, b^T \, E_{22}^T +$$ $$
   E_{11} \, B_1 \, B_1^T \, E_{21}^T +
    \frac{(E_{12} + E_{21}^T) \, b \, b^T \, E_{22}^T}{2} +
    \frac{(E_{12} + E_{21}^T) \, b \, B_1^T \, E_{21}^T}{2},$$
$$H_{zz}=E_{22} \, a + a^T \, E_{22}^T + c^T \, c + E_{21} \, B_2
\, c +
    c^T \, B_2^T \, E_{21}^T +
E_{21} \, B_1 \, b^T \, E_{22}^T + $$ $$E_{21} \, B_1 \, B_1^T \,
E_{21}^T +
   E_{22} \, b \, b^T \, E_{22}^T + E_{22} \, b \, B_1^T \,
E_{21}^T. $$
   \bigskip

Here $A$, $B_1$, $B_2$, $C_1$, $C_2$ are known and the unknowns
are $a$, $b$, $c$ and  for $E_{11}$, $E_{12}$, $E_{21}$ and
$E_{22}$.

  We very much wish that these inequalities are convex in
the unknowns (so that numerical solutions will be reliable). But
our key inequality above is not convex in the unknowns.

\subsubsection{The key question}
\label{sec:convQ} Is
  there is a set of noncommutative {\it convex} inequalities
with an equivalent set of solutions? \\

This is a question in  algebra not in numerics and the  answer
after a lot of work is yes. The path to success is: {\it
\begin{enumerate}
   \item
   \label{convQ:chgV}
Firstly, one must eliminate unknowns, change variables and get a
new set of inequalities $\cE$.
   \item
   \label{convQ:conv}
   Secondly, one must check that $\cE$
   is ``convex" in the unknowns.
\end{enumerate}
}

  This outline transcends our example and applies to very many
situations. The second issue of this is becoming reasonably
understood, for as we saw earlier,  a convex polynomial with real
coefficients has degree two or less, so these are trivial to
identify. While the level of generality of the theory we have
presented in this paper is less than we now require, to wit,
polynomials with indeterminates
  as coefficients and matrices with  polynomial entries;
this does not add radically different structure, see discussion in
\S \ref{sec:symbcoef}. The first issue,  changing  variables,
  is still a collection  of isolated tricks,
with which mathematical theory has not caught up. For the
particular problem in our example we shall not derive the solution
since it is long. However, we do state the classical answer in the
next subsection.

\subsubsection{Solution to the Problem}
The textbook solution is as follows, due to Doyle-Glover-
Kargonekar-Francis. It appeared in \cite{DGKF89} which won the
1991 annual prize for the best paper to appear in an IEEE journal.
Roughly speaking it was deemed the best paper in electrical
engineering in that year.

We denote
$$
DGKF_X  :=  (A-B_2 C_1)'X  + X ( A-B_2 C_1)
  + X (\gamma^{-2}B_1B_1' - B_2^{-1}B_2')X $$
$$
DGKF_Y  :=  A^\times Y + Y{A^\times}' + Y (\gamma^{-2}C_1'C_1 -
C_2'C_2)Y,
$$
where $A^\times:=A-B_1C_2$.

\begin{thm}  \cite{DGKF89}
There is a system  solving the control problem if there exist
solutions
$$X \geq 0 \ \ \ \ \ \ \
and \ \ \ \ \ \ \ Y > 0$$ to inequalities the
$$ DGKF_Y \leq 0 \ \ \mbox{and} \ \  DGKF_X \leq  0$$
which satisfy the {coupling condition}
$$X - Y^{-1} < 0.$$
This is if and only if provided $Y > 0$ is replaced by
  $Y \geq 0 $ and $Y^{-1}$
is interpreted correctly.
\end{thm}

This set of inequalities while not usually convex in $X,Y$ are
convex in the new variables
  $W = X^{-1}$ and $Z = Y^{-1}$, since
  $DGKF_X $ and $DGKF_Y$ are linear in them and
  $X - Y^{-1}= W^{-1} - Z$ has second derivative
  $ 2 W^{-1} H W^{-1}H W^{-1} $
which is non negative in $H$ for each $ W^{-1}= X > 0$.
These
inequalities are also equivalent to LMIs which we do not write
down.

\subsubsection{Numerics and symbolics}
 A single Riccati
inequality is much more special than an LMI and numerical solvers
for Riccatis are faster and handle bigger matrices.
This survey obviously has not aimed at numerics, but at  algebraic
precursors to using numerics.

The mathematics here aims toward
helping an engineer who writes a
toolbox which other engineers will use for designing systems,
like control systems.
What goes in such toolboxes is algebraic formulas
like the DGKF inequalities above with matrices $A,B, C$
unspecified and reliable numerics for solving them when a user
does specify $A,B, C$ as matrices.
A user who designs a controller
for a helicopter puts in the mathematical systems model for his
helicopter and puts in matrices, for example, $A$ is a particular
$R^{8 \times 8}$ matrix etc.
Another user who designs a satellite
controller might have a 50 dimensional state space and of course
would pick completely different $A,B,C$.
Essentially any matrices
of any compatible dimensions can occur and our claim
that our algebraic
formulas are convex in the ranges we specify must be true.

The toolbox designer faces two completely different tasks.
One is manipulation of algebraic inequalities;
 the other is numerical solutions.
 Often the first is far more daunting
 since the numerics is handled by some standard
 package.
 Thus there is a great need for algebraic theory.

\subsection{Engineers need generality}
\label{sec:symbcoef}
To make exposition palatable in this paper we have
refrained from generality which does not have much effect on
mathematical structure.
However, to embrace linear systems
problems we need more general theorems.
A level of generality
which most linear systems problems require is to work with
polynomials $p$ in two classes of variables
  $p(a,x)$
  where we shall be interested in
  matrix convexity in $x$ over ranges of the
  variable $a$.
  Describing this setup fully takes a while,
as one can see in
  \cite{CHSY03} where it is worked out.
An engineer might look at
\cite{CHSprept}, especially the first part which describe a
computational noncommutative algebra attack on convexity, it seems to
be the most intuitive read on the subject at hand.
Here we try to indicate the idea.
In private notes of Helton and Adrian Lim one shows that second
derivatives of $p(a,x)$ in $x$ determine convexity in $x$ and that
convexity in the $x$ variable on
  some ``open set" of $a,x$
implies that $p$ has degree 2 or less in $x$.
From this we get \\

{\it  If $P(a,x)$ is a symmetric $ d \times d$ matrix with
polynomial entries $p_{ij}(a,x)$, then convexity in $x$ for all
$X$ and all $A$ satisfying some strict algebraic inequality of the
form  $g(A) > 0$,
  implies
each $p_{ij}$ has degree 2 or less}.
\\

We obtain this from
  the following argument.
We shall test $P(a,x)$ by plugging in tuples $A$ and $X$ of $n
\times n$ matrices for $a$ and $x$. First note that matrix
convexity of $P$ in $X$ through a range of $A,X$ implies that the
diagonals $p_{ii}$  must have this property. Thus they have degree
2 or less in $x$.
Consider how the Hessian $ F(t):=P_{xx}(A,t
X)[H]$ in $x$ scales with a scalar parameter $t$.
The matrix function
being convex implies its diagonals are convex.
Thus as we saw above,
 $ F_{kk}(t) $ is independent of $t$ for all $k$. Apply $F(t)$
to the vectors $v= column\;(\pm v_1, \cdots,\pm v_d)$ in
$\Rnn$ and use that $v^T F(t) v \geq 0$ for all $t$, to get
that for each $i,j$ the entries $F_{ij}$ satisfy
$$
   v_i^T F_{ii}(t) v_i +
   v_j^T F_{jj}(t) v_j \ \  \geq \ \
\pm  (v_i^T F_{ij}(t) v_j)^2.
$$
This implies by letting $t \to \infty$ that the degree of $v_i^T
F_{ij}(t) v_j$ in $t$ is 0, which implies the same for
$F_{ij}(t)$. To this point we have that all
polynomials in $P_{xx}(a, x)[h]$ are independent of $X$
  whenever matrix tuples $A$ from an open set
  $\{A: g(A) > 0 \}$ are plugged in.
This  is independent of the size $n$ of the matrices we plug in,
so all polynomials in $P_{xx}(a, x)[h]$ are 0, algebraically
speaking. Thus all polynomials in $P(a, x)[h]$ have degree 2 in
$x$ or less. The engineering conclusion from all of this is
formulated below.
\bigskip

\subsection{Conclusion}
\begin{enumerate}
\item Many
  linear systems problems which are ``dimension free"
  readily reduce   to noncommuting inequalities
  on $d \times d$ matrices of polynomials
  of the form  $P(a,x) \leq 0$.
  These do so as in the \S \ref{sec:HinftyProb}
example, or even after simplifying solving and substituting they
yield a matrix of polynomials.
\item
  If such $P(A,X)$ is $X$-convex on the set
  of $n \times n$ matrix tuples $A$ satisfying a strict polynomial
  inequality $g(A) > 0$
   and on all $X$ (regardless of dimension $n$),
  then $P(a,x)$ has degree 2 in $x$,
  as we saw in \S \ref{sec:symbcoef}.
  Alas, $P$ is surprisingly simple.
\item \label{conc:testp}
Assume  a  $d \times d$ matrix of polynomials $P(a,x) $ has degree
2 in x. There are tests (not perfect) to see where in the $a$
variable $P(X,A)$ is negative semi-definite for all $X$.
Equivalently, to see where $P$ is convex in $x$.
\item \label{conc:LMI}
   Convexity and the  degree 2 property
   imply   $P(a,x) \leq 0$
   can be expressed as
   an LMI.
Often the LMI can be constructed with coefficients which are
noncommutative polynomials (dimension free formula).
  See proof below.
\end{enumerate}

This very strong conclusion is bad news for engineers and we
emphasize that it does not preclude transformation  to convexity
or convexity for  dimension dependent problems.

\subsubsection{Tests for convexity and the making of an LMI}

\noindent {Here we shall sketch of the proof of
Conclusions (\ref{conc:testp}) and (\ref{conc:LMI})}.
\ We use methods not  described earlier
in this paper, but despite that restrict  our presentation to be
only  a brief
outline.
For proofs in detail see  \cite{CHSY03} or more generally
\cite{HMPprepta}.

Suppose
$q(a)[h]$ is a symmetric polynomial in $a,h$ which is homogeneous
of  degree 2 in $h$, then $q$  being quadratic in $h$, can be
represented as
\begin{equation}
\label{eq:repM} q(a)[h] = V(a)[h]^T M(a) V(a)[h]
\end{equation}
where $M$ is  matrix of noncommutative polynomials in $a$,
and $V$ is
a vector each entry of which is a monomial of the form $h_j m(a)$
where $m(a)$ is a monomial in $a$.
We can choose the representation
so that no monomial repeats.
A key is Theorem 10.10 and Lemma 9.4
in \cite{CHSY03} which imply
\begin{lem}
Let $q(a,h), g(a)$ be polynomials in the free algebra
with $q$ purely quadratic in $h$.
Then $q(A)[H] \geq 0$ for  $g(A) > 0$ and all $H$ is equivalent to
$M(A) \geq 0$ for $g(A) > 0$.
\end{lem}
%
We shall apply this by representing
$ P_{xx}(a,x)[h] = V(a)[h]^T M(a) V(a)[h]$. The $x$- Hessian being
quadratic in  $x$ satisfies $P_{xx}(a, x)[h]$ is independent of
$x$.
From the lemma we have
$P_{xx}(A,X)[H] \geq 0$ for  $g(A)>  0$
and all $X$ is equivalent to $M(A) \geq 0$ for $g(A)>  0$.

\bigskip

{\it Two tests for positivity} as mentioned in Conclusion
(\ref{conc:testp}) follow.
\begin{enumerate}
\item
  The test in \cite{CHSY03} is: \ \
take the symbolic noncommutative $M(a)= L(a)^TD(a) L(a)$
decomposition of $M(a)$. This gives $D(a)$ a matrix with diagonal
or $2 \times 2$ block diagonal entries which are nc rational
functions in $a$. $M(A) \geq 0 $ if and only if $D(A) \geq 0$, so
checking positivity of the functions on diagonal $D(a)$ is a test
for where $p(a,x)$ is convex. \item
  Here is another test.
If the Positivstellensatz holds (despite a failure of the strict
positivity hypothesis), then
\begin{equation}
\label{eq:posssM} M(a,x)\in \Sigma^2 + \sum r_j^T g r_j + \sum
t_{ij}^T ( C^2 - x_j^T x_j ) t_{ij}.
\end{equation}
Computing the terms $r_j, t_{ij}$ and the sums of squares
component gives an algebraic test.
\end{enumerate}

\medskip

{\it Conversion to LMIs}, namely, Conclusion
(\ref{conc:LMI}).
Denote by  $P^{I}(a,x)$
the terms in $P(a,x)$  with $x$ degree exactly  one,
respectively  $P^{II}(a,x)$with $x$ degree exactly  two.
\begin{enumerate}
\item
  We now give
  quite a practical numerical algorithm for producing an LMI,
  under no offensive assumptions.
However, we do not get formulas which are polynomials in the
symbol $a$. Once $n \times n$ matrices $A$ are given with $M(A)$
positive semi-definite we can compute numerically its Cholesky
decomposition $M(A)= L(A)^T D(A)^{\frac 1 {2}}  D(A)^{\frac 1 {2}}
L(A)$ (actually any square root of $M(A)$ will do). Then we have
$$
P^{II}(A,X)= V(A)[X]^T L(A)^T D(A)^{\frac 1 {2}} \    D(A)^{\frac
1 {2}}L(A) V(A)[X]
$$
and taking
$$ \cL(A)[X]:
= \left(%
\begin{array}{cc}
  p^0(A,X) +  p^I(A,X) & \ \ D(A)^{\frac 1 {2}}L(A) V(A)[X] \\
   V(A)[X]^T L(A)^TD(A)^{\frac 1 {2}} & \ \ - I \\
\end{array}%
\right)
$$
which produces $\cL(A)[X]$ whose Schur
complement equals $P(A,X)$
and which produces a matrix inequality
\begin{equation}
\label{eq:equivIneq}
 \{ X: P(A,X) \leq 0 \} \ = \  \{ X: \cL(A)[X] \leq 0 \}.
\end{equation}
 The  entries of $\cL(A)[X]$  are linear in
scalar unknowns $X_{lm}$ and have $n \times n$ matrix
coefficients. This is standard input to the LMI numerical solvers
prevalent today.
\item
Another recipe which produces algebraic
formulas for solution the following.
  Continue with item (2) above.
The   terms  $P^{II}(a,x)$ in $P(a,x)$  with $x$ degree exactly
  two   can be  represented by
  $M(a)$ as in (\ref{eq:repM}).
 From the Positivstellens\"atz (\ref{eq:posssM}) for $M(a)$
  and the fact that linear terms are trivial
to handle, we can easily build an algebraic expression for a
matrix $\cL(a)[x]$  with polynomial entries which are linear in
$x$ whose Schur complement equals $P(a,x)$. Moreover, for any
fixed $A$ satisfying $g(A) > 0$, the solution sets to our favorite
matrix inequality $P>0$ and the  LMI based on $\cL$ are the same,
as in
(\ref{eq:equivIneq}).
This completes the proof that,
if the Positivstellens\"atz (\ref{eq:posssM}) for
$P^{II}(a,x)$ exists, then a LMI which is polynomial in $a$
exists.
\end{enumerate}
As a side remark, for the degree 2 and other properties of matrix
valued polynomials we could use weaker  hypotheses allowing
  coupling of  $a$ and $x$
  (as in done in private Helton- Lim notes for polynomials), these
probably work by the same  argument, basically the argument  in
\cite{HM04conv}).

\subsection{Keep going}
This subject of noncommutative real algebraic geometry and its
geometric offshoots on convexity is a child of the $21^{st}$
century. Understanding the relationship between Convex MIs and
LMIs was a core motivation  for its developments. When we look at
the two basic techniques in \S (\ref{sec:convQ}) what we have done
in this paper bears successfully on issue \ref{convQ:conv}. But
nothing has been said about issue \ref{convQ:chgV}. Nick
Slinglend's UCSD thesis in progress makes a start in that
direction.

This physical section has focused on ``dimension free" problems.
What about dimension dependent ones? In these problems the
variables commute.
There the behavior is quite different;
as we saw in \S \ref{sec:LMIs and SoS}
there is an extra constraint
beyond convexity to have equivalence to  an LMI.

\section{A guide to literature}

While classical semi-algebraic geometry has developed over the
last century through an outpouring of seemingly countless papers,
the thrust toward a noncommutative semi-algebraic geometry is
sufficiently new that we have attempted to reference the majority
of papers directly on the subject here in this survey.  This
non-discriminating approach
 is not entirely good news for the student,
so in this section we provide some guidance to the
more readable references.

The Functional Analysis book by Riesz and Nagy \cite{RN} is a
class in itself. For a historical perspective on the evolution of
the spectral theorem the reader can go directly to Hilbert's book
\cite{Hilbert4} or the German Encyclopedia article by Hellinger
and Toeplitz \cite{HT}. Reading von Neumann in original \cite{vN1}
is still very rewarding.

The many facets of matrix positivity, as applied to function
theory and systems theory, are well exposed in the books by
Agler-McCarthy \cite{AM}, Foias-Frazho \cite{FF} and
Rosenblum-Rovnyak \cite{RR}. The monograph of Constantinescu
\cite{C} is entirely devoted to the Schur algorithm.

For the classical moment problem Akhiezer's text \cite{A} remains
the basic reference, although having a look at Marcel Riesz
original articles \cite{MRiesz}, Carleman's quasi-analytic
functions \cite{Carleman}, or at the continued fractions monograph
of Perron \cite{Perron} might bring new insights. Good surveys of
the multivariate moment problems are Berg \cite{Berg1} and Fuglede
\cite{Fuglede}. Reznick's memoir \cite{Reznick1} exploits in a
novel and optimal way the duality between  moments and positive
polynomials.

For real algebraic geometry, including the logical aspects of the
theory, we refer to the well circulated texts \cite{BCR, Jacobson,
Marshall} and the recent monograph by Prestel and Delzell
\cite{PD}; the latter offers an elegant and full access to a wide
selection of aspects of positive polynomials.
 For new results in algorithmic real (commutative) algebra see \cite{BPR};
 all recent articles of Lasserre
contain generous recapitulations and
reviews of past articles
devoted to applications of sums of squares and moments to
optimization. Scheiderer's very informative survey
\cite{Scheiderer2} is centered on sums of squares decompositions.
Parrilo's thesis \cite{parThesis} is a wonderful exposition of
many new areas of application which he discovered.

An account of one of the most systematic and elegant ways for
producing LMIs for engineering problems is the subject of the book
\cite{SIG97}. The condensed version we heartily recommend is their
15 page paper \cite{SI95}.\\

\noindent {\bf Software:}

Common semi-definite programming packages are
 \cite{Sturm99}{\it  SeDuMi} and
 {\it LMI Toolbox} \cite{GNLC95}.

Semi-algebraic geometry  packages are {\it SOS tools}
\cite{PPSP04} and {\it GloptiPoly}
 \cite{Henrion1}.

For symbolic computation in a free $*$- algebra see {\it
NCAlgebra} and {\it NCGB} (which requires Mathematica) \cite{HSM05}.

{\small

\printindex

\end{document}